\documentclass[12pt]{article}
\usepackage{amsfonts,amsmath,amsxtra}
\usepackage{amssymb}
\usepackage{bm}
 \usepackage[utf8]{inputenc}
\usepackage{xcolor}

\def\Id{\operatorname{Id}}
\def\a{\alpha}
\newcommand\An[1]{A^{(#1)}}
\def\C{\mathbb{C}}
\def\b{\beta}

\def\blambda{\bm{\lambda}}
\def\beq{\begin{equation}}
\def\eeq{\end{equation}}
\def\beqq{\begin{equation*}}
\def\eeqq{\end{equation*}}
\def\d{d}

\def\eps{\epsilon}

\def\g{\mathfrak{g}}
\def\gl{\mathfrak{gl}}
\def\GL{\mathbf{GL}}
\def\SO{\mathbf{SO}}
\def\SL{\mathbf{SL}}
\def\Sp{\mathbf{Sp}}
\def\h{\mathfrak{h}}
\def\i{\hat{i}}
\def\j{\hat{j}}
\def\k{\bar{k}}

\def\Mat{\operatorname{Mat}}
\def\n{\mathfrak{n}}
\def\N{N_{*}}
\def\p{\mathbf{p}}
\def\Psiold{\tilde{\Psi}}
\def\Psinew{\Psi}
\def\R{\mathbb{R}}
\def\Re{\mathrm{Re}\,}
\def\H{\mathcal {H}}

\newcommand\Ln[1]{\Lambda^{(#1)}}

\def\so{\mathfrak{so}}
\def\sp{\mathfrak{sp}}
\def\ss{\bar{s}}
\def\t{\mathbf{t}}

\def\tgamma{\tilde{\gamma}}
\def\tdelta{\tilde{\delta}}

\def\tr{\tilde{r}}

\def\ts{\bar{s}}

\newcommand\tXn[1]{\tilde{X}^{(#1)}}
\newcommand\ttXn[1]{\bar{X}^{(#1)}}

\def\U{\mathrm{U}}
\def\V{{V}}
\def\ve{\varepsilon}
\def\vf{\varphi}
\newcommand\Xn[1]{X^{(#1)}}
\def\x{\bm{x}}
\def\Z{\mathbb{Z}}
\def\ZZ{\mathrm{Z}}
\def\z{{\bm z}}
\def\wr{v^R}
\def\wl{v^L}
\def\wwr{\mathrm{\omega}^R}
\def\wwl{\mathrm{\omega}^L}
\def\w{\bar{w}}
\def\Ad{{\rm Ad}}

 \def\ga{\gamma}

   \def\wwrm{\hat{{\omega}}^R}
   
   \def\wwlm{\hat{{\omega}}^L}
    \def\wwwlm{\bar{\hat{{\omega}}}^L}
\newtheorem{lemma}{Lemma}[section]
\newtheorem{proposition}{Proposition}[section]
\newtheorem{corollary}{Corollary}[section]
\newtheorem{theorem}{Theorem}[section]

\newcommand{\rf}[1]{(\ref{#1})}

\begin{document}
	\begin{center}
		{\bf \large Mellin-Barnes presentations for Whittaker wave functions}
			\bigskip
			
	{\bf S. Kharchev$^{\,\star\,\natural}$,\, S. Khoroshkin$^{\,\star\,\circ}$,
	}\medskip\\
	$^\star${\it Institute for Theoretical and Experimental Physics, Moscow, Russia;}\smallskip\\
	$^\natural${\it
	Institute for Information Transmission Problems RAS (Kharkevich Institute),Bolshoy Karetny per. 19, Moscow, 127994, Russia;}\smallskip\\
$^\circ${\it National Research University Higher School of Economics, Moscow, Russia.}
\end{center}
\begin{abstract}\noindent  We obtain certain Mellin-Barnes integrals which present Whittaker wave functions related to classical real split forms of simple complex Lie groups.
	\end{abstract}		
\bigskip
\begin{center} {\it Keywords}: Whittaker function, Mellin transform, Lusztig parametrization
	\end{center}
\section{Introduction} {\bf 1}. In this paper we obtain certain Mellin-Barnes integrals which present Whittaker functions related to classical real split forms of simple complex Lie groups.

 Whittaker functions originally appeared as solutions of a special differential equation of hypergeometric type \cite{WW}. They were then realized as eigenfunctions of Laplace operators of the group
 $\GL(2,\R)$, see e.g. \cite{V}. This construction was generalized to real semisimple groups and led to a series of significant researches in representation theory, see e.g. \cite{J,Sch}. B.Kostant related the group theory of Whittaker functions with a family of Toda integrable systems \cite{T}.

 Whittaker functions admit several integral presentations.  Particular examples were obtained more than forty years before, see e.g. \cite{Bu}.
 S.Kharchev and D.Lebedev in \cite{KL} found integral presentation of $\GL(n,\R)$ Whittaker wave functions using the machinery of inverse scattering method. Then in \cite{GKL} the same presentation was obtained by calculating matrix elements in certain infinite dimensional "Gelfand-Zetlin" representations. A.Givental found quite different integral presentation for the same Whittaker function using geometric arguments. A.Gerasimov  et al. then realized in \cite{GKLO,GLO}
 that Givental construction can be reformulated as the description of the matrix element in principal series representation using Gauss decomposition and Lusztig coordinates \cite{L} on nilpotent subgroup.
  Moreover, Whittaker wave functions for all classical split real groups were described in \cite{GLO} as certain integrals over positive cone in corresponding maximal nilpotent subgroup.
 Recently both presentations were generalized to a quantum group setting \cite{SS} using the machinery of cluster mutations. In particular, the completeness and orthogonality of $q$-versions of $\GL(n,\R)$ Whittaker wave functions were proved there.

 We start with the presentation of Whittaker wave function as of special matrix element  \cite{GLO}
 \beq\label{00}	\Psinew_\lambda(x) =e^{-(\rho,x)}(\wl_{i\lambda-\rho},\exp{(-x)}\wr_{i\lambda-\rho}),
 \eeq
  see \rf{5} and \rf{5a} for precise notations, and rewrite this matrix element as Barnes integral.  Our technique is rather elementary. It contains  four ingredients: use of Lusztig coordinates, Berenstein-Zelevinsky transform, Plancherel formula for Mellin transform and linear algebra matrix calculations. Lusztig coordinates are convenient in the description of Whittaker functions by several reasons. They separate Lusztig positive cone $N_*$, which is original space of integration of the matrix element in consideration. Cluster type mutations between different Lusztig chats enable us to avoid the use of formulas for the action of the Lie algebra and observe elegant expressions for invariant forms and vectors.

 The definition \rf{00} of the matrix element exploit two special vectors in generalized principal series;
  they are called left and right Whittaker vector and the Whittaker functions is the matrix element of the Cartan flow related to these vectors.
  The right Whittaker vector has a simple direct expression in Lusztig coordinates, see Proposition \ref{prop2}. The definition of the left Whittaker vector requires the conjugation by the longest element $w_0$ of the Weyl group and the calculation of Gauss coordinates of the new matrix. This induces a birational map of the nilpotent subgroup $N$ which in slightly different setting was studied by A.Berenstein and A.Zelevinsky \cite{BZ} and was crucial for their description of Lusztig coordinates via 'generalised minors' -- matrix elements in fundamental representations. We find another expressions for this birational map (we call it BZ transform) just by elementary linear algebra calculations using the induction by the rank of the group.  This is the main point of the construction and we suspect that the formulas for BZ maps which we found will serve elsewhere. The precise description of BZ maps is given in Theorems 3.1, 4.1, 5.1 and 6.1

  The rest is the application of the convolution theorem for Mellin transform, which we use in a form of Plancherel formula
   \cite{Tt}. It gives a presentation of Whittaker functions for classical split real groups $\GL(n,\R)$, $\SO(n,n)$, $\SO(n,n+1)$ and
   $\Sp(2n,\R)$ (which present root system of $A,D,B$ and $C$ series) by means of Barnes integrals.
   These presentations can be equivalently rewritten as Mellin transforms of Whittaker functions. See Section \ref{section6}.

 For the group $\GL(n,\R)$ we have now two different integral presentation of Whittaker function by Mellin--Barnes integral:
 Kharchev-Lebedev formulas \cite{KL,GKL} and  \rf{I1} of the present paper.
 They are quite different. We do not know how to derive one from another. On the other hand, Mellin transform for $\GL(n,\R)$ Whittaker function was calculated by E.Stade \cite{St}. Our formula \rf{M2} also differs from that of \cite{St}.
 
  The following are the main results of the paper.
 \bigskip

 {\bf 2}. Let $\x=(x_1,\ldots,x_n)$, $x_i\in\R$ be the coordinates on the Cartan subalgebra $\h$ of $\gl(n,\R)$,
 $\lambda_i\in\R$ be dual coordinates on $\h^*$, $\blambda=(\lambda_1,\ldots,\lambda_n)$,
 so that the element $h(\x)\in\h$ is given
 by the diagonal matrix $h=\sum_{i=1}^n x_ie_{ii}$; and $(\x,\blambda)=\sum_{k=1}^nx_k\lambda_k$
 . Let $\Psinew_{\blambda}(\x)$ be  the Whittaker wave function for $\GL(n,\R)$, defined by the relation \rf{5a}.
 \smallskip

 {\bf Theorem \ref{t1}} {\em The function $\Psinew_{\blambda}(\x)$ is given by the integral}
 \beq\label{I1}
 \begin{split}
	\Psinew_{\blambda}(\x)&=\frac{e^{  -i(\x,\blambda)}}{(2\pi i)^\d} 
	\int\limits_{C}\exp\sum_{1\leq k\leq n
	}({\gamma}_{1,n+2-k}-\gamma_{1,n+1-k})x_{k}
	\\&
	\prod_{1\leq k<l\leq n}{\Gamma}(\gamma_{kl}-\gamma_{k+1,l}+i(\lambda_k-\lambda_{n+k-l+1}))
	\Gamma(\gamma_{kl}-\gamma_{k+1,l+1})d\gamma_{k,l}.
\end{split}\eeq

 Here  we set $\gamma_{k,l}=0$ if the pair $(k,l)$ does not satisfies the condition $1\leq k<l\leq n$;
 $\d=n(n-1)/2$ is the dimension of the maximal unipotent subgroup of $\GL(n,\R)$. The integration cycle $C$ is a deformation of the imaginary plain
 $\Re \gamma_{k,l}=0$ into the domain  $D\subset\C^{\d}$ of the analyticity of the integrand. For instance one can use iterated integration, which
 starts with
  integration over $\gamma_{n-1,n}$ the over $\gamma_{n-2,n-1}$, then over $\gamma_{n-2,n}$ etc., which respects the conditions $\Re \gamma_{k,l}>\Re \gamma_{k+1,m}$
   for all admissible triples $k,l,m$.
   \bigskip

   Let $\x=(x_1,\ldots,x_n)$, $x_i\in\R$ be the coordinates on the Cartan subalgebra $\h$ of $\so(n,n)$,  and
 $\blambda=(\lambda_1,\ldots,\lambda_n)$, $\lambda_i\in\R$ be dual coordinates on $\h^*$, so that the element $h(\x)\in\h$ is given
 by the diagonal matrix $h=\sum_{i=1}^n x_i(e_{ii}-e_{2n+1-i,2n+1-i})$ and $(\x,\blambda)=\sum_{k=1}^nx_k\lambda_k$. Let $\Psinew_{\blambda}(\x)$
 be the Whittaker wave function for $\SO(n,n)$  defined by the relation \rf{5a}.\smallskip

 {\bf Theorem \ref{t2}} {\em The function $\Psinew_{\blambda}(\x)$ is given by the integral}
 \beqq
 \begin{split}
	\Psinew_{\blambda}(\x)&= \frac{e^{  -i(\x,\blambda)}}{(2\pi i)^\d}\int\limits_{C}	
	\exp {H}(\bm{x},\bm{\gamma})
	\prod_{k=1}^{n-2}\frac{\Gamma(\gamma_{k,k+1}+\delta_{k,k+1}+2i\lambda_k)}{\Gamma(\gamma_{k,n}
		\!+\delta_{k,n}\!-\gamma_{k+1,n}\!-\delta_{k+1,n}\!+
		2i\lambda_k)} \\
	 &\!\!\!\prod_{1\leq k<l\leq n} \!\!\!\Gamma({\xi}_{k,l}+i(\lambda_k-\!\lambda_l))\Gamma({\eta}_{k,l}+i(\lambda_k+\!\lambda_l))
	\Gamma(\gamma_{k,l}-\gamma_{k+1,l})\Gamma(\delta_{k,l}-\delta_{k+1,l})d\gamma_{k,l}d\delta_{k,l}.
\end{split}\eeqq

 Here we set $\gamma_{k,l}=\delta_{k,l}=0$ if the pair $(k,l)$ does not satisfies the condition $1\leq k<l\leq n$,
  $\d=n(n-1)$ is the dimension of the maximal unipotent subgroup of $\SO(n,n)$.
 \beq\label{I3}
	H(\bm{x},\bm{\gamma})=\!\sum_{j<n}(\gamma_{1,j}+\delta_{1,j})(x_j- x_{j-1})+\gamma_{1,n}  (x_n-x_{n-1}) -\delta_{1,n}(x_{n-1}+x_n),\eeq
\beq
	\label{I4}
	\begin{split}
		&\xi_{i,j}=-\gamma_{i+1,j+1}-\delta_{i+1,j+1}+\gamma_{i+1,j}+\delta_{i,j},\qquad 1\leq i<j<n,\\	
		&\eta_{i,j}=\gamma_{i,j+1}+\delta_{i,j+1}-\gamma_{i+1,j}-\delta_{i,j},\qquad 1\leq i<j<n,\\
		&\xi_{i,n}=\gamma_{i,n}-\delta_{i+1,n},\qquad
		\eta_{i,n}=\delta_{i,n}-\gamma_{i+1,n}\qquad 1\leq i<n,	\end{split}
	\eeq
	 The integration cycle is a deformation of the imaginary plain
 $\Re \gamma_{k,l}=\Re \delta_{k,l}=0$ into nonempty domain  $D\subset\C^{\d}$ of the analyticity of the integrand, which is described by the relations
 \beqq
 	\Re \gamma_{i,j}>0, \ \Re \delta_{i,j}>0,\ Re\, \xi_{i,j}>0,\ \Re \eta_{i,j}>0,
 	\qquad 1\leq i<j\leq n.
 \eeqq
\smallskip

 Let $\x=(x_1,\ldots,x_n)$, $x_i\in\R$ be the coordinates on the Cartan subalgebra $\h$ of $\so(n+1,n)$, and
 $\blambda=(\lambda_1,\ldots,\lambda_n)$, $\lambda_i\in\R$ be dual coordinates on $\h^*$, so that the element $h(\x)\in\h$ is given
 by the diagonal matrix $h=\sum_{i=1}^n x_i(e_{ii}-e_{2n+2-i,2n+2-i})$ and  $(\x,\blambda)=\sum_{k=1}^nx_k\lambda_k$. Let $\Psinew_{\blambda}(\x)$ be
  the Whittaker function for $\SO(n+1,n)$ defined by the relation \rf{5a}.\smallskip

 {\bf Theorem \ref{t3}} {\em The function $\Psinew_{\blambda}(\x)$ is given by the integral}
 \begin{align*}
  &\Psinew_{\blambda}(\x)= \frac{e^{  -i(\x,\blambda)}}{(2\pi i)^\d}\int\limits_{C}
	\exp {H}(\bm{x},{\bm{\gamma}})\Gamma(\gamma_n+2i\lambda_n)\Gamma(\gamma_n)d\gamma_n\cdot \\
&
\ \, \prod_{k=1}^{n-1}{\Gamma(\gamma_{k,k+1}+\delta_{k,k+1}+2i\lambda_k)}
\Gamma(\gamma_k-\gamma_{k+1})d\gamma_k\cdot \\
&\!\!\prod_{1\leq k<l\leq n}\!\! \Gamma({\xi}_{k,l}+i(\lambda_k-\lambda_l))\Gamma({\eta}_{k,l}+i(\lambda_k+\lambda_l))
\Gamma(\gamma_{k,l}-\gamma_{k+1,l})\Gamma(\delta_{k,l}-\delta_{k+1,l})d\gamma_{k,l}d\delta_{k,l}.
\notag\end{align*}	
Here we set $\gamma_{k,l}=\delta_{k,l}=0$ if the pair $(k,l)$ does not satisfies the condition $1\leq k<l\leq n$, and $\gamma_k=0$ if $k=1$ or $k=n+1$, $\d=n^2$ is the dimension of the maximal unipotent subgroup of $\SO(n+1,n)$
\beq\label{I5}
	H(\bm{x},{\bm{\gamma}})=\sum_{j=2}^n(\gamma_{1,j}+\delta_{1,j})(x_{j}-x_{j-1}) -\gamma_1 x_n.\eeq
\beq
\label{I6}
\begin{split}
&\xi_{i,j}=-\gamma_{i+1,j+1}-\delta_{i+1,j+1}+\gamma_{i+1,j}+\delta_{i,j},\qquad 1\leq i<j<n,\\	
&\eta_{i,j}=\gamma_{i,j+1}+\delta_{i,j+1}-\gamma_{i+1,j}-\delta_{i,j},\qquad 1\leq i<j<n,\\
&\xi_{i,n}=\gamma_{i+1,n}+\delta_{i,n}-\gamma_{i+1},\qquad
\eta_{i,n}=-\gamma_{i+1,n}-\delta_{i,n}+\gamma_i, \qquad 1\leq i<n.	\end{split}
\eeq
	 The integration cycle is a deformation of the imaginary plain
$\Re \gamma_{k,l}=\Re \delta_{k,l}=\Re \gamma_k=0$ into nonempty domain  $D\subset\C^{\d}$ of the analyticity of the integrand, which is described by the relations
\beqq
\Re \gamma_{i,j}>0, \ \Re \delta_{i,j}>0,\ \ Re\, \xi_{i,j}>0,\ \Re \eta_{i,j}>0,
\ 1\leq i<j\leq n, \ \ \Re \gamma_k>0,\ k=1,...,n.
\eeqq
\bigskip

 Let $\x=(x_1,\ldots,x_n)$, $x_i\in\R$ be the coordinates on the Cartan subalgebra $\h$ of $\sp(2n,\R)$,  and
$\blambda=(\lambda_1,\ldots,\lambda_n)$, $\lambda_i\in\R$ be dual coordinates on $\h^*$, so that the element $h(\x)\in\h$ is given
by the diagonal matrix $h=\sum_{i=1}^n x_i(e_{ii}-e_{2n+1-i,2n+1-i})$ and  $(\x,\blambda)=\sum_{k=1}^nx_k\lambda_k$. Let $\Psinew_{\blambda}(\x)$
be the Whittaker function for $\Sp(2n,\R)$ defined by the relation \rf{5a}.
\smallskip

{\bf Theorem \ref{t4}}. {\em The function $\Psinew_{\blambda}(\x)$ is given by the integral}
\begin{align*}
&\Psinew_{\blambda}(\x)= \frac{e^{  -i(\x,\blambda)}}{(2\pi i)^\d}\int\limits_{C} \!\!
\exp {H}(\bm{x},\bm{\gamma})
\prod_{k=1}^{n-1}\frac{\Gamma(\gamma_{k,k+1}+\delta_{k,k+1}+2i\lambda_k)}
{\Gamma(2\gamma_k-2\gamma_{k+1}+2\lambda_k)}\cdot
\\
&
\prod_{k<l} \Gamma({\xi}_{k,l}+i(\lambda_k-\lambda_l))\Gamma({\eta}_{k,l}+
i(\lambda_k+\lambda_l))
\Gamma(\gamma_{k,l}-\gamma_{k+1,l})\Gamma(\delta_{k,l}-\delta_{k+1,l})d\gamma_{k,l}
d\delta_{k,l}\cdot \\
& \prod_{k=1}^n{\Gamma(\gamma_k-\gamma_{k+1}+i\lambda_k)}\Gamma(\gamma_k-\gamma_{k+1})d\gamma_k.
\end{align*}	

Here we set $\gamma_{k,l}=\delta_{k,l}=0$ if the pair $(k,l)$ does not satisfies the condition $1\leq k<l\leq n$, and $\gamma_k=0$ if $k=1$ or $k=n+1$, $\d=n^2$,
\beq\label{I7}
H(\bm{x},{\bm{\gamma}})=\sum_{j=2}^n(\gamma_{1,j}+\delta_{1,j})(x_{j}-x_{j-1}) -2\gamma_1 x_n.\eeq
\beq
\label{I8}
\begin{split}
	&\xi_{i,j}=-\gamma_{i+1,j+1}-\delta_{i+1,j+1}+\gamma_{i+1,j}+\delta_{i,j},\qquad 1\leq i<j<n,\\	
	&\eta_{i,j}=\gamma_{i,j+1}+\delta_{i,j+1}-\gamma_{i+1,j}-\delta_{i,j},\qquad 1\leq i<j<n,\\
	&\xi_{i,n}=\gamma_{i+1,n}+\delta_{i,n}-2\gamma_{i+1},\qquad 1\leq i<n,\\
	&\eta_{i,n}=-\gamma_{i+1,n}-\delta_{i,n}+2\gamma_i,	\end{split}
\eeq
 The integration cycle is a deformation of the imaginary plain
$\Re \gamma_{k,l}=\Re \delta_{k,l}=\Re \gamma_k=0$ into nonempty domain  $D\subset\C^{\d}$ of the analyticity of the integrand, which is described by the relations
\beqq
\Re \gamma_{i,j}>0, \ \Re \delta_{i,j}>0,\ \ Re\, \xi_{i,j}>0,\ \Re \eta_{i,j}>0,
\ 1\leq i<j\leq n, \ \ \Re \gamma_k>0,\ k=1,...,n.
\eeqq

\section{Generalities}
\subsection{Whittaker vectors and Whittaker wave functions}

Let $G$ be a split real form of a reductive group over $\C$, $B_\pm$ two opposite Borel subgroups of $G$,
$N_\pm$ their maximal nilpotent subgroups and  $H=B_+\cap B_-$ the Cartan subgroup. Let $D(G)$ be the ring of bi-invariant differential operators on $G$. Denote by ${\mathfrak b}_\pm$,
$\n_\pm$ and $\h$ the corresponding Lie algebras. Let $\Delta_\pm\in\h^*$ be the systems of positive and negative roots of
$\g$, $\Pi\subset \Delta_+$ be a subsystem of simple roots and $W$ the Weyl group of $G$. Denote by $G_0$ the big Bruhat cell $G_0=N_-HN_+$. It is  dense open in $G$.

For each index $i$ of a simple root $\alpha_i\in\Pi$ denote by $e_i=e_{\alpha_i}$, $f_i=f_{\alpha_i}=e_{-\alpha_i}$ and $h_i=h_{\alpha_i}$ the corresponding Chevalley generators of $\g$, so that
\begin{equation}\label{L2a}\alpha_i(h_i)=2\end{equation}
and $e_i,f_i$ and $h_i=\alpha_i^\vee$ are standard generators of the embedded Lie algebra $\mathfrak{sl}_2$:
$$[h_i,e_i]=2e_i,\qquad [h_i,f_i]=-2f_i,\qquad [e_i,f_i]=h_i.$$

 Let $\zeta_\pm : N_\pm \to \C^* $ be  nondegenerate characters, defined by the relations \begin{equation}\label{0}\zeta_+(\exp(te_{i})):=e^{-t},\qquad \zeta_-(\exp(tf_{i})):=e^{-t}
 \end{equation} for all simple roots $\alpha_i$. In this paper   Whittaker function  is a analytical function $\vf$ on  $ G_0$ satisfying the conditions
 \begin{align}\label{1}
 &\vf(g) \qquad \text{is an eigenfunction for any}\qquad D\in  D(G);\\ \label{1a}
  &\vf( n_-gn_+)=\zeta_-(n_-)\vf(g)\zeta_+(n_+).\end{align}
  for any $g\in G_0$, $n_-\in N_-$, $n_+\in N_+$. The condition \rf{1} implies that Whittaker functions are completely determined by their restriction to Cartan subgroup $H$. B.~Kostant noticed \cite{K} that the restriction of
  the action of the center ${\mathrm Z}(\g)$ of universal enveloping algebra $\U(\g)$ to the space of Whittaker functions can be identified with Hamiltonians of the Toda chain related to the root system $\Delta=\Delta_+\sqcup\Delta_-$.

 Whittaker functions can be constructed as matrix elements between a pair of dual Whittaker vectors. Let $V$ be a representation
 of Lie algebra $\g$, such that its restriction to ${\mathfrak b}_+$ admits an extension to representation of the Borel group
 $B_+$ compatible with $\g$-module structure on $V$, that is $bgb^{-1}u=\Ad_b(g)u$ for any $u\in V$, $g\in\g$ and
 $b\in B_+$\footnote{$(\g,B)$- module in Harish-Chandra terminology}. A vector $v\in V$ is called right Whittaker vector
 if
 $$e_{i}v=-v\qquad \text{for any}\qquad \alpha_i\in\Pi.$$ Let also $V'$ be a representation of Lie algebra $\g$, such that its restriction to
 ${\mathfrak b}_-$ admits an extension to representation of the Borel group $B_-$ compatible with $\g$-module structure on $V'$.
 A vector $v'\in V'$ is called left Whittaker vector if
 $$f_{i}v'=-v'\qquad \text{for any}\qquad \alpha_i\in\Pi.$$

  Assume now that $V$ and $V'$ are dual
 to each other, that is there is a nondegenerate pairing $(\,,\, ): V'\otimes V\to\C$
 such that $( xu',u)+( u',xu)=0$ for any $u'\in V'$, $u\in V$ and $x\in\g$. Then the matrix element
  \begin{equation}\label{2}\vf(g)= (v',gv)\end{equation}
  is well defined for any $g\in G_0$ satisfies the condition \rf{1a} .

  If in addition the representations $V$ and $V'$ are quiasi-simple, that is the center $\ZZ(\g)$ acts on them by scalar operators, then the Whittaker functions \rf{2} are eigenfunctions of generalized Toda Hamiltonians. It is common to use for this aim the representations of $G$ induced from Borel subalgebra $B_-$. Namely for any $\mu\in\h^*$ let $\chi_\mu:B_-\to\C^*$ be a one-dimensional representation (character) of $B_-$  defined by the relation
  \begin{equation*} \chi(e^{h}n)=e^{\mu(h)}\qquad\text{for any}\ h\in\h\ \text{and}\ n\in N_-.	
  	\end{equation*}
  Denote by $\V_\mu$ the space of analytical functions on $G_0$ satisfying the condition
  \begin{equation}\label{4} f(bg)=\chi_\mu(b)f(g)\qquad \text{for any}\qquad b\in\, B_-, g\in G_0
  	\end{equation}
The element of the group $B_+$ act on the functions from $\V_\mu$ by the right shifts, $b\cdot f(g)=f(gb)$ for $g\in G_0$ and $b\in B_+$ and the elements of $\g$ act by infinitesimal right shifts,
\beq \label{4a} x f(g)=\frac{d}{d t}f(ge^{tx})|_{t=0}.\eeq
 This induced module is quasisimple, see \cite{Zh}, and can be regarded to wide extent as a representation of non-unitary principal series.
If $\mu+\nu=-2\rho$ then the pairing
\begin{equation*} \langle f,g\rangle =\int_{N_+}f(n)g(n)dn,\qquad f\in\V_\mu, g\in\V_{\nu}
	\end{equation*}
is formally invariant.  Here $\rho=\frac{1}{2}\sum_{\alpha\in\Delta_+}{ \alpha}$. The pairing  is invariant under condition of the convergency of the integral. Here $dn$ is invariant measure on $N_+$. We use instead sesquilinear pairing
\begin{equation}\label{4b} (f,g)=\int_{\N}\bar{f}(n){g}(n)dn,\qquad f\in\V_\mu, g\in\V_{\nu}
\end{equation}
where $\N\subset N_+$ is Lusztig positive cone, see \rf{NL}. It is  invariant when $\nu+\bar{\mu}=-2\rho$ and $f$ and $g$ rapidly vanish at the boundary of $\N$.

In the following we investigate right Whittaker vector $\wr_\mu$ in the space $\V_{\mu}$, left Whittaker vector
$\wl_{\nu}$ in $\V_{\nu}$, where $\nu+\bar{\mu}=-2\rho$ and the Whittaker wave function
\begin{equation}\label{5}
	\Psiold_{\mu,\nu}(x)=e^{-(\rho,x)}(\wl_{\nu},\exp{(-x)}\wr_\mu)
	\end{equation}
 Here $x\in\h$ is an element of Cartan subalgebra. When $\mu=\nu$ is a parameter of unitary principal series, $\mu=\nu=i\lambda-\rho$, where $\lambda$ is real, that is $\lambda(h_k)\in\R$ for all simple roots $\alpha_k$, see \rf{L2a}, the integral in the RHS of \rf{5} definitely converges  and has the form
 \begin{equation}\label{5a}
 	\Psinew_\lambda(x):= \Psiold_{i\lambda-\rho, i\lambda-\rho}(x)=e^{-(\rho,x)}(\wl_{i\lambda-\rho},\exp{(-x)}\wr_{i\lambda-\rho})
 	\end{equation}
 The functions $\Psiold_{\mu,\nu}(x)$ and $\Psinew_\lambda(x)$ are eigenfunctions for a family of
  Toda
 Hamiltonians $\H_k$,
 \begin{equation*}
 	\H_k \Psinew_\lambda(x)=e^{-(\rho,x)}C_ke^{(\rho,x)}\Psinew_\lambda(x)= e^{-(\rho,x)}(\wl_{i\lambda-\rho},\exp{(-x)}C_k\wr_{i\lambda-\rho}),
 	\end{equation*}
 where $C_k$ are generators of the ring $\ZZ(\g)$.
 We will also call $\Psinew_\lambda(x)$ Whittaker wave function despite it differs from the restriction to Cartan subgroup
 of the Whittaker function $(\wl_{i\lambda-\rho},g\wr_{i\lambda-\rho})$ by the normalizing factor $e^{-(\rho,x)}$ chosen for the agreement with Toda Hamiltonians.

The Whittaker function in  a form of matrix coefficient \rf{5a} and \rf{4b} was studied in \cite{GKLO, GLO}. The paper \cite{GKLO} contains integral presentations of Whittaker functions generalizing Givental formula \cite{G} for $\gl_n$. The integration over positive cone $N_*$ implies the important property of this construction: the function $\Psi_{\lambda}(x)$ rapidly decreases in the region
 \beqq
 \h_*=\{h\in\h^*, (\alpha,h)>0 \qquad \text{for all}\ \ \alpha\in\Delta_+\}.  \eeqq
 For $G=\GL(n,\R)$ this Whittaker function is known to be symmetric on parameters $\bm{\lambda}$, see \cite{Si,SS}.

 \subsection{Lusztig coordinates}
 In this subsection we describe Lusztig parametrization of the group $N_+$ in slightly different notation.

 Let $w_0$ be the longest element of the Weyl group $W$ and
 \begin{equation}\label{L1}w_0=s_{i_1}s_{i_{2}}... s_{i_{N-1}}s_{i_N}\end{equation}
 be its reduced decomposition. Here $s_{j}$ is the simple reflection in $\h^*$ corresponding to the root $\alpha_j$. We associate to \rf{L1} the following normal (or convex in other terminology) ordering  of the set $\Delta_+$:
 \begin{equation}\label{L2}
 	\gamma_1=\a_{i_{1}}, \gamma_2=s_{i_1}(\a_{i_2}),
 	\gamma_{3}= s_{i_{1}} s_{i_{2}}(\alpha_{i_{3}}),\ \ldots \,
 	\gamma_N=s_{i_1}s_{i_{2}}\!\cdots s_{i_{N-1}}(\alpha_N).
 	%
 \end{equation}
A normal orderings $<$ of the system $\Delta_+$ of positive roots of a reductive Lie algebra over $\C$ is  characterized by the condition
\beq\label{L2b}\alpha<\alpha+\beta<\beta \qquad\text{or}\qquad \beta<\alpha+\beta<\alpha\eeq
if $\alpha$, $\beta$,  $\alpha+\beta$ are all in $\Delta_+$.
The rule \rf{L2} establishes a bijective correspondence between
 reduced decompositions of $w_0$ and normal orderings of positive roots \cite{Zh2}.

   Elementary transformations of reduced decompositions are performed by means of braid group relations
 $$\underbrace{s_is_j\cdots}_{n_{i,j}+1}=\underbrace{s_js_i\cdots}_{n_{i,j}+1},\qquad i\not=j,
 $$
 where $n_{i,j}=a_{i,j}a_{j,i}+1$ and $a_{i,j}$ is the entry of  Cartan matrix, are reformulated into changes of normal orderings inside subsystems of the rank two:
 $$\begin{array}{ccl}
 \alpha,\a+\b,\beta\to \beta,\a+\b,\alpha&\text{if}& \a,\b\in{\mathrm A}_2\\
  \alpha,\a+\b,\a+2\b,\beta\to \beta,\a+2\b,\a+\b,\alpha&\text{if}& \a,\b\in{\mathrm B}_2\\
  \alpha,\a\!+\b,2\a\!+3\b,\a\!+2\b,\a\!+3\b,\beta\to \beta,\a\!+3\b,\a\!+2\b,2\a\!+3\b,\a\!+\b,\alpha&\text{if}& \a,\b\in{\mathrm G}_2
  \end{array}$$
  and
  $$\alpha,\beta\to \beta,\alpha\qquad\text{if}\qquad \a,\b\in{\mathrm A}_1\times {\mathrm A}_1.$$
   The reverse of the normal ordering of positive roots does not destroys its defining property \rf{L2b} and thus is the normal ordering as well. Moreover the reverse respects the above transformations of normal orderings and thus defines an involutive automorphism of the root system preserving the subsystem of positive roots. Thus it is induced by an automorphism $\theta$ of Dynkin diagram \footnote{ We further preserve the notation $\theta$ for the corresponding automorphism of Lie algebra $\g$ and group $G$.}. In particular this means that the last root
   $\gamma_N$ is simple,  as well as the first root $\gamma_1=\alpha_{i_1}$, and
   \begin{equation*}\gamma_N=\theta(\alpha_{i_N})\end{equation*}

  Following Lusztig \cite{L} we associate to each reduced decomposition \rf{L1} (or, equivalently, to the related normal ordering \rf{L2}) the group element $X(\t)\in N_+$
  \begin{equation}\label{L3} X(\t)=\exp\left(t_{\gamma_1}e_{i_1}\right)\cdots  \exp\left(t_{\gamma_N}e_{i_N}\right) 	
  	\end{equation}
  The correspondence \rf{L3} establishes a birational isomorphism of the varieties $N_+$ and $\R^{|\Delta_+|}$. Denote by $\N\subset N_+$, $\N\sim \R_+^{|\Delta_+|}$ the open subset of $N$ defined by the conditions $t_\gamma>0$ for all $\gamma\in\Delta_+$
  \beq\label{NL}
   \N=\{ X(\t)\ |\ t_\gamma>0, \gamma\in\Delta_+\}
  \eeq
  The definition of $\N$ does not depend on the choice of reduced decomposition of $w_0$, see \rf{L4}, \rf{L5}, \rf{L6}.
  The passage to another reduced decomposition defines the involutive transition map
  $$X(\t')=X(\t)$$
  where
  \begin{equation}\label{L4}
  	t'_\a=\frac{t_\a t_{\a+\b}}{t_a+t_\b},\qquad t'_{\a+\b}=t_\a+t_\b,\qquad t'_\b=\frac{t_\b t_{\a+\b}}{t_a+t_\b}
  	\end{equation}
  for the changes of the normal order
  $$ \ldots \a,\a+\b,\b,\ldots \to \ldots \b,\a+\b,\a, \ldots \quad \text{or}  \quad \ldots \b,\a+\b,\a,\ldots \to \ldots \a,\a+\b,\b, \ldots$$  of $\mathrm{A}_2$ subsystem, see \cite{L};
  \begin{equation}\label{L5}
  t'_\a=\frac{t_{\a+2\b} t^2_{\a+\b}t_{\a}}{\pi_2},\qquad
    	t'_{\a+\b}=\frac{\pi_2}{\pi_1},\qquad
    	 t'_{\a+2\b}=\frac{\pi_1^2}{\pi_2},\qquad
    t'_\b=\frac{t_\b t_{\a+\b}t_{\a+2\b}}{\pi_1},
    \end{equation}
    where
    $$\pi_1= t_\b t_{\a+2\b}+(t_\b+t_{\a+\b})t_\a, \qquad \pi_2=t_\b^2 t_{\a+2\b}+(t_\b+t_{\a+\b})^2t_\a$$
  for the changes of the normal order
   \begin{equation*}\begin{split}
    &\ldots, \a,\a+\b,\a+2\b, \b,\ldots \to \ldots \b,\a+2\b, \a+\b,\a, \ldots\qquad \text{or}\\  &\ldots,\b,\a+2\b, \a+\b,\a,\ldots \to \ldots \a,\a+\b,\a+2\b, \b,\ldots
    \end{split}\end{equation*} of $\mathrm{B}_2$ subsystem, see \cite[Theorem 3.1]{BZ}; and
    \begin{equation}\label{L6}\begin{split}
      & t'_\a=\frac{t_{\a} t^3_{\a+\b}t^2_{2\a+3\b}t_{\a+2\b}^3 t_{\a+3\b}}{\pi_3},\qquad
        	t'_{\a+\b}=\frac{\pi_3}{\pi_2},\qquad
        	 t'_{2\a+3\b}=\frac{\pi_2^3}{\pi_3\pi_4},\\
        & t'_{\a+2\b}=\frac{\pi_4}{\pi_1\pi_2},\qquad
        t'_{\a+3\b}=\frac{\pi_1^3}{\pi_4}\qquad
        t'_{\b}=\frac{t_\b t_{\a+3\b}t_{\a+2\b}^2t_{2\a+3\b}t_{\a+\b}}{\pi_1},\end{split}
        \end{equation}
        where
  \begin{align*}
  \pi_1=& t_{\b}t_{\a+3\b}t_{\a+2\b}^2t_{2\a+3\b}+t_{\b}t_{\a+3\b}(t_{\a+2\b}+t_{\a+\b})^2t_{\a}+(t_{\b}+t_{\a+2\b})t_{2\a+3\b}t_{\a+\b}^2t_{\a},\\
  \pi_2=& t_{\b}^2t_{\a+3\b}^2t_{\a+2\b}^3t_{2\a+3\b}+t_{\b}^2t_{\a+3\b}^2(t_{\a+2\b}+t_{\a+\b})^3t_{\a}+(t_{\b}+t_{\a+2\b})^2t_{2\a+3\b}^2t_{\a+\b}^3t_{\a}+\\ + & t_{\b}t_{\a+3\b}t_{2\a+3\b}t_{\a+\b}^2t_{\a}(3t_{\b}t_{\a+2\b}+2t_{\a+2\b}^2+2t_{\a+2\b}t_{\a+\b}+2t_{\b}t_{\a+\b}),\\  \pi_3=&  t_{\b}^3t_{\a+3\b}^2t_{\a+2\b}^3t_{2\a+3\b}+t_{\b}^3t_{\a+3\b}^2(t_{\a+2\b}+t_{\a+\b})^3t_{\a}+(t_{\b}+t_{\a+2\b})^3t_{2\a+3\b}^2t_{\a+\b}^3t_{\a}+\\
  +& t_{\b}^2t_{\a+3\b}t_{2\a+3\b}t_{\a+\b}^2t_{\a}(3t_{\b}t_{\a+2\b}+3t_{\a+2\b}^2+3t_{\a+2\b}t_{\a+\b}+2t_{\b}t_{\a+\b}),\\
  \pi_4=& t_{\b}^2t_{\a+3\b}^2t_{\a+2\b}^3t_{2\a+3\b}\left[t_{\b}t_{\a+3\b}t_{\a+2\b}^3t_{2\a+3\b}+2t_{\b}t_{\a+3\b}(t_{\a+2\b}+t_{\a+\b})^3t_{\a}+\right. \\ +&\left.(3t_{\b}t_{\a+2\b}+3t_{\a+2\b}^2+3t_{\a+2\b}t_{\a+\b}+2t_{\b}t_{\a+\b})t_{2\a+3\b}t_{\a+\b}^2t_{\a}\right]+\\
  + & t_{\a}^2(t_{\b}t_{\a+3\b}(t_{\a+2\b}+t_{\a+\b})^2+(t_{\b}+t_{\a+2\b})t_{2\a+3\b}t_{\a+\b}^2)^3
  \end{align*}
   for the changes of the normal order
     \begin{equation*}\begin{split}
      &\ldots, \a,\a\!+\b,2\a\!+3\b,\a\!+2\b,\a\!+3\b, \b,\ldots \to \ldots \b,\a\!+3\b,\a\!+2\b,
      2\a\!+3\b, \a\!+\b,\a, \ldots;\\  &\ldots,\b,\a\!+3\b,\a\!+2\b,2\a\!+3\b, \a\!+\b,\a,\ldots
      \to \ldots \a,\a\!+\b,2\a\!+3\b,\a\!+2\b,\a\!+3\b, \b,\ldots
      \end{split}\end{equation*} of $\mathrm{G}_2$ subsystem, see \cite[Theorem 3.1]{BZ}.

      Let $(,)$ be $W$-invariant  bilinear form on $\h^*$. We use the common notation
      $$\gamma^\vee=\frac{2\gamma}{(\gamma,\gamma)},\qquad \gamma^\vee\in\h^*$$
      for coroots. When identifying $\h$ and $\h^*$ by means of the form $(\,,\,)$ we denote them by $h_\gamma$, so that $\mu(h_\gamma)=(\mu, \gamma^\vee)$. This notation is in agreement with \rf{L2a}. Let $\mu\in\h^*$ be a weight (an integer weight), that is $\mu(h_i)\in\Z$  for any simple root $\alpha_i$. Set
       \begin{equation}\label{L7} \t^\mu=\prod_{\gamma\in\Delta_+}t_\gamma^{(\mu,\gamma^\vee)}
       \end{equation}
       \begin{lemma}\label{lemma1} The product \rf{L7} does not depend on the choice of the reduced decomposition of $w_0$.
       \end{lemma}
Thus $\t^\mu$ is a well defined rational function on $N_+$.

{\bf Proof}. This is a direct consequence of \cite[Theorem 4.3]{BZ}, which states that the matrix element
 $\Delta_k(g)=(v_k^+,gv_k^-)$ (generalized minor), restricted to $N_+$, admits a presentation
 \begin{equation*}
 \Delta_k(X(t))=\t^{\omega_k}
 \end{equation*}
  where $(\omega_k,\alpha_i^\vee)=\delta_{ik}$ for each simple root $\alpha_i$. Here $v_k^+$ and $v_k^-$ are highest and lowest weight vectors of the fundamental representation $V_{\omega_k}$. Another way to see that is to notice that the transition maps \rf{L4}, \rf{L5} and \rf{L6} leave invariant the following  monomials:
  \begin{align*}& \t^\a=t_\a t_{\a+\b} && \text{and}&&\t^\b= t_\b t_{\a+\b} &&\text{for}\qquad \rf{L4},\\
   & \t^\a=t_\a t_{\a+\b}^2t_{\a+2\b} && \text{and}&& \t^\b=t_\b t_{\a+2\b}t_{\a+\b} &&\text{for}\qquad \rf{L5},\\
   & \t^\a=t_\a t_{\a+\b}^3t_{2\a+3\b}^{2}t_{\a+2\b}^{3}t_{\a+3\b}^{} && \text{and}&& \t^\b=t_\b t_{\a+3\b}^{}t_{\a+2\b}^{2}t_{2\a+3\b}^{} t_{\a+\b}^{} &&\text{for}\qquad \rf{L6}
\end{align*}
which ensure the invariance of the product \rf{L7} under all transition maps.
\begin{lemma} \label{lemma2} (see {\rm\cite[Proposition 2.1]{GLO}}) The invariant measure $dn$ on $N_+$ is
\begin{equation}\label{L9}
dn=\t^{\rho}\prod_{\gamma\in\Delta_+}\frac{dt_\gamma}{t_\gamma}
\end{equation}
\end{lemma}
{\bf Proof} \cite{GLO}. First one checks by direct calculation that the measure
$$\frac{d\t}{\t}=\prod_{\gamma\in\Delta_+}\frac{dt_\gamma}{t_\gamma}$$
 is invariant with respect to transition maps  \rf{L4}, \rf{L5} and \rf{L6}. Thus the measure $\t^\rho\frac{d\t}{\t}$ is invariant as well. Second we check that this measure is invariant with respect to multiplication of the nilpotent matrix $X(\t)$ by group element $\exp\left(se_{\alpha}\right)$ from the right, where $\alpha$ is arbitrary simple root.
 \begin{equation}\label{L10}X(\t)\to X(\t)\cdot \exp\left(se_{\alpha}\right)
 \end{equation}
To this end  we choose a normal ordering which ends by the simple root $\theta(\alpha)$ (such ordering surely  exists, see \cite{Zh2}). In these Lusztig coordinates the map \rf{L10} becames a translation $t_{\theta(\alpha)}\to t_{\theta(\alpha)}+s$, the measure in the right hand side of  \rf{L9} factorises to the product
$\omega\wedge dt_{\theta(\alpha)}$, where $\omega$ does not depend on $t_{\theta(\alpha)}$ and thus is invariant under  that translation .
\subsection{Structure of Whittaker vectors}
To describe Whittaker vectors, we need some more invariants of transition maps \rf{L4}--\rf{L6}.
\begin{proposition}\label{prop1} ${}$

\noindent
{\em (i)} The sum $t_\a+t_{\a+\b}+t_\b$ is invariant with respect to transition map \rf{L4};\\
{\em (ii)} The sums $t_\a+t_{\a+2\b}$ and $t_{\a+\b}+t_\b$ are invariant with respect to transition map \rf{L5};\\
{\em (iii)} The sums $t_\a+t_{2\a+3\b}+t_{\a+3\b}$ and $t_{\a+\b}+t_{\a+2\b}+t_\b$ are invariant with respect to transition map \rf{L6}.
\end{proposition}
{\bf Proof}. Direct Maple check. \hfill{$\square$}
 \begin{corollary}\label{corollary1} The sum $\!\sum\limits_{\gamma\in\Delta_+}\!t_\gamma$ does not depend on the choice of coordinates $t_\gamma$.
 \end{corollary}
 Moreover, the sums of coordinates over the roots of the same length are invariant as well due to statements (ii) and (iii) of Proposition \ref{prop1}.

 Since Whittaker vectors are  functions on $G_0$ from the spaces $\V_\mu$ and $\V_{\nu}$, they are completely determined by their restrictions to the subgroup $N_+$. Denote respectively the restriction  $\wr_\mu(X(\t))$ of the functions $\wr_\mu(g)$ to $N_+$ by $\wwr_\mu(\t)$  and the restriction  $\wl_{\nu}(X(\t))$ of the functions $\wl_{\nu}(g)$ to $N_+$ by $\wwl_{\nu}(\t)$

\begin{proposition}\label{prop2} The right Whittaker vector $\wr_\mu$ is given by the function
$\wwr_\mu(\t)$ on $N_+$
\begin{equation*}
 \wwr_\mu(\t)=\exp\Big(-\sum_{\gamma\in\Delta_+}t_\gamma\Big)
\end{equation*}
\end{proposition}
{\bf Proof}. Choose a simple root $\alpha$. The group element $g\exp(se_\alpha)$ has the same Gauss components from $N_-$ and $H$. Thus it is enough to prove that for any simple root $\alpha$
\begin{equation*}
\wr_\mu(X(\t)\exp(se_\alpha))=e^{-s}\wr_\mu(X(\t))
\end{equation*}
Again we choose a normal ordering which ends by $\theta(\alpha)$. Then    $X(\t)\exp(se_\alpha)$ has the same coordinates as $X(\t)$ except $t_{\theta(\alpha)}$ which shifts to $t_\theta(\alpha)+s$. Then  the sum
$\sum_{\gamma\in\Delta_+}t_\gamma$ changes to $s+\sum_{\gamma\in\Delta_+}t_\gamma$  and the function
$\wwr_\mu(\t)$ transmits to $e^{-s}\wwr_\mu(\t)$ \hfill{$\square$}
\medskip

Note also that  for any real number $\epsilon$ the function $\wwr_\mu(\epsilon\t)$ gives rise to a Whittaker vector with a renormalized character
$\zeta^\epsilon_+$, see \rf{0}. In particular, the function
 $\bar{v}_\mu^R(g)$ on $g\in G_0$, where $g=YTX(t)$, $Y\in N_-$, $T\in H$, $X(\t)\in N_+$
\begin{equation*}
\bar{v}_\mu^R(g)= \chi_\mu(H) \wwr_\mu(X(-\t)) \end{equation*}
 is the Whittaker vector from the space $V_\mu$ with respect to the character $\zeta_+^{-1}$, that is
 $$\bar{v}_\mu^R(g\exp(s e_\a))=e^s\cdot\bar{v}_\mu^R(g)$$
for each simple root $\alpha$.

 We use in the following the standard lifts $\ss_{\a_i}$ of simple reflections $s_{\a_i}\in W$ to the group elements which we denote by the same letter,
\begin{equation}\label{W2a}\ss_{\a_i}=\exp(e_i)\exp(-f_i)\exp(e_i)\end{equation}
 so that the products $\w=\ss_{\a_{i_1}}\cdots \ss_{\a_{i_l}}$ do not depend on the choice of reduced decomposition of any element $w\in W$. Then this  lift to $G$ of the longest element group element $w_0$
has the property
\begin{equation}
\label{W4} \w_0e_\a \w_0^{-1}=\Ad_{\w_0}(e_\a)=-f_{-w_0(\a)},\qquad  \w_0f_\a \w_0^{-1}=\Ad_{\w_0}(f_\a)=-e_{-w_0(\a)}
\end{equation}
 for each simple root $\alpha$. This fact is a direct consequence of Tits's result
 \cite[Proposition 2.1 (3)]{T} which states that if for a simple root $\alpha_i$ and $w\in W$ there exists
  a simple root $\alpha_j$ such that  $ws_{\alpha_i}w^{-1}=s_{\alpha_j}$, then the same is true for their lifts. See e.g. \cite[Lemme 4.9]{BB}.
\begin{proposition}\label{prop3} The function
\begin{equation*}
\wl_{\nu}(g)=\bar{v}_\nu^R(g\w_0)\end{equation*}
is the left Whittaker vector in the space $\V_{\nu}$ and character $\zeta_-$.
\end{proposition}
{\bf Proof}. By the construction, the function $u(g)=\bar{v}_{\nu}^R(g\w_0)$ belongs to the space $\V_{\nu}$. Choose a simple root $\alpha_i$ and consider the function $u(g\exp(\ve f_\a))$. We have
\begin{equation*}\begin{split}
u(g\exp(s f_i))= \bar{v}_{\nu}^R(g\exp(s f_i)\w_0)=&\bar{v}_{\nu}^R(g\w_0\exp(-s e_i))=\\
&\exp(-s)\cdot \bar{v}_{\nu}^R(g\w_0))= \exp(-s)\cdot u(g).\qquad \qquad \hfill{\square}\end{split}
\end{equation*}
The left Whittaker vector, as a function on $G_0$, is completely determined by its restriction $\wwl_{\nu}$ to $N_+$. One of the  goals of this work is a proper description of this restriction, applicable for the study of Mellin transform. The proposition \ref{prop3} describes this restriction as follows. We start with the group element
$X(-\t)\w_0$. Take its Gauss decomposition
$$ X(-\t)\w_0=Y\cdot T\cdot \bar{X},\qquad Y\in N_-,\ T\in H,\ \bar{X}\in N_+.$$
 Let $p_\gamma(\t)$ be Lusztig coordinates of the element $\bar{X}\in N_+$.
 Then the restriction of the left Whittaker vector to $N_+$ is given by the function
 \begin{equation*}
 \wwl_{\nu}(\t)=\chi_{\nu}(T)\cdot\exp\Big(-\sum_{\gamma\in\Delta_+}p_\gamma(\t)\Big)
\end{equation*}
 Denote by  $g_\pm$ and $g_0$ the parts of Gauss decomposition of the group element $g$,
 $$g=g_-g_0g_+,\qquad\qquad g_\pm\in N_\pm,\ g_0\in H,$$
 and by $g_{0+}$ the product $g_0g_+$.
We see that the problem of writing precise  expressions for the left Whittaker vector leads to the study of  birational isomorphism $\tau$ of the manifold $N_+$ and the map $\sigma:N_+\to H$ given by
\begin{equation*}
\tau(X(\t))= \left(X(-\t)\w_0\right)_+,\qquad \sigma(X(\t))=\left(X(-\t)\w_0\right)_0
\end{equation*}
Another possibility, is to use for the construction of left Whittaker vector slightly different maps
$\tilde{\tau}:N_+\to N_+$  and  $\tilde{\sigma}:N_+\to H$ given by
\begin{equation*}
\tilde{\tau}(X(\t))= \left(X^{-1}(\t)\w_0\right)_+,\qquad \sigma(X(\t))=\left(X^{-1}(\t)\w_0\right)_0.
\end{equation*}
Below we present explicit description of the maps $\tau$ and $\sigma$ in Lusztig coordinates and use this description for the derivation of Mellin transforms of Whittaker functions.
\bigskip

The  maps $\tilde{\tau}$ and $\tau$ are closely related to the birational transform $\eta_{\w_0}$ by Berenstein and Zelevinsky. The latter plays the crucial role in their derivation of the factorized expression of Lusztig coordinates via generalized minors \cite[Theorem 1.4]{BZ}. The map $\eta_{\w_0}$ is
 \begin{equation}\label{W8}
 \eta_{\w_0}(X)= \left(\w_0X^T\right)_+
\end{equation}
where ${}^T$ is an anti-automorphism of $G$, trivial on $H$, and satisfying the relations  $$(\exp(te_\alpha))^T=\exp(tf_\alpha)), \qquad\text{and}\qquad
(\exp(tf_\alpha))^T=\exp(te_\alpha))$$ for each simple root $\alpha$.
 Compute first $\w_0X(\t)\w_0^{-1}$. We have
\begin{equation}\label{W10}\w_0X(\t)\w_0^{-1}=\w_0\prod\limits_{k}^{\leftarrow}\exp\left(t_{\gamma_k}f_{\alpha_{i_k}}\right)\w_0^{-1}=
\prod\limits_{k}^{\leftarrow}\exp\left(-t_{\gamma_k}e_{-w_0(\alpha_{i_k})}\right)\end{equation}
This expression coincides with $(\theta(X))^{-1}(\t)$.
We then have
\begin{equation}\label{W9}\eta_{\w_0}(X)=(\theta(X)\w_0)_+= \tilde{\tau}(\theta(X))
\end{equation}
so that the maps $\eta_{\w_0}$ and $\tilde{\tau}$ differ by Dynkin automorphism $\theta$. The relation
\rf{W10} implies also slightly more complicated relation between $\eta_{\w_0}$ and $\tau$,
\begin{equation}\label{W11}\eta_{\w_0}(X)= {\tau}(\iota\theta(X))
\end{equation}
where $\iota$ is the birational automorphism of the manifold $N_+$, reversing the order of the product in
Lusztig presentation \rf{L3} of the group element of $N_+$. Note that the map $\iota$ defined first for a given normal ordering, respects the transition maps \rf{L4}--\rf{L6} and thus does not depend on a choice of Lusztig coordinates.

Having in mind relations \rf{W9} and \rf{W11} we further use the name Berenstein--Zelevinsky (BZ) transform for the birational map $\tau: N_+\to N_+$ and for the corresponding change of Lusztig coordinates. The map $\sigma$ will be referred as (Cartan) twist.

\subsection{Using Plancherel formula}
 The convolution property for the Mellin transform admits the following reading.
 Let $F_1(s)$ and $F_2(s)$ be Mellin transforms of integrable functions $f_1(t)$ and $f_2(t)$,
 $$ F_i(s)=\int_0^\infty f_i(t)t^s\frac{dt}{t}.$$
 Let $S_1$ and $S_2$ be the strips of analyticity of $F_1(s)$ and $F_2(s)$,
 $$S_i:\ \a_i<\Re s<b_i.$$
 Assume that the intersection of strips $-S_1$ and $S_2$ is nonzero and both functions $F_1(s)$ and $F_2(s)$ rapidly decrease when $s$ goes to $\pm i\infty$
 along the contour $C=\{\Re s=c\}\subset -S_1\cap S_2$.
 Then
 \beq \label{Pl1} \int_0^\infty f_1(t)f_2(t)\frac{dt}{t}=\frac{1}{2\pi i}\int_{C}F_1(-s)F_2({s})ds\eeq
 Here is a check:
 \beqq \begin{split}& \frac{1}{2\pi i}\int\limits_{c-i\infty}^{c+i\infty}F_1(-s)F_2(s)ds=
\frac{1}{2\pi i}\int\limits_{c-i\infty}^{c+i\infty} F_1(-s)ds\int\limits_0^\infty f_2(t)t^s\frac{dt}{t}=\\
\frac{1}{2\pi i}\int\limits_0^\infty&  f_2(t)\!\frac{dt}{t}\int\limits_{c-i\infty}^{c+i\infty}t^sF_1(-s)ds=
\frac{1}{2\pi i}\int\limits_0^\infty f_2(t)\frac{dt}{t}\!\int\limits_{-c-i\infty}^{-c+i\infty}t^{-s}F_1(s)ds=\!\int\limits_0^\infty f_1(t)f_2(t)\frac{dt}{t}.
\end{split}\eeqq
The first equality is due to the inclusion $C\subset S_2$, the second uses the decreasing at infinities,  the third is the change of variables $s\to -s$,
the fourth is the Mellin inversion theorem due to
the inclusion $-C\subset S_1$.  Since Mellin transform of the function $\bar{f}_1(t)$ equals to $\bar{F}_1(\bar{s})$ , we can rewrite the convolution property
\rf{Pl1} in a form of Plancherel formula
\beq \label{Pl2}
 \int_0^\infty \bar{f}_1(t)f_2(t)\frac{dt}{t}=\frac{1}{2\pi i}\int_{C}\bar{F}_1(-\bar{s})F_2({s})ds
\eeq
under the same assumptions on the functions $F_1(s)$ and $F_2(s)$.

 Denote by $\wwlm_{\nu}(\bm{\gamma})$ and $\wwrm_\mu(\bm{\gamma})$ the Mellin transforms of
  $\t^\rho\wwl_{\nu}(\t)$ and $\wwr_\mu(\t)$ as of functions on $N_*$,
   \begin{equation*}
   \begin{split}& \wwrm_\mu(\bm{\gamma})=\int_{\t>0}\wwr_\mu(\t) \t^{\bm{\gamma}}\frac{d\t}{\t},\qquad
    \wwlm_{\nu}(\bm{\gamma})=\int_{\t>0}\t^\rho\wwl_{\nu}(\t) \t^{\bm{\gamma}}\frac{d\t}{\t}
   \end{split}
   \end{equation*}
   Here $\t^{\bm{\gamma}}$ means the product
   $$\t^{\bm{\gamma}}=\prod_{\gamma\in\Delta_+}(t_\gamma)^\gamma,$$
   where $(t_\gamma)^\gamma$ is the power of the variable $t_\gamma$, and $\t^\rho$ is defined by \rf{L7}.
We now use the Plancherel formula \rf{Pl2} to rewrite Whittaker wave functions \rf{5} and \rf{5a} in terms of
$\wwlm_{\nu}(\bm{\gamma})$ and $\wwrm_\mu(\bm{\gamma})$.
First we note that by definition \rf{4}, \rf{4a} of the action of the Lie algebra $\g$ in the space $V_\mu$, the restriction of the function
$\exp{(-x)}\wr_\mu$ for any $x\in\h$ is given by the relation
\beqq \begin{split} \exp{(-x)}\wr_\mu|_{N_+}=\wwr_\mu(\t) \exp{(-x)}|_{N_+}=\\
\exp{(-x)}\left(\exp{(x)} \wwr_\mu(\t)\exp{(-x)}\right)= e^{-(x,\mu)}e^{\mathrm {ad}_x}\left(\wwr_\mu(\t)\right)\end{split}\eeqq
 Mellin transform turns first order differential operators $\mathrm {ad}_x$ into operators of multiplications on functions $\tilde{H}(\x,\gamma)$, see
 \rf{gl31}, \rf{d20a}, \rf{b17}, \rf{c12}. Then by \rf{Pl2} and Lemma \ref{lemma2} we have
 \beq\label{Pl4}\begin{split}\tilde{\Psi}_{\mu,\nu}(x)=&e^{-(\rho+\mu,x)}\int_{N_*}\wwl_\nu(\t)e^{\mathrm {ad}_x}\left(\wwr_\mu(\t)\right)\frac{d\t}{\t}= \\ &
 \frac{e^{-(\rho+\mu,x)}}{(2\pi i)^d}\int_C \wwwlm_{\nu}(-\bar{\bm{\gamma}})\exp \tilde{H}(x,\bm{\gamma})\wwrm_\mu(\bm{\gamma})d\bm{\gamma}
 \end{split}\eeq
and
\beq \label{Pl3}
{\Psi}_{\lambda}(x)=
 \frac{e^{-i(\lambda,x)}}{(2\pi i)^d}\int_C \wwwlm_{i\lambda-\rho}(-\bar{\bm{\gamma}})\exp \tilde{H}(x,\bm{\gamma})\wwrm_{i\lambda-\rho}(\bm{\gamma})d\bm{\gamma}
\eeq
Here $d$ is the dimension of $N_+$, the contour $C$ is a  deformation of imaginary plane $\Re \bm{\gamma}=0$ into the intersection of strips of analyticity of
$\wwrm_\mu(\bm{\gamma})$ and  $\wwwlm_{\nu}(-\bar{\bm{\gamma}})$ under the assumption of nonemptiness
of their intersections and their vanishing on imaginary infinities.
\renewcommand\w[1]{w_0^{(#1)}}
\newcommand\ww[1]{\bar{w}_0^{(#1)}}
\section{$\GL(n,\R)$} \label{section2}
\subsection{ BZ transform}
\setcounter{equation}{0}
 Let $v_i$ be a fixed basis of $\R^n$ and $e_{i,j}\in \Mat_{n\times n}$ be  matrix units, $e_{i,j}(v_k)=\delta_{j,k}v_i$. Denote by $\ve_i\in\h^*$ the basic elements of $\h$ defined by the condition $\ve_i(e_{j,j})=\delta_{i,j}$ and by $s_i$ the group elements, related to simple reflection $s_{\ve_{i}-\ve_{i+1}}$, that is $s_i=(e_{i,i+1}-e_{i+1,i})+\sum_{j\not=i,i+1} e_{j,j}$ as element of  $\Mat_{n\times n}$. We use the following decomposition of the Weyl group element $w_0=\w{n}$,
 \beqq 
  \w{n}= s_{n-1}(s_{n-2}s_{n-1})\cdots (s_1s_2\cdots s_{n-1}).
 \eeqq
 Under the
 so that the corresponding ordering \rf{L2} of positive roots is
 \beqq
 \ve_{n-1,n},\ \ve_{n-2,n},\ve_{n-2,n-1},\ \ve_{n-3,n},\ve_{n-3,n-1},\ve_{n-3,n-2},\cdots,\ \ve_{1,n},\ve_{2,n}\cdots \ve_{1,2},\eeqq {where}  $\ve_{i,j}=\ve_i-\ve_j$.
Let $\Xn{n}(\t)$ be a group element of the subgroup $N_+$, given by \rf{L3},
 \beqq
 \begin{split}
 	\Xn{n}(\t)=\exp(t_{n-1,n}e_{n-1,n})\ \exp(t_{n-2,n}e_{n-2,n-1})\exp(t_{n-2,n-1}e_{n-1,n})\ \cdots \\
 \exp(t_{1,n}e_{1,2})\exp(t_{1,n-1}e_{2,3})\cdots \exp(t_{1,2}e_{n-1,n}).\end{split}
 \eeqq
 Here we regard $e_{i,i+1}$ as Chevalley generators of Lie algebra $\gl_n$. In fundamental representation
 they are given by the matrices
\beqq
\begin{split}
 \Xn{n}(\t)=(\Id+t_{n-1,n}e_{n-1,n})\ (\Id+t_{n-2,n}e_{n-2,n-1})(\Id+t_{n-2,n-1}e_{n-1,n})\ \cdots \\
  (\Id+t_{1,n}e_{1,2})(\Id+t_{1,n-1}e_{2,3})\cdots (\Id+t_{1,2}e_{n-1,n}).
  \end{split}
  \eeqq
   where $\Id\in\Mat_{n\times n}$ is the identity matrix.
    We are going to solve the following equation
   \beq \label{gl5}
  T\Xn{n}(\p)= \left(\Xn{n}(-\t)\ww{n}\right)_{0+}
   \eeq
 where
 \beq\label{gl6}
 \begin{split}
  \Xn{n}(\p)=(\Id+p_{n-1,n}e_{n-1,n})\ (\Id+p_{n-2,n-1}e_{n-2,n})(\Id+p_{n-2,n-1}e_{n-1,n})\ \cdots \\
   (\Id+p_{1,n}e_{1,2})(\Id+p_{1,n-1}e_{2,3})\cdots (\Id+p_{1,2}e_{n-1,n}).
   \end{split}
 \eeq
 and $T$ is a diagonal matrix. In other word, we have to find out the unknown variables $p_{i,j}$ as function of $t_{i,j}$ and the diagonal matrix $T$ depending on $p_{i,j}$ (and thus on $t_{i,j}$) such that relation \rf{gl6} holds.
  Denote further in this section
  $$\hat{k}=k+1-i.$$
 \begin{theorem}\label{theorem1}{\em{a)}} BZ map $\tau$, see \rf{W8} and\rf{gl5}, has the form
 	\begin{equation}\label{gl18}
 		p_{i,j}=\frac{1}{t_{i, i+\j }}\frac{t_{i-1,i+\j-1}}{t_{i-1,i+\j}}\frac{t_{i-2,i+\j-1}}{t_{i-2,i+\j}}\cdots
 		\frac{t_{1,i+\j-1}}{t_{1,i+\j}}
 	\end{equation}
{\em{b)}} The twisting matrix $T\in H$ reads
\begin{equation}\label{gl19}\begin{split}
T= t_{1,n}\cdots t_{1,2}\cdot e_{1,1}+&\frac{t_{2,n}\cdots t_{2,3}}{t_{1,2}}\cdot e_{2,2}+
		\frac{t_{3,n}\cdots t_{3,4}}{t_{1,3}t_{2,3}}\cdot e_{3,3}+\ldots\\
		+&\frac{t_{i,n}\cdots t_{i,i+1}}{t_{1,i}\cdots t_{i-1,i}}\cdot e_{i,i}+\ldots
		+\frac{1}{t_{1,n}\cdots t_{n-1,n}}\cdot e_{n,n}.
		\end{split}\end{equation}
{ \em{c)}} BZ transform preserves the measure $\dfrac{d\t}{\t}$,
 $$ \prod_{i<j}\frac{dt_{i,j}}{t_{i,j}}=\prod_{i<j}\frac{dp_{i,j}}{p_{i,j}}.$$	
 	\end{theorem}

 \noindent {\bf Remark 1}. BZ transform $\tau$ is involutive, that is the inverse relation
 \begin{equation}\label{gl20}
 t_{i,j}=\frac{1}{p_{i, i+\j }}
 \frac{p_{i-1,i+\j-1}}{p_{i-1,i+\j}}\frac{p_{i-2,i+\j-1}}{p_{i-2,i+\j}}\cdots
 \frac{p_{1,i+\j-1}}{p_{1,i+\j}}
 \end{equation}
 has the same form as \rf{gl18}. Indeed, the map $\tau$: $\Xn{n}(\t)\to \Xn{n}(\p)$ is defined by the relation $\Xn{n}(\p)=\left(\Xn{n}(-\t)\ww{n}\right)_+$, Since $\left(\ww{n}\right)^2=(-1)^{n+1}\mathrm{Id}$, this is equivalent to $\Xn{n}(-\t)=\left(\Xn{n}(\p)w_0\right)_+$. This means that we can interchange in \rf{gl18} $t_{i,j}$ to $-p_{i,j}$ and $p_{i,j}$ to $-t_{i,j}$. After cancelation of signs we get \rf{gl20}. \hfill{$\square$}
 \medskip

 \noindent {\bf Remark 2}. The Cartan twist $T$ can be as well  written in shorthand notation as
 $T(\t)=T(\p^{-1})$, that is
 \begin{equation}\label{gl21}\begin{split}
 T= \frac{1}{p_{1,n}\cdots p_{1,2}}\cdot e_{1,1}+&\frac{p_{1,2}}{p_{2,n}\cdots p_{2,3}}\cdot e_{2,2}+
 		\frac{p_{1,3}p_{2,3}}{p_{3,n}\cdots p_{3,4}}\cdot e_{3,3}+\ldots\\
 		+&\frac{p_{1,i}\cdots p_{i-1,i}}{p_{i,n}\cdots p_{i,i+1}}\cdot e_{i,i}+\ldots
 		+{p_{1,n}\cdots p_{n-1,n}}\cdot e_{n,n}.
 		\end{split}\end{equation}
Indeed, $T$ is specified by the condition
$$ T\Xn{n}(\p)=\left(\Xn{n}(-\t)\ww{n}\right)_{0+}.$$
Again, since $\left(\ww{n}\right)^2=(-1)^{n+1}\mathrm{Id}$, this is equivalent to
$$T^{-1}\Xn{n}(-\t)= (-1)^{n+1}\left(\Xn{n}(\p)\ww{n}\right)_{0+},$$
which implies the equality $T^{-1}(-\p)=(-1)^{n+1}T(\t)$. The cancelation of signs gives \rf{gl21}.
\hfill{$\square$}
\begin{corollary}\label{corollary2}
	 The restriction of left Whittaker vector to $N_+$ is given by the function
\begin{equation}\label{gl29}
	 \wwl_{\nu}(\t)=\t^{\nu}\cdot\exp\Big(-\sum_{i<j}p_{i,j}\Big)=\p^{-\nu}\cdot
\exp\Big(-\sum_{i<j}p_{i,j}\Big)
	\end{equation}
where $p_{i,j}$ are given by \rf{gl18}
\end{corollary}
The formula \rf{gl29} is well known, see \cite{GLO} an references therein. However, we got it in a form convenient for the study of Mellin transform. Further on we find analogous form for other classical groups.
\bigskip

The proof of Theorem \ref{theorem1} is based on inductive calculation of BZ transform $\tau$ and related map $\sigma$. Namely, the group element $\Xn{n}(-\t)$ admits a factorization
  \beqq
  \Xn{n}(-\t)=\tXn{n-1}(-\t)\cdot \An{n}(-\t)
  \eeqq
  where
   \beqq
    \An{n}(-\t)=
     (\Id-t_{1,n}e_{1,2})(\Id-t_{1,n-1}e_{2,3})\cdots (\Id-t_{1,2}e_{n-1,n}).
     \eeqq
   and
  \beqq
  \begin{split}
   \tXn{n-1}(-\t)=(\Id-t_{n-1,n}e_{n-1,n})\ (\Id-t_{n-2,n}e_{n-2,n-1})(\Id-t_{n-2,n-1}e_{n-1,n})\ \cdots \\
    (\Id-t_{2,n}e_{1,2})(\Id-t_{2,n-1}e_{2,3})\cdots (\Id-t_{2,3}e_{n-2,n-1}).
    \end{split}
    \eeqq
 represents a group element of the unipotent subgroup of embedded $\GL(n-1)$,
 $$  \tXn{n-1}(-\t)=\left(\begin{array}{c|c}1&0\\ \hline 0&\begin{array}{ccc}1&*&*\\0&\ddots&*\\0& 0& 1
 \end{array}\end{array}
 \right)$$
The matrix $\Xn{n}(\p)$ has the same structure,
 \begin{equation*}
  \Xn{n}(\p)=\tXn{n-1}(\p)\cdot \An{n}(\p).
  \end{equation*}
 In the induction step we first express Lusztig coordinates of $\An{n}(\p)$ via that of  $\An{n}(-\t)$;
 compute the input of $\An{n}(-\t)$ into Cartan twist $T$ and reduce the rest of calculations to the computation of BZ transform and twist of $\GL(n-1)$ matrix $\ttXn{n-1}(-t)$, which is given by a certain gauge transform of $\tXn{n-1}(-\t)$.
 \bigskip

 Denote by $\Ln{n}$ the diagonal matrix with diagonal entries
 \beq\label{gl22}
 \Ln{n}_{1,1}=t_{1,n}t_{1,n-1}\cdots t_{1,2},\qquad \Ln{n}_{j,j}=\frac{1}{t_{1,\j+1}}=\frac{1}{t_{1,n-j+2}},\ j>1
   \eeq
 and define the variables $\tilde{t}_{i,j}$, $1< i<j\leq n$ by the relation
 \beq\label{gl23}
 \tilde{t}_{i,j}=t_{i,j}\cdot\frac{t_{1,j-i+2}}{t_{1,j-i+1}}
 \eeq
 \begin{proposition}\label{prop4} {\em{a)}} The parameters $p_{1,j}$, $j=2,...,n$ of $\An{n}(\p)$ are equal to
 \beq\label{gl24}
 p_{1,j}=\frac{1}{t_{1,n-j+2}}\eeq
{\em{b)}} Cartan twist $T$ and the matrix $\tXn{n-1}(\p)$ satisfy the relation
\beq\label{gl25}
T \tXn{n-1}(\p)=\Ln{n}\cdot\left(\tXn{n-1}(-\tilde{t})\ww{n-1}\right)_{0+}
\eeq
 \end{proposition}

  {\bf Proof } of part a) of Proposition \ref{prop4}.
  Matrix element $A_{i,j}$ of upper triangular unipotent matrix $\An{n}(-\t)$
   is
   \beq\label{gl10}A_{i,j}^{(n)}= A_{i,j}^{(n)}(-\t)=(-1)^{j-i}t_{1,\j+1}t_{1,\j+2}\cdots t_{1,\i-1},t_{1,\i} \qquad i<j
   \eeq
    Due to the structure \rf{gl24} of the matrix $\tXn{n-1}(-\t)$ the first row of  $\Xn{n}(-\t)$ coincides with that of $\An{n}(-\t)$ and is equal to
     $$ ((-1)^{n-1}A_{1,n}^{(n)},(-1)^{n-2}A_{2,n}^{(n)},\ldots A_{1,1}^{(n)})=(t_{1,n}t_{1,n-1}\cdots t_{1,2},\
     \ldots\ ,t_{1,3}t_{1,2}, t_{1,2},1).$$
   Since the right multiplication by $\ww{n}=e_{1,n}-e_{2,n}+\ldots +(-1)^{n-1}e_{n,1}$ leaves the first row of the matrix $\Xn{n}(\t)$ invariant, permuting its entries, the parameters of $\An{n}(\p)$ and matrix element $T_{1,1}$ can be completely determined  from the parameters of $\An{n}(-\t)$. Namely,
  $T_{1,1}=t_{1,n}t_{2,n}\cdots t_{n-1,n}$ and we have two presentations for the first row of the matrix $\An{n}(\p)$:
 $$\left( 1,p_{1,n},p_{1,n}p_{1,n-1}, \cdots ,p_{1,n}p_{2,n}\!\cdots\! p_{1,2}                  \right)\ = \    \left(1,\ \frac{1}{t_{1,2}},\frac{1}{t_{1,2}t_{1,3}},\  \dots\ , \frac{1}{t_{1,2}t_{1,3}\cdots t_{1,n}}\right)$$
  so that the parameters $p_{i,n}$ of the matrix $\An{n}(\p)$ are equal to
 $
 p_{1,j}=\dfrac{1}{t_{1,n+2-j}}
 $.
 \hfill{$\square$}

 Let $\ww{n-1}$ be the  element of embedded group $\GL(n-1)$  representing the corresponding longest element of the Weyl group,  $\ww{n-1}=e_{1,1}+e_{2,n}-e_{3,n-1}+\ldots +(-1)^ne_{n,2}$. Let  $U$ be the matrix
 \beq\label{gl12}U=	\left(\ww{n-1}\right)^{-1}\An{n}(-\t)\ww{n}\left(\An{n}(\p)\right)^{-1}\eeq
  The following  lemma is the crucial technical step in the  calculation of the transformation
   \rf{gl5}. It states that the matrix $U$ has the block structure
   \beqq
   U=\left(
   \begin{array}{c|c}*&0\\ \hline *& \Lambda' \end{array}
   \right)
   \eeqq
   where $\Lambda'$ is $(n-1)\times(n-1)$ diagonal matrix. More precisely
   \begin{lemma}\label{lemma3}  Matrix elements $U_{i,j}$ equal zero for $i\not=j$, $j>1$.  Matrix element $U_{1,1}$ equals $\left(p_{1,2}p_{1,3}\cdots p_{1,n}\right)^{-1}=t_{1,2}t_{1,3}\ldots t_{1,n}$. The element $U_{i,i}$ equals $p_{1,\i+1}=t_{1,i}^{-1}$ for $i>1$.
    	\end{lemma}
    The proof	 of Lemma  \ref{lemma3} is given in Appendix A.

 {\bf Proof} of part b) of Proposition \ref{prop4}.
 Denote by ${\Lambda'}^{(n)}=\ww{n-1}\Lambda^{(n)}\left(\ww{n-1}\right)^{-1}$ the diagonal matrices  with nonzero entries
 \begin{equation*}
{\Lambda'}^{(n)}_{1,1}=t_{1,2}t_{1,3}\cdots t_{1,n},\qquad
{\Lambda'}^{(n)}_{i,i}=\frac{1}{t_{1,i}},\ i>1
  \end{equation*}
 Since $\An{n}(\p)\in N$, the relation \rf{gl5} can be equivalently rewritten as
 $$T \tXn{n-1}(\p)=\left( \tXn{n-1}(-\t) \An{n}(-\t)\ww{n}\left(\An{n}(\p)\right)^{-1}\right)_{0+}$$
 or
 \begin{equation*}
 \begin{split}T \tXn{n-1}(\p)=&\left(\tXn{n-1}(-\t)\ww{n-1}U\right)_{0+}=\left(\tXn{n-1}(-\t)\ww{n-1}{\Lambda'}^{(n)}V\right)_{0+}=\\ & \Ln{n}\left(\ttXn{n-1}(-\t)\ww{n-1}V\right)_{0+}, \end{split}\end{equation*}
where
\begin{equation}\label{gl26}V=\left({\Lambda'}^{(n)}\right)^{-1}U\qquad\text{ and}\qquad
\ttXn{n-1}(-\t)=\left(\Ln{n}\right)^{-1}\tXn{n-1}(-\t)\Ln{n}.\end{equation}
The matrices $Y=\ttXn{n-1}(-\t)\ww{n-1}$ and $V$ have the following block structure due to Lemma \ref{lemma3}:
\beqq
 Y=\left(\begin{array}{c|c} 1&0\\ \hline 0&\tilde{Y}\end{array}\right),
\qquad V=\left(\begin{array}{c|c} 1&0\\ \hline Z&\mathrm{Id}_{n-1}\end{array}\right)
\eeqq
where $\tilde{Y}$ and $\mathrm{Id}_{n-1}$ are $(n-1)\times (n-1)$ matrices and $Z$ is $1\times(n-1)$ matrix. Then the product $YV$ also has a block structure and can be rewritten as
$$YV=\left(\begin{array}{c|c} 1&0\\ \hline \tilde{Y}Z&\mathrm{Id}_{n-1}\end{array}\right)\cdot\left(\begin{array}{c|c} 1&0\\ \hline 0&\tilde{Y}\end{array}\right)$$
The first factor of the latter product is unipotent lower triangular matrix so that the upper triangular part of $YV$ coincides with upper triangular part of $Y$. Thus we proved the equality
\beqq
T \tXn{n-1}(\p)=\Ln{n}\left(\ttXn{n-1}(-{\t})\ww{n-1}\right)_{0+}
\eeqq
where the matrix $\ttXn{n-1}(-{\t})$ is given by the relation \rf{gl26}.

The last step is the computation of the conjugation in \rf{gl26}.  To perform it we note that the exponent in the products \rf{gl19}  containing the variable $t_{i,j}$ is
 \beq\label{gl27a}
 \exp\left(-t_{i,j}e_{n-j+i,n-j+i+1}\right)=\left(1-t_{i,j}e_{n-j+i,n-j+i+1}\right)
 \eeq
 which means that the conjugation \rf{gl26} results to multiplication of the coefficient at
 $e_{n-j+i,n-j+i+1}$ by ${\Ln{n}_{n-j+i+1,n-j+i+1}}/{\Ln{n}_{n-j+i,n-j+i}}$ which is due to
 \rf{gl22}, the rescaling
 $$t_{i,j}\to t_{i,j}\frac{t_{1,j-i+2}}{t_{1,j-i+1}}. $$
 ${}$\hfill{$\square$}

 {\bf Proof} of Theorem \ref{theorem1} follows from the inductive application of Proposition \ref{prop4}. At the first step we find the variables $p_{1,j}$, see \rf{gl24} and the input $\Ln{n}$ of the first step to the diagonal matrix $T$. In particular, we find that $T_{1,1}=t_{1,2}\cdots t_{1,n}$. Then we pass to the second step, where we deal with $(n-1)$ square matrix but with rescaled by \rf{gl23} matrix elements $\tilde{t}_{i,j}$. Here after the corresponding shift of indices  we find out the next portion of variables,
  $$
  p_{2,j}=\frac{1}{\tilde{t}_{2,n-j+3}}= \frac{1}{{t}_{2,n-j+3}}\frac{{t}_{1,n-j+2}}{{t}_{1,n-j+3}} $$
  due to \rf{gl23} and \rf{gl24}. The Cartan matrix $T$ gains the new income, equal to
  $\Ln{n-1}$ with matrix entries
  \begin{align*}
  \Ln{n-1}_{1,1}=&1,\qquad \Ln{n-1}_{2,2}=\tilde{t}_{2,3}\cdots \tilde{t}_{2,n}=
  {t}_{2,3}\cdots {t}_{2,n}\cdot\frac{t_{1,n}}{t_{1,2}},\\
  \Ln{n-1}_{j,j}=&\frac{1}{\tilde{t}_{2,n-j+3}}=
  \frac{1}{t_{2,n-j+3}}\frac{t_{1,n-j+2}}{t_{1,n-j+3}} ,\qquad j>2\end{align*}
  and 	the renormalization of the variables $\tilde{t}_{i,j}$ for new $(n-2)\times(n-2)$ task,
  $$\tilde{t}_{i,j}\to \tilde{\tilde{t}}_{i,j}=\tilde{t}_{i,j}\frac{\tilde{t}_{2,j-i+3}}{\tilde{t}_{2,j-i+2}}.$$
 Following this procedure we get both a) and b) statements of the Theorem.
  The part c) may be observed as follows. Denote by $\Omega_t$ and $\Omega_p$ the skew forms $\Omega_t=\wedge_{i<j}\frac{dt_{i,j}}{t_{i,j}}$ and $\Omega_p=\wedge_{i<j}\frac{dp_{i,j}}{p_{i,j}}$. They admit the factorizations
  $\Omega_t=\Omega'_t\wedge\Omega''_t$, and $\Omega_p=\Omega'_p\wedge\Omega''_p$, where
  $$\Omega'_t=\wedge_{1<j}\frac{dt_{1,j}}{t_{1,j}},\qquad \Omega''_t=\wedge_{1<i<j}\frac{dt_{i,j}}{t_{i,j}},$$
  and analogously for $\Omega'_p$ and $\Omega''_p$
   Looking at the induction step and relation \rf{gl24} we see that the skew forms $\Omega'_t$ and $\Omega'_p$ coincide up to sign. But then the relations \rf{gl23} say that in the wedge product
   $$\Omega'\wedge_{1<i<j}\frac{d\tilde{t}_{i,j}}{\tilde{t}_{i,j}}$$ the renormalization fractions $\frac{t_{1,j-i+2}}{t_{1,j-i+1}}$ should be regarded as constants which so do not contribute to the wedge product, so that
   $$\Omega''=\wedge_{1<i<j}\frac{d\tilde{t}_{i,j}}{\tilde{t}_{i,j}}$$
   and we may further use the same equalities for the next induction steps. \hfill{$\square$}

   \subsection{Whittaker function} \label{section2.2}
   Using the relation  \rf{gl27a} we immediately describe the action of Cartan
   generators on the right Whittaker vector. By  definition,
   $\exp(\sum_k - x_k e_{k,k})\cdot f(g)=f(g\exp(\sum_k - x_k e_{k,k}))$
    for any function $f:G\to\C$, so the vector $\exp(-\sum_k x_k e_{k,k})\wr_\mu$
    is presented by the  function  $\exp(-\sum_k \mu_k x_k)\wwr_\mu(\t')$, where
\beq\label{gl30}
t'_{i,j}=t_{i,j}\cdot\exp\left(x_{n-j+i}-x_{n-j+i+1}\right)
\eeq
  Denote by $\wwlm_{\nu}(\bm{\ga})$ and $\wwrm_\mu(\bm{\ga})$ the Mellin transforms of the functions on $N_+$
  $\t^\rho\wwl_{\nu}(\t)$ and $\wwr_\mu(\t)$,
   \begin{equation*}
   \begin{split}& \wwrm_\mu(\bm{\ga})=\int_{\t>0}\wwr_\mu(\t) \t^{\mathbf{\ga}}\frac{d\t}{\t}:=
   \int_{t_{i,j}>0} \wwr_\mu(\t)\prod_{i<j}t_{i,j}^{\ga_{i,j}}\frac{dt_{i,j}}{t_{i,j}},\\
   & \wwlm_{\nu}(\bm{\ga})=\int_{\t>0}\t^\rho\wwl_{\nu}(\t) \t^{\mathbf{\ga}}\frac{d\t}{\t}:=
   \int_{t_{i,j}>0} \wwl_{\nu}(\t)\prod_{i<j}t_{i,j}^{\ga_{i,j}+j-i}\frac{dt_{i,j}}{t_{i,j}}.
   \end{split}
   \end{equation*}
 Here $\rho=-\sum_{k=1}^nk\ve_k$.
 Due to \rf{gl30}the action of Cartan subgroup on the right Whittaker vector $\wr_\mu$ transforms
    the function $\wwrm_\mu(\bm{\ga})$ to the product
     \beqq
   \exp\sum_k- x_k e_{k,k}: \ \wwrm_\mu(\bm{\ga})\mapsto
   \exp(-\sum_k \mu_k x_k)\cdot\exp \tilde{H}(\x,\bm{\ga})\cdot \wwrm_\mu(\bm{\ga})
   \eeqq
   where
   \beq   \label{gl31}
   \tilde{H}(\x,\bm{\ga})=\sum_{i<j}\ga_{i,j}(x_{n-j+i+1}-x_{n-j+i})
   \eeq

   Due to Proposition \ref{prop2} we have:
   \beqq
   \wwrm_\mu(\bm{\gamma})=\prod_{i<j}\Gamma(\gamma_{i,j}).
   \eeqq
   For the calculation of $\wwlm_{\nu}(\bm{\gamma})$ we pass in the integral
   $$\int_{t_{i,j}>0}\exp \Big(-\sum_{i<j}p_{i,j}\Big)\cdot\t^{\nu'}\prod_{i<j}t_{i,j}^{\ga_{i,j}}\frac{d\t}{\t},$$
   where $$\nu'_i=\nu_i-i$$
   from the integration variables $t_{i,j}$ to $p_{i,j}$.
   Substitution of \rf{gl20} gives the relation
   $$\t^{\nu}\prod_{i<j}t_{i,j}^{\ga_{i,j}}:=\prod_{i<j}t_{i,j}^{\ga_{i,j}+\nu'_i-\nu'_j}=
   \prod_{i<j}p_{i,j}^{-\vf_{i,j}+\nu'_j-\nu'_i},$$
   where
   \beqq 
   \vf_{i,j}= -\sum_{k>i}\ga_{k,k+\hat{j}-1}+\sum_{k\geq  i}\ga_{k,k+\hat{j}}
   \eeqq
    This implies the relation
   \beqq
   \wwlm_{\nu}(\bm{\gamma})=\prod_{i<j}\Gamma(-\vf_{i,j}+\nu'_j-\nu'_i),\qquad
   \wwwlm_\nu(-\bar{\bm{\gamma}})=\prod_{i<j}\Gamma(\vf_{i,j}+{\nu}'_i-{\nu}'_j)
   \eeqq
   Then by \rf{Pl4} and \rf{Pl3} we have
   \beqq 
  \tilde{ \Psi}_{\mu,\nu}(\x)= \frac{e^{-(x,\rho+\mu)}}{(2\pi i)^d}\!\int_{C}
   \exp \tilde{H}(\x,\bm{\ga}) \prod_{i<j}
   {\Gamma}(\vf_{i,j}+{\nu}'_i-{\nu}'_j)\Gamma(\ga_{i,j})d\ga_{i,j}.
   \eeqq
    where now $\nu'+\bar{\mu}=-\rho$,  $d=n(n-1)/2$ and the contour $C$ is a  deformation of
     imaginary plane $\Re\gamma_{ij}=0$ to the strip of analyticity of the integrand.
     In particular, for $\mu=i\lambda-\rho$ we have
    \beq \label{gl36}
    { \Psi}_{\lambda}(\x)= \frac{e^{-i(x,\lambda)}}{(2\pi i)^d}\int_{C}
    \exp \tilde{H}(\x,\bm{\ga}) \prod_{k<l}
    {\Gamma}(\vf_{k,l}+i({\lambda}_k-{\lambda}_l))\Gamma(\ga_{i,j})d\ga_{i,j}.
    \eeq

Performing the following change of variables in the integral in \rf{gl36}:
\beqq
\tgamma_{i,j}= \gamma_{i,j}+\gamma_{i+1,j+1}+\ldots +\gamma_{n-j+i,n}.\eeqq
we finally arrive to
\begin{theorem}\label{t1}
\beqq
\begin{split}
	\Psi_{\lambda}(\x)=&\frac{e^{-i(\lambda,x)}}{(2\pi i)^d}\int_C
	\exp\sum_{1\leq k\leq n
	}({\gamma}_{1,n+2-k}-\gamma_{1,n+1-k})x_{k}
	\\&
	\prod_{k<l}{\Gamma}(\gamma_{kl}-\gamma_{k+1,l}+i(\lambda_k-\lambda_{n+k-l+1}))
	\Gamma(\ga_{kl}-\gamma_{k+1,l+1})d\gamma_{k,l}
\end{split}\eeqq
\end{theorem}
Here we assume that $\gamma_{k,l}=0$ unless $1\leq k<l\leq n$.
 The integration cycle $C_n$ is a deformation of the imaginary plain
$\Re \gamma_{k,l}=0$ into nonzero strip   $D\subset\C^{\d}$ of the analyticity of the integrand, which
 can be described by inequalities
 \beqq \Re \gamma_{k,l}>0, \Re \gamma_{k,l}>\Re \gamma_{k+1,l} \eeqq
for all admissible pairs $(k,l)$ of indices.

\def\wb{\bar{w}}

   \setcounter{equation}{0}

   \section{$\SO(n,n)$}
   \subsection{BZ transform}
   The split real form of the group $\SO(2n,\C)$ is the group $\SO(n,n)\subset \SL(2n,\R)$ preserving symmetric form $\sum_{i=1}^{2n} x_iy_{\hat{i}}$, where as before $\hat{i}=2n+1-i$. Gauss decomposition is induced from that of $\GL(2n,\R)$. Positive roots are  $\ve_i\pm\ve_j$ for $1\leq i<j\leq n$, where
    $\ve_k\in\h^*$ are now defined by the condition $\ve_k( e_{j,j}-e_{\hat{j},\hat{j}})=\delta_{j,k}$.
   We denote Chevalley generator of the Lie algebra $\g=\so(2n)$ by $e_i$ and $f_i$, $i=1,\ldots,n-2$ and $e_{n-1}^\pm$, $f_{n-1}^\pm$. Here
   \begin{align*} &e_i=e_{i,i+1}-e_{2n-i,2n+1-i},&&  f_i=e_{i+1,i}-e_{2n+1-i,2n-i},&&
   i=1,...,n-2,\\
   & e_{n-1}^+=e_{n-1,n}-e_{n+1,n+2},&& f_{n-1}^+=e_{n,n-1}-e_{n+2,n+1},&&\\
   & e_{n-1}^-=e_{n-1,n+1}-e_{n,n+2},&& f_{n-1}^+=e_{n+1,n-1}-e_{n+2,n}.&&
   \end{align*}
   Denote by $s_i$ $i=1,...,n-1$ and $s_{n-1}^\pm$ the corresponding generators of the Weyl group and
    their lifts to the group $\SO(n,n)$ according to \rf{W2a}.
    We choose the following normal ordering of the system $\Delta_+$ of positive roots:
    \beqq
    \begin{split}&(\ve_{n-1}+\ve_n,\, \ve_{n-1}-\ve_n),\ (\ve_{n-2}+\ve_{n-1},\ve_{n-2}+\ve_{n},\ve_{n-2}-\ve_{n},\ve_{n-2}-\ve_{n-1}),\cdots,\\
    &(\ve_{1}+\ve_{2},\ve_{1}+\ve_{3},\cdots,\ve_{1}+\ve_{n},\ve_{1}-\ve_{n},\cdots, \ve_{1}-\ve_{3},\ve_{1}-\ve_{2}).
    \end{split}
     \eeqq
    It corresponds by \rf{L2} to the following reduced decomposition of the longest element $w_0^{(n)}$ of the Weyl group $W$:
    \beqq
    w_0^{(n)}=({s}^-_{n-1}s_{n-1}^+)\,(s_{n-2}s_{n-1}^+{s}_{n-1}^-s_{n-2})\,
       \cdots\,(s_{1}s_2...s_{n-2}{s}^{\eps(n-1)}_{n-1}s^{\eps(n)}_{n-1}s_{n-2}...s_2s_1).
    \eeqq
   where $\eps(k)=(-1)^k$. Note that as element of $\SO(n,n)$,
   \beq\label{d4}\begin{split}w^{(n)}_0=&-\sum_{k=1}^{2n}e_{k,\hat{k}}\qquad\qquad\qquad\qquad\qquad\text{for even}\ \ n,\\
  \wb_0^{(n)}=&\sum_{k\not=n,n+1}e_{k,\hat{k}}+\left(e_{n,n}+e_{n+1,n+1}\right)\qquad\text{for odd}\ \ n.
   \end{split}\eeq
     Denote by $t_{i,j}$ the Lusztig parameter corresponding to the root $\ve_i-\ve_j$, and  by $s_{i,j}$ the Lusztig parameter corresponding to the root $\ve_i+\ve_j$. Here $1\leq i<j\leq n$. Then the group element $\Xn{n}(\t)$ looks as
   \begin{align}\label{d3}
    \Xn{n}(\t)=&\left(\exp(s_{n-1,n}e_{n-1}^{-})\exp(t_{n-1,n}e_{n-1}^{+})\right)\cdot \notag\\
    	&\qquad\qquad\ldots\\
    	\left( \exp(s_{2,3}\right.& \left. e_2)\cdots
    	\exp(s_{2,n}e_{n-1}^{\eps(n-2)})\exp(t_{2,n}e_{n-2}^{\eps(n-1)})\cdots
    	\exp(t_{2,4}e_{3})\exp(t_{2,3}e_2)\right)\cdot	\notag \\
    	 \left( \exp(s_{1,2}\right.& e_1)\exp(s_{1,3}e_2)\cdots \left.
    	 	\exp(s_{1,n}e_{n-1}^{\eps(n-1)})\exp(t_{1,n}e_{n-1}^{\eps(n)})\cdots
    	 	\exp(t_{1,3}e_{2})\exp(t_{1,2}e_1)\right).	\notag
    	\end{align}
   BZ transform for $\SO(n,n)$ is the solution of the equation \rf{gl5}, where in $\Xn{n}(\p)$ we use the notation
   $p_{i,j}$ for Lusztig parameter, corresponding to the root $\ve_i-\ve_j$, and $q_{i,j}$ for the Lusztig parameter, corresponding to the root $\ve_i+\ve_j$,
    \begin{align*}
        \Xn{n}(\p)=&\left(\exp(q_{n-1,n}e_{n-1}^{-})\exp(p_{n-1,n}e_{n-1}^{+})\right)\ \cdots\notag\\
        	 \left( \exp(q_{1,2}\right.& e_1)\exp(q_{1,3}e_2)\cdots \left.
        	 	\exp(q_{1,n}e_{n-1}^{\eps(n-1)})\exp(p_{1,n}e_{n-1}^{\eps(n)})\cdots
        	 	\exp(p_{1,3}e_{2})\exp(p_{1,2}e_1)\right).	\notag
        	\end{align*}
        	Set in addition $s_{j,j}=t_{j,j}=p_{j,j}=q_{j,j}=u_{j,j}=1$, and for all
        	  $i,j$, $1\leq i< j\leq n$ put
        	\begin{equation*}
        	u_{i,j}=\left\{\begin{array}{clll}t_{i,j}+s_{i,j},& j< n,&&\\
        		t_{i,j},& j=n, & n-i \ \text{odd},&\ (\text{i.e.} \ \eps(n-i)=-1),\\
        		s_{i,j},& j=n, & n-i\  \text{even},&\ (\text{i.e.} \ \eps(n-i)=1)	
        	\end{array}\qquad r_{i,j}=\frac{u_{i,j-1}}{u_{i,j}s_{i,j-1}}.
        	\right.
        	\end{equation*}
        	\begin{theorem} a)\label{theorem3} BZ map for $SO(n,n)$ looks as
        		\begin{equation}\label{d5}\begin{split}
        	& p_{i,j}=r_{i,j}\cdot\prod_{k<i}\frac{t_{k,j-1}}{s_{k,j-1}}\prod_{k<i}\frac{s_{k,j}}{t_{k,j}},\qquad
        	 \qquad    	
        	q_{i,j}=r_{i,j}\cdot\prod_{k<i}\frac{t_{k,j-1}}{s_{k,j-1}}\prod_{k\leq i}\frac{s_{k,j}}{t_{k,j}},\qquad j<n\\ &
        	p_{i,n}=	r_{i,n}\cdot\prod_{k<i}\frac{t_{k,n-1}}{s_{k,n-1}}
        	 \prod_{k<i}\left(\frac{t_{k,n}}{s_{k,n}}\right)^{\eps(n-i)}\!\!\!\!\!\!\!\!\!\!\!\!\!,\qquad
        	  q_{i,n}= r_{i,n}\cdot\prod_{k<i}\frac{t_{k,n-1}}{s_{k,n-1}}
        	          	\prod_{k\leq i}\left(\frac{s_{k,n}}{t_{k,n}}\right)^{\eps(n-i)}\!\!\!\!\!\!\!\!\!\!\!\!.      	        
        	\end{split}	
        	\end{equation}				
   b) The twisting matrix $T\in H$ is $T=\sum_{k=1}^n\left(T_{k}e_{k,k}+T_k^{-1}e_{\hat{k},\hat{k}}\right)$, where
   \beqq
    T_k=\prod_{j<k}\frac{s_{j,k}}{t_{j,k}}\cdot\prod_{j>k}s_{k,j}t_{k,j}
   \eeqq
   c) BZ transform preserves the measure $\dfrac{d\t}{\t}$:
   $$\prod_{i<j}\frac{dt_{i,j}ds_{i,j}}{t_{i,j}s_{i,j}}= \prod_{i<j}\frac{dp_{i,j}dq_{i,j}}{p_{i,j}q_{i,j}} .$$				
	\end{theorem}
  We supply theorem \ref{theorem3} with remarks identical to those related to Theorem \ref{theorem1}. Their proofs are similar
  without any troubles with signs since here $\left(\ww{n}\right)^2=1$. Namely

\noindent {\bf Remark 1}. BZ transform $\tau$ is involutive, that is the inverse relations
\begin{equation}\label{d8}\begin{split}
& t_{i,j}=\tr_{i,j}\cdot\prod_{k<i}\frac{p_{k,j-1}}{q_{k,j-1}}\prod_{k<i}\frac{q_{k,j}}{p_{k,j}},\qquad
\qquad    	
s_{i,j}=\tr_{i,j}\cdot\prod_{k<i}\frac{p_{k,j-1}}{q_{k,j-1}}\prod_{k\leq i}\frac{q_{k,j}}{p_{k,j}},\qquad j<n\\ &
t_{i,n}=	\tr_{i,n}\cdot\prod_{k<i}\frac{p_{k,n-1}}{q_{k,n-1}}
\prod_{k<i}\left(\frac{p_{k,n}}{q_{k,n}}\right)^{\eps(n-i)}\!\!\!\!\!\!\!\!\!\!\!\!\!,\qquad
s_{i,n}= \tr_{i,n}\cdot\prod_{k<i}\frac{p_{k,n-1}}{q_{k,n-1}}
\prod_{k\leq i}\left(\frac{q_{k,n}}{p_{k,n}}\right)^{\eps(n-i)}\!\!\!\!\!\!\!\!\!\!\!\!.      	        	        	
\end{split}	
\end{equation}	
have the same form as \rf{d5}. Here
\begin{equation*}
v_{i,j}=\left\{\begin{array}{clll}p_{i,j}+q_{i,j},& j< n,&&\\
p_{i,j},& j=n, & n-i \ \text{odd},&\ (\text{i.e.} \ \eps(n-i)=-1),\\
q_{i,j},& j=n, & n-i\  \text{even},&\ (\text{i.e.} \ \eps(n-i)=1)	
\end{array}\qquad \tr_{i,j}=\frac{v_{i,j-1}}{v_{i,j}q_{i,j-1}}.
\right.
\end{equation*}
\medskip

\noindent {\bf Remark 2}. The Cartan twist $T$ can be as well  written  as
$T(\t)=T(\p^{-1})$, that is $T=\sum_{k=1}^n\left(T_{k}e_{k,k}+T_k^{-1}e_{\hat{k},\hat{k}}\right)$, where
\beqq
 T_k=\prod_{j<k}\frac{q_{j,k}}{p_{j,k}}\cdot\prod_{j>k}\frac{1}{q_{k,j}p_{k,j}}
\eeqq

\begin{corollary}\label{corollary3}
	The restriction of the left Whittaker vector to $N_+$ is given by the function
	\begin{equation*}
	 \wwl_{\nu}(\t)=\t^{\nu}\cdot\exp\Big(-\sum_{i<j}(p_{i,j}+q_{i,j})\Big)=\p^{-\nu}\cdot\exp\Big(-\sum_{i<j}(p_{i,j}+q_{i,j})\Big)
	\end{equation*}
	where $p_{i,j}$ are given by \rf{gl5}  and $\t^{\nu}$ for a weight $\nu=\sum_{k=1}^n \nu_k\ve_k$ means the product
	\beqq
	\t^{\nu}=\prod_{i<j} t_{i,j}^{\nu_i-\nu_j}s_{i,j}^{\nu_i+\nu_j}
	\eeqq
\end{corollary}
\bigskip

The proof of Theorem \ref{theorem3} follows the same scheme as for $\GL(n)$.
The matrix
$\Xn{n}(-\t)$ admits a factorization
$\Xn{n}(-\t)=\tXn{n-1}(-\t)\cdot \An{n}(-\t)$,
where
\beqq
\begin{split}
\An{n}(-\t)=
 & \exp(-s_{1,2} e_1)\exp(-s_{1,3}e_2)\cdots
\exp(-s_{1,n}e_{n-1}^{\eps(n-1)})\cdot\\ &\exp(-t_{1,n}e_{n-1}^{\eps(n)})\cdots
\exp(-t_{1,3}e_{2})\exp(-t_{1,2}e_1)
\end{split}\eeqq
and
 \begin{align*}
\tXn{n-1}(-\t)=&\left(\exp(-s_{n-1,n}e_{n-1}^{-})\exp-(t_{n-1,n}e_{n-1}^{+})\right)\cdot \notag\\
&\qquad\qquad\ldots\\
\left(- \exp(s_{2,3}\right.& \left. e_2)\cdots
\exp(-s_{2,n}e_{n-1}^{\eps(n-2)})\exp(-t_{2,n}e_{n-2}^{\eps(n-1)})\cdots
\exp(-t_{2,4}e_{3})\exp(-t_{2,3}e_2)\right)\cdot	\notag
\end{align*}
represents a group element of the unipotent subgroup of embedded $\SO(n-1,n-1)$,
$$  \tXn{n-1}(-\t)=\left(\begin{array}{c|c|c}1&0&0\\ \hline 0&\begin{array}{ccc}1&*&*\\0&\ddots&*\\0& 0& 1
\end{array}&0\\ \hline 0&0&1\end{array}
\right)$$
The matrix $\Xn{n}(\p)$ has the same structure,
$
\Xn{n}(\p)=\tXn{n-1}(\p)\cdot \An{n}(\p)$.
Denote by $\Ln{n}$ the diagonal matrix $\Ln{n}=\sum_{k=1}^n(\Ln{n}_{k,k}e_{k,k}+{\Ln{n}_{k,k}}^{-1}e_{\hat{k},\hat{k}})$, where
\beq\label{d15}\Ln{n}_{1,1}=\prod_{k=1}^{n-1}t_{1,k}s_{1,k}\qquad \Ln{n}_{k,k}=\frac{t_{1,k}}{s_{1,k}},\qquad 1<k<n,\qquad \Ln{n}_{n,n}=
 \left(\dfrac{t_{1,n}}{s_{1,n}}\right)^{\eps(n-1)}	
\eeq
and define the variables $\tilde{t}_{i,j}$ and $\tilde{s}_{i,j}$, $1< i<j\leq n$ by the relations
\beq\label{d16}\begin{split}
&\tilde{t}_{i,j}=t_{i,j}\cdot\frac{s_{1,j-1}}{t_{1,j-1}}\frac{t_{1,j}}{s_{1,j}},\qquad j<n,\qquad \tilde{t}_{i,n}=t_{i,n}\cdot\frac{s_{1,n-1}}{t_{1,n-1}}\left(\frac{t_{1,n}}{s_{1,n}}\right)^{\eps(n-1)}\\
&\tilde{s}_{i,j}=s_{i,j}\cdot\frac{s_{1,j-1}}{t_{1,j-1}}\frac{t_{1,j}}{s_{1,j}},\qquad j<n,\qquad \tilde{s}_{i,n}=s_{i,n}\cdot\frac{s_{1,n-1}}{t_{1,n-1}}\left(\frac{t_{1,n}}{s_{1,n}}\right)^{\eps(n)}
\end{split}
\eeq
\begin{proposition}\label{prop5} {\em{a)}} The parameters $p_{1,j}$ and $q_{1,j}$, $j=2,...,n$ of $\An{n}(\p)$ are equal to
	\beq\label{d17} p_{1,j}=\frac{u_{1,j-1}}{u_{1,j}s_{1,j-1}}, \qquad
	q_{1,j}=\frac{u_{1,j-1}s_{1,j}}{u_{1,j}s_{1,j-1}t_{1,j}}
\eeq
	{\em{b)}} Cartan twist $T$ and the matrix $\tXn{n-1}(\p)$ satisfy the relation
	\beqq
	T \tXn{n-1}(\p)=\Ln{n}\cdot\left(\tXn{n-1}(-\tilde{\t})\ww{n-1}\right)_{0+}
	\eeqq
\end{proposition}
{\bf Proof} of part a) of Proposition \ref{prop5} consists as before in comparison of the first rows of matrices $B^{(n)}=\An{n}(-\t)\ww{n}$
and $\An{n}(\p)$. We have
 \begin{equation*}\begin{split}&B^{(n)}_{1,j}=(-1)^{j+1}s_{1,i+1}\cdots s_{1,j-1}u_{1,j}, \qquad B^{(n)}_{1\hat{j}}=-(s_{1,2}\cdots s_{1,i+1})(t_{1,j+1}\cdots t_{1,n}),\qquad j<n,\\
 &B^{(n)}_{1,n}=s_{1,2}\cdots s_{1,n},\qquad B^{(n)}_{1\hat{n}}=s_{1,2}\cdots s_{1,n-1}t_{1,n}
 \end{split}
 \end{equation*}
By this we see first that the first diagonal entry $\Ln{n}_{1,1} $ in the Gauss decomposition of the right hand side of \rf{gl5} equals to $B^{(n)}_{1,1}=\prod_{j>1}t_{1,j}s_{1,j}$ and since the first row of the left hand side of
\rf{gl5}coincides with the first row of $\An{n}(\p)$,
$$\An{n}_{1,j}(\p)= q_{1,i+1}\cdots q_{1,j-1}v_{1,j}, \qquad\qquad\qquad \An{n}_{1\hat{j}}(\p)=(-1)^j(q_{1,2}\cdots q_{1,i+1})(p_{1,j+1}\cdots p_{1,n}),$$
 the variables $p_{1,j}$ and $q_{1,j}$ then can be found via ratios of coefficients $B^{(n)}_{1,k}$. Thus we get \rf{d17}. \hfill{$\square$}

 For the proof of part b) we again need the crucial technical lemma which says that  the matrix
 $U=	\left(\ww{n-1}\right)^{-1}\An{n}(-\t)\ww{n}\left(\An{n}(\p)\right)^{-1}$, see \rf{gl12} has the block structure
 \beqq
 U=\left(\begin{array}{c|c|c}*&0&0\\ \hline *&\begin{array}{ccc}*&0&0\\0&\ddots&0\\0& 0& *
 	\end{array}&0\\ \hline *&*&*\end{array}
 \right).
 \eeqq
  and specializes its diagonal entries.
 \begin{lemma}\label{lemma4} Nondiagonal matrix elements $U_{i,j}$  ($i\not=j$) equal zero if $j>1$ and $i<2n$.
 	Matrix element $U_{1,1}$ equals $\prod_{j>1} t_{1,j}s_{1,j}$. Matrix element $U_{k,k}$ equals $\dfrac{s_{1,k}}{t_{1,k}}$ for $1<k\leq n$. Matrix element $U_{\hat{k},\hat{k}}$ equals $\dfrac{t_{1,k}}{s_{1,k}}$ for $1<k\leq n$	
 	\end{lemma}
 The proof of Lemma \ref{lemma4} is sketched in Appendix B.
\bigskip

As well as in the proof of Proposition \rf{prop4}, Lemma \ref{lemma4} implies the equality
\beqq
T \tXn{n-1}(\p)=\Ln{n}\left(\ttXn{n-1}(-{\t})\ww{n-1}\right)_{0+}
\eeqq
where the diagonal matrix $\Ln{n}$ is given in \rf{d15} and
\beq\label{d19}\ttXn{n-1}(-\t)=\left(\Ln{n}\right)^{-1}\tXn{n-1}(-\t)\Ln{n}.\eeq
Then the structure of the group element $\Xn{n}(\t)$, see \rf{d3} says that the parameters
 $t_{i,j}$ and $s_{i,j}$ are the coefficients at $e_{j-1}$ for $j<n$; $t_{i,n}$ is the coefficient at $e_{n-1}^{\eps(n-i+1)}$, and $s_{i,n}$ is the coefficient at $e_{n-1}^{\eps(n-i)}$ in Lusztig presentation of $X(-\t)$. This enables us to rewrite the conjugation \rf{d19} as the change of variables \rf{d16} and
finish the proof of Proposition \ref{prop5}. Then the proof of part b) of Theorem \ref{theorem3} follows by induction on $n$. The inductive proof of part c) is analogous to that of Theorem \ref{theorem1}.
\hfill{$\square$}
 \subsection{Whittaker function} \label{section3.2}
Using the arguments of the end of the previous subsection, we  describe the action of Cartan
generators on the right Whittaker vector. Namely, the vector
 $\exp(\sum_{k=1}^n  x_k (e_{\hat{k},\hat{k}}-e_{k,k}))\cdot\wr_\mu$
is presented by the  function  $\exp(-\sum_k \mu_k x_k)\wwr_\mu(\t')$, where
\beq\label{d20} \begin{split}&t'_{i,j}=t_{i,j}\cdot\exp\left(x_{j}-x_{j-1} \right),\qquad\qquad\qquad
 s'_{i,j}=s_{i,j}\cdot\exp\left(x_{j}-x_{j-1} \right),\qquad j<n,\\
&t'_{i,n}=t_{i,n}\cdot\exp\left(\eps(n-i+1)x_{n}-x_{n-1} \right),\qquad
s'_{i,n}=s_{i,n}\cdot\exp\left(\eps(n-i)x_{n}-x_{n-1} \right).
\end{split}\eeq
Denote by $\wwlm_{\nu}(\bm{\ga})$ and $\wwrm_\mu(\bm{\ga})$ the Mellin transforms of the functions on $N_+$
$\wwl_{\nu}(\t)$ and $\wwr_\mu(\t)$,
\begin{equation*}\begin{split}& \wwrm_\mu(\bm{\ga})=\int_{\t>0}\wwr_\mu(\t) \t^{\mathbf{\ga}}\frac{d\t}{\t}:=
\int\limits_{\substack{t_{i,j}>0,\\ s_{i,j}>0}} \wwr_\mu(\t)\prod_{i<j}t_{i,j}^{\ga_{i,j}}s_{i,j}^{\delta_{i,j}}\frac{dt_{i,j}}{t_{i,j}}\frac{ds_{i,j}}{s_{i,j}},\\
& \wwlm_{\nu}(\bm{\ga})=\int_{\t>0}\t^\rho\wwl_{\nu}(\t) \t^{\mathbf{\ga}}\frac{d\t}{\t}:=
\int\limits_{\stackrel{t_{i,j}>0,}{ s_{i,j}>0}} \wwl_{\nu}(\t)\prod_{i<j}t_{i,j}^{\ga_{i,j}+j-i}s_{i,j}^{\delta_{i,j}+2n-i-j}\frac{dt_{i,j}}{t_{i,j}}.
\frac{ds_{i,j}}{s_{i,j}}\end{split}
\end{equation*}
Here $\rho=\sum_{k=1}^n (n-k)(\ve_{k,k}-\ve_{\hat{k},\hat{k}})$.
Due to \rf{d20} the action of Cartan subgroup on the right Whittaker vector $\wr_\mu$ transforms
the function $\wwrm_\mu(\bm{\ga})$ to the product
$
\exp(-\sum_k \mu_k x_k)\cdot\exp \tilde{H}(\bm{x},\bm{\gamma})\wwrm_\mu(\bm{\ga})
$,
 where
\beq\label{d20a}\begin{split}
\tilde{H}(\bm{x},\bm{\gamma})=\!\sum_{i<j<n}(\ga_{i,j}+\delta_{i,j})(x_j-& x_{j-1})+(\ga_{n-1,n}+\delta_{n-2,n}+\ldots)    (x_n-x_{n-1}) +\\ &(\delta_{n-1,n}+\gamma_{n-2,n}+\delta_{n-3,n}+\ldots)(-x_{n-1}-x_n),\end{split}\eeq

 Due to Proposition \ref{prop2} we have:
\beqq
\wwrm_\mu(\bm{\gamma})=\prod_{i<j}\Gamma(\gamma_{i,j})\Gamma(\delta_{i,j}).
\eeqq
For the calculation of $\wwlm_{\nu}(\bm{\gamma})$ we pass in the integral
$$\int_{t_{i,j},s_{i,j}>0}\exp \sum_{i<j}-(p_{i,j}+q_{i,j})\cdot\prod_{i<j}t_{i,j}^{\ga_{i,j}+\nu'_i-\nu'_j}
s_{i,j}^{\delta_{i,j}+\nu'_i+\nu'_j}\frac{dt_{i,j}}{t_{i,j}}\frac{ds_{i,j}}{s_{i,j}}$$
from the integration variables $t_{i,j}$ and $s_{i,j}$ to $p_{i,j}$ and $q_{i,j}$. Here
$$\nu'=\nu+\rho,\qquad \nu'_i=\nu_i+n-i.$$
Substitution of \rf{d8} gives the relation
$$\prod_{i<j}t_{i,j}^{\ga_{i,j}}
s_{i,j}^{\delta_{i,j}}=\prod_{1\leq i<j\leq n}p_{i,j}^{\phi_{i,j}}q_{i,j}^{\psi_{i,j}}\cdot\prod_{1\leq i<j<n}\left(p_{i,j}+q_{i,j}\right)^{\theta_{i,j}}$$
where
\begin{equation}\label{D15}\begin{split}
\varphi_{i,j}&=-\delta_{i,j}+\sum_{k>i}\left(\gamma_{k,j+1}+\delta_{k,j+1}- \gamma_{k,j}-\delta_{k,j}\right)\\
\psi_{i,j}&=\ \delta_{i,j}-\sum_{k>i}\left(\gamma_{k,j+1}+\delta_{k,j+1}- \gamma_{k,j}-\delta_{k,j}\right)-\gamma_{i,j+1}-\delta_{i,j+1},\\
\theta_{i,j}&=-\delta_{i,j}-\gamma_{i,j}+\delta_{i,j+1}+\gamma_{i,j+1},
\end{split}\end{equation}
if $j<n$ and $\theta_{i,n}=0$,
\begin{equation} \label{D15a} \begin{split}
\varphi_{i,n}&=\!\left\{\!\!\begin{array}{ll}-\delta_{i,n}+\sum_{k=1}^{n-i-1}(-1)^{k}(\gamma_{i+k,n}-\delta_{i+k,n}),& \ve(n-i)=1 \\
-\gamma_{i,n}+\sum_{k=1}^{n-i-1}(-1)^{k}(\delta_{i+k,n}-\gamma_{i+k,n}),& \ve(n-i)=-1 \end{array}\right.\\
\psi_{i,n}&=\!\left\{\!\!\begin{array}{ll}-\gamma_{i,n}+\sum_{k=1}^{n-i-1}(-1)^{k}(\delta_{i+k,n}-\gamma_{i+k,n}),& \ve(n-i)=1 \\
-\delta_{i,n}+\sum_{k=1}^{n-i-1}(-1)^{k}(\gamma_{i+k,n}-\delta_{i+k,n}),& \ve(n-i)=-1
\end{array}\right. \end{split}
\end{equation}
This means that the Mellin transform of the left Whittaker vector is described by the integral
$$
\int_{p_{i,j},q_{i,j}>0}\prod_{i<j} p_{i,j}^{{\varphi}_{i,j}-\nu'_i+\nu'_j}q_{i,j}^{{\psi_{i,j}-\nu'_i-\nu'_j}}(p_{i,j}+q_{i,j})^{\theta_{i,j}}\frac{dp_{i,j}}{p_{i,j}}\frac{dq_{i,j}}{q_{i,j}}
$$
Integration over the variables $p_{i,n}$ and $q_{i,n}$ produced the product
\beq\label{D15b}\prod_{i=1}^{n-1}\Gamma(\vf_{i,n}-\nu_i+\nu_n)\Gamma(\psi_{i,n}-\nu'_i-\nu'_n)\eeq
 of corresponding $\Gamma$ functions.
For the calculation of integrals over the variables $p_{i,j}$ and $q_{i,j}$ with $j<n$ we use the integral
\beq\label{int}
\int_{x,y>0} x^{a-1} y^{b-1}(x+y)^ce^{-x-y}dxdy=\frac{\Gamma(a)\Gamma(b)\Gamma(a+b+c)}{\Gamma(a+b)}.
\eeq
Its evaluation is based on the change of variables $u=x+y,\ t=\dfrac{x}{x+y}$.
This gives the product over $1<i<j<n$ of the factors
\beq\label{D16} \prod_{1\leq i<j\leq n}\Gamma (\varphi_{i,j}-\nu'_i+\nu'_j)\Gamma(\psi_{i,j}-\nu'_i-\nu'_j)\frac{\Gamma(\varphi_{i,j}+\theta_{i,j}+\psi_{i,j}-2
	\nu_i)}{\Gamma(\varphi_{i,j}+\psi_{i,j}-2\nu'_i)}.\eeq
Due to \rf{D15}
 $$\frac{\Gamma(\varphi_{i,j}+\theta_{i,j}+\psi_{i,j}-2
 	\nu_i)}{\Gamma(\varphi_{i,j}+\psi_{i,j}-2\nu'_i)}=
\frac{\Gamma(-\gamma_{i,j}-\delta_{i,j}-2\nu'_i)}{\Gamma(-\gamma_{i,j+1}-\delta_{i,j+1}-2\nu'_i)}
$$
which implies cancelations of ratios in the product \rf{D16} so that it becomes equal to
\begin{equation}\label{D17}
\prod_{1\leq i<j< n}\Gamma (\varphi_{i,j}-\nu'_i+\nu'_j)\Gamma(\psi_{i,j}-\nu'_i-\nu'_j)
\prod_{i=1}^{n-2}\frac{\Gamma(-\gamma_{i,i+1}-\delta_{i,i+1}-2\nu'_i)}{\Gamma(-\gamma_{i,n}-\delta_{i,n}-2\nu'_i)}
 .
\end{equation}
Using \rf{D15b} we arrive to the following answers
\begin{equation}\label{D19}
\wwlm_{\nu}(\bm{\gamma})=\prod_{i=1}^{n-2}\ \frac{\Gamma(-\gamma_{i,i+1}-\delta_{i,i+1}-2\nu'_i)}{\Gamma(-\gamma_{i,n}-\delta_{i,n}-2\nu'_i)}\cdot \!\!\!\!
\prod_{1\leq i<j\leq n}\Gamma({\varphi}_{i,j}-\nu'_i+\nu'_j)\Gamma({\psi}_{i,j}-\nu'_i-\nu'_j).
\end{equation}	
and
	\beq\label{D20}\begin{split}
		\Psi_{\lambda}(\x)= \frac{e^{-i(\bm{\lambda},\x)}}{(2\pi i)^d}\int\limits_{C}
		\exp {\tilde{H}}(\x,\bm{\gamma})
	 \prod_{k=1}^{n-2}\frac{\Gamma(\gamma_{k,k+1}+\delta_{k,k+1}+2i\lambda_k)}{\Gamma(\gamma_{k,n}
	 	+\delta_{k,n}+
		2i\lambda_k)}\\
	\prod_{k<l} \Gamma(-{\varphi}_{k,l}+i(\lambda_k-\lambda_l))\Gamma(-{\psi}_{k,l}+i(\lambda_k+\lambda_l))
	\Gamma(\gamma_{k,l})\Gamma(\delta_{k,l})d\gamma_{k,l}d\delta_{k,l}
\end{split}\eeq
 Here $\tilde{H}(\bm{x},\bm{\gamma})$ is given in \rf{d20a}, $\psi_{i,j}$ and $\vf_{i,j}$ are given in \rf{D15} and \rf{D15a}, $d=n(n-1)$ and the contour $C$ is a deformation of the imaginary plane $\Re \gamma_{ij}=\Re \delta_{ij}=0$
 into the strip of analyticity of the integrand.

Perform now the following  change of variables:
	\beqq
	\begin{split} &\gamma_{i,j}=\tgamma_{i,j}-\tgamma_{i+1,j},\qquad
		\delta_{i,j}=\tdelta_{i,j}-\tdelta_{i+1,j},\qquad 1\leq i<j < n,\\
		&\gamma_{i,n}=\left\{\begin{array}{ll}
			\tdelta_{i,n}-\tdelta_{i+1,n},& \ve(n-i)=1,\\
			\tgamma_{i,n}-\tgamma_{i+1,n},& \ve(n-i)=-1,
			\end{array}\right.	\\
		&\delta_{i,n}=\left\{\begin{array}{ll}
			\tgamma_{i,n}-\tgamma_{i+1,n},& \ve(n-i)=1,\\
			\tdelta_{i,n}-\tdelta_{i+1,n},& \ve(n-i)=-1,
		\end{array}\right.
	\end{split}\eeqq
	that is
	\beqq
	 \begin{split}
		&\tgamma_{i,j}=\gamma_{i,j}+\ldots+\gamma_{j-1,j},\qquad \tdelta_{i,j}=\delta_{i,j}+\ldots+\delta_{j-1,j},
		\qquad
		1\leq i<j<n,\\
		&\tgamma_{i,n}=\delta_{in}+\gamma_{i+1,n}+...+\gamma_{n-1,n}\qquad\text{if}\ n-i\ \text{is even},\\
		&\tgamma_{i,n}=\gamma_{in}+\delta_{i+1,n}+...+\gamma_{n-1,n}\qquad\text{if}\ n-i\ \text{is odd},\\
		&\tdelta_{i,n}=\gamma_{in}+\delta_{i+1,n}+...+\delta_{n-1,n}\qquad\text{if}\ n-i\ \text{is even},\\
		&\tdelta_{i,n}=\delta_{in}+\gamma_{i+1,n}+...+\delta_{n-1,n}\qquad\text{if}\ n-i\ \text{is odd},\\	
	\end{split}
	\eeqq	
	In this variables we have
	\beqq
	\begin{split}
		&\vf_{i,j}=\tgamma_{i+1,j+1}+\tdelta_{i+1,j+1}-\tgamma_{i+1,j}-\tdelta_{i,j},\qquad 1\leq i<j<n,\\	
		&\psi_{i,j}=-\tgamma_{i,j+1}-\tdelta_{i,j+1}+\tgamma_{i+1,j}+\tdelta_{i,j},\qquad 1\leq i<j<n,\\
		&\vf_{i,n}=-\tgamma_{i,n}+\tdelta_{i+1,n},\qquad
		\psi_{i,n}=-\tdelta_{i,n}+\tgamma_{i+1,n}\qquad 1\leq i<n,	\end{split}
	\eeqq
and
	$H(\bm{x},\bm{\gamma})=\!\sum_{j<n}(\tgamma_{1,j}+\tdelta_{1,j})(x_j- x_{j-1})+\tgamma_{1,n}  (x_n-x_{n-1}) +\tdelta_{1,n}(-x_{n-1}-x_n)$.
Then we have
 \begin{theorem}\label{t2} The function $\Psinew_{\blambda}(\x)$ is given by the integral
	\beqq
	\begin{split}
		\Psinew_{\blambda}(\x)&= \frac{e^{  -i(\x,\blambda)}}{(2\pi i)^\d}\int\limits_{C_n}	
		\exp {H}(\bm{x},\bm{\gamma})
		\prod_{k=1}^{n-2}\frac{\Gamma(\gamma_{k,k+1}+\delta_{k,k+1}+2i\lambda_k)}{\Gamma(\gamma_{k,n}
			\!+\delta_{k,n}\!-\gamma_{k+1,n}\!-\delta_{k+1,n}\!+
			2i\lambda_k)} \\
		&\!\!\!\prod_{1\leq k<l\leq n} \!\!\!\Gamma({\xi}_{k,l}+i(\lambda_k-\!\lambda_l))\Gamma({\eta}_{k,l}+i(\lambda_k+\!\lambda_l))
		\Gamma(\gamma_{k,l}-\gamma_{k+1,l})\Gamma(\delta_{k,l}-\delta_{k+1,l})d\gamma_{k,l}d\delta_{k,l}.
	\end{split}\eeqq
\end{theorem}
Here we set $\gamma_{k,l}=\delta_{k,l}=0$ if the pair $(k,l)$ does not satisfies the condition $1\leq k<l\leq n$,
$\d=n(n-1)$,
$H(\bm{x},\bm{\gamma})=\!\sum_{j<n}(\gamma_{1,j}+\delta_{1,j})(x_j- x_{j-1})+\gamma_{1,n}  (x_n-x_{n-1}) -\delta_{1,n}(x_{n-1}+x_n),$ 
\beqq
\begin{split}
	&\xi_{i,j}=-\gamma_{i+1,j+1}-\delta_{i+1,j+1}+\gamma_{i+1,j}+\delta_{i,j},\qquad 1\leq i<j<n,\\	
	&\eta_{i,j}=\gamma_{i,j+1}+\delta_{i,j+1}-\gamma_{i+1,j}-\delta_{i,j},\qquad 1\leq i<j<n,\\
	&\xi_{i,n}=\gamma_{i,n}-\delta_{i+1,n},\qquad 1\leq i<n,\\
	&\eta_{i,n}=\delta_{i,n}-\gamma_{i+1,n}\qquad 1\leq i<n,	\end{split}
\eeqq
The integration cycle is a deformation of the imaginary plain
$\Re \gamma_{k,l}=\Re \delta_{k,l}=0$ into nonempty domain  $D\subset\C^{\d}$ of the analyticity of the integrand, which is described by the relations
$
\Re \gamma_{i,j}>0$, \ $\Re \delta_{i,j}>0$,\ $Re\, \xi_{i,j}>0$,\ $\Re \eta_{i,j}>0$,
\  $ 1\leq i<j\leq n$.


\setcounter{equation}{0}
\section{$\SO(n+1,n)$}
\subsection{BZ transform}
The split real form of the group $\SO(2n+1,\C)$ is the group $\SO(n+1,n)\subset \SL(2n,\R)$ preserving symmetric form $\sum_{i=1}^{n} (x_iy_{\hat{i}}+x_{\hat{i}}y_i)+x_ny_n$, where $\hat{i}=2n+2-i$. Gauss decomposition is induced from that of $GL(2n+1,\R)$. Positive roots are elements $\ve_i\pm\ve_j$ for $1\leq i<j\leq n$ and $\ve_i$, where
$\ve_k\in\h^*$, $k=1,...,n$ are  defined by the condition $\ve_k( e_{j,j}-e_{\hat{j}\hat{j}})=\delta_{j,k}$.
We denote Chevalley generator of the Lie algebra $\g=\so(2n+1)$ by $e_i$ and $f_i$, $i=1,\ldots,n$  Here
\begin{align*} &e_i=e_{i,i+1}-e_{2n-i,2n+1-i},&&  f_i=e_{i+1,i}-e_{2n+1-i,2n-i},&&
i=1,...,n-1,\\
& e_{n}=\sqrt{2}(e_{n,n+1}-e_{n+1,n+2}),&& f_{n}=\sqrt{2}(e_{n+1,n}-e_{n+2,n+1}).
\end{align*}
Denote by $s_i$ $i=1,...,n$  the corresponding generators of the Weyl group and by $\ss_i$
their lifts to the group $\SO(n+1,n)$ according to \rf{W2a}.
We choose the following normal ordering of the system $\Delta_+$ of positive roots:
\beqq
\begin{split}&\ve_n,\ (\ve_{n-1}+\ve_n,\, \ve_{n-1},\ve_{n-1}-\ve_n),\ \cdots, \\
	&(\ve_{1}+\ve_{2},\ve_{1}+\ve_{3},\cdots,\ve_{1}+\ve_{n},\ve_1,\ve_{1}-\ve_{n},\cdots, \ve_{1}-\ve_{3},\ve_{1}-\ve_{2}).
\end{split}
\eeqq
It corresponds by \rf{L2} to the following reduced decomposition of the longest element $w_0^{(n)}$ of the Weyl group $W$:
\beq\label{b2}w_0^{(n)}=s_n({s}_{n-1}s_ns_{n-1})\,(s_{n-2}s_{n-1}s_n{s}_{n-1}s_{n-2})\,
\cdots\,(s_{1}s_2...s_{n-1}{s}_{n}s_{n-1}...s_2s_1).
\eeq
 Note that as element of $\SO(n+1,n)$ ,
\beqq
\bar{w}^{(n)}_0=(-1)^n\left(\sum_{k=1}^{2,n}e_{k\hat{k}}+e_{n+1,n+1}\right)
\eeqq
 Denote by $t_{i,j}$ Lusztig parameter corresponding to the root $\ve_i-\ve_j$, and  by $s_{i,j}$ Lusztig parameter corresponding to the root $\ve_i+\ve_j$, and by $t_k$ the parameter, corresponding to $\ve_k$. Here $1\leq i<j\leq n$, $k=1,...,n$.  Then the group element $\Xn{n}(\t)$ looks as
\begin{align}\label{b3}
\Xn{n}(\t)=&\exp( t_ne_n)\left(\exp(s_{n-1,n}e_{n-1})\exp(t_{n-1}e_n)\exp(t_{n-1,n}e_{n-1})\right)\cdot \notag\\
&\qquad\qquad\ldots\\
\left( \exp(s_{1,2}\right.& e_1)\exp(s_{1,3}e_2)\cdots \left.
\exp(s_{1,n}e_{n-1})\exp t_1e_n\exp(t_{1,n}e_{n-1})\cdots
\exp(t_{1,3}e_{2})\exp(t_{1,2}e_1)\right).	\notag
\end{align}
BZ transform for $\SO(n+1,n)$ is the solution of the equation \rf{gl5}, where in $\Xn{n}(\p)$ we use the notation
$p_{i,j}$ for Lusztig parameter, corresponding to the root $\ve_i-\ve_j$,  $q_{i,j}$ for the Lusztig parameter, corresponding to the root $\ve_i+\ve_j$, and $p_i$ for the parameters, corresponding to $\ve_i$,	
\begin{align}
\Xn{n}(\p)&=\exp( p_ne_n)\left(\exp(q_{n-1,n}e_{n-1})\exp(p_{n-1}e_n)\exp(p_{n-1,n}e_{n-1})\right)\cdot \notag\\
&\qquad\qquad\ldots\label{b4a}\\
 \exp(q_{1,2}& e_1)\exp(q_{1,3}e_2)\cdots
\exp(q_{1,n}e_{n-1})\exp p_1e_n\exp(t_{1,n}e_{n-1})\cdots
\exp(p_{1,3}e_{2})\exp(p_{1,2}e_1).	\notag
\end{align}
We set 	$s_{j,j}=t_{j,j}=p_{j,j}=q_{j,j}=p_j=u_{j,j}=1$, and for all
$i,j$, $1\leq i< j\leq n$ put
\begin{equation}\label{b5a}
u_{i,j}=t_{i,j}+s_{i,j}, \qquad r_{i,j}=\frac{u_{i,j-1}}{u_{i,j}s_{i,j-1}}.
\end{equation}
\begin{theorem} a)\label{theorem5} BZ map for $SO(n,n)$ looks as
	\begin{equation}\label{b5}
	 p_{i,j}=r_{i,j}\cdot\prod_{k<i}\frac{t_{k,j-1}}{s_{k,j-1}}\prod_{k<i}\frac{s_{k,j}}{t_{k,j}},
	\qquad    	
	q_{i,j}=r_{i,j}\cdot\prod_{k<i}\frac{t_{k,j-1}}{s_{k,j-1}}\prod_{k\leq i}\frac{s_{k,j}}{t_{k,j}},\qquad
	p_i=\frac{u_{i,n}}{t_is_{i,n}}\prod_{k<i}	\frac{t_{kn}}{s_{kn}}        		
	\end{equation}				
	b) The twisting matrix $T\in H$ is $T=e_{n+1,n+1}+\sum_{k=1}^n\left(T_{k}e_{k,k}+T_k^{-1}e_{\hat{k},\hat{k}}\right)$, where
	\beqq
	T_k=t_k^2\prod_{j<k}\frac{s_{j,k}}{t_{j,k}}\cdot\prod_{j>k}s_{k,j}t_{k,j}
	\eeqq
	c) BZ transform preserves the measure $\frac{d\t}{\t}$:
	$$\prod_{i<j}\frac{dt_{i,j}ds_{i,j}}{t_{i,j}s_{i,j}}\prod_k\frac{dt_k}{t_k}= \prod_{i<j}\frac{dp_{i,j}dq_{i,j}}{p_{i,j}q_{i,j}}\prod_k\frac{dp_k}{p_k} .$$				
\end{theorem}	
Again, we have the same remarks with the same proofs:

\noindent {\bf Remark 1}. BZ transform $\tau$ is involutive, that is the inverse relation
		\begin{equation}\label{b7}
	 t_{i,j}=\tilde{r}_{i,j}\cdot\prod_{k<i}\frac{p_{k,j-1}}{q_{k,j-1}}\prod_{k<i}\frac{q_{k,j}}{p_{k,j}},
	\qquad    	
	s_{i,j}=\tilde{r}_{i,j}\cdot\prod_{k<i}\frac{p_{k,j-1}}{q_{k,j-1}}\prod_{k\leq i}\frac{q_{k,j}}{p_{k,j}},\qquad
	t_i=\frac{v_{i,n}}{p_iq_{i,n}}\prod_{k<i}	\frac{p_{kn}}{q_{kn}}        		
	\end{equation}
has the same form as \rf{d5}. Here
\begin{equation*}
v_{i,j}=p_{i,j}+q_{i,j},\qquad \tr_{i,j}=\frac{v_{i,j-1}}{v_{i,j}q_{i,j-1}}.
\end{equation*}
\medskip
\noindent {\bf Remark 2}. The Cartan twist $T$ can be as well  written  as
$T(\t)=T(\p^{-1})$, that is $T=e_{n+1,n+1}+\sum_{k=1}^n\left(T_{k}e_{k,k}+T_k^{-1}e_{\hat{k},\hat{k}}\right)$, where
\beqq
T_k=\frac{1}{p_k^2}\prod_{j<k}\frac{q_{j,k}}{p_{j,k}}\cdot\prod_{j>k}\frac{1}{q_{k,j}p_{k,j}}
\eeqq

\begin{corollary}\label{corollary4}
	The restriction of left Whittaker vector to $N_+$ is given by the function
	\begin{equation*}
	 \wwl_{\nu}(\t)=\t^{\nu}\cdot\exp\Big(-\sum_{i<j}(p_{i,j}+q_{i,j})-\sum_kp_k\Big)=\p^{-\nu}\cdot\exp\Big(-\sum_{i<j}(p_{i,j}+q_{i,j})-\sum_kp_k)\Big)
	\end{equation*}
	where $p_{i,j}$  and $p_k$ are given by \rf{b5}  and $\t^{\nu}$ for a weight $\nu=\sum_{k=1}^n \nu_k\ve_k$ means the product
	\beqq
	\t^{\nu}=\prod_{i<j} t_{i,j}^{\nu_i-\nu_j}s_{i,j}^{\nu_i+\nu_j}\prod_kt_k^{2\nu_k}
	\eeqq
\end{corollary}
\bigskip

The proof of Theorem \ref{theorem5} follows the same scheme as before.
The matrix
$\Xn{n}(-\t)$ admits a factorization
$\Xn{n}(-\t)=\tXn{n-1}(-\t)\cdot \An{n}(-\t)$,
where
\begin{equation}\label{b13a}
\begin{split}
	\An{n}(-\t)=
	& \exp(-s_{1,2} e_1)\exp(-s_{1,3}e_2)\cdots
	\exp(-s_{1,n}e_{n-1})\exp (-t_1e_n)\cdot\\ &\exp(-t_{1,n}e_{n-1})\cdots
	\exp(-t_{1,3}e_{2})\exp(-t_{1,2}e_1)
\end{split}\end{equation}
and
\begin{align}
\tXn{n-1}(\t)=&\exp( t_ne_n)\left(\exp(s_{n-1,n}e_{n-1})\exp(t_{n-1}e_n)\exp(t_{n-1,n}e_{n-1})\right)\cdot \notag\\
&\qquad\qquad\ldots\label{b13b}\\
\left( \exp(s_{2,3}\right.& e_2)\cdots \left.
\exp(s_{2,n}e_{n-1})\exp t_2e_n\exp(t_{2,n}e_{n-1})\cdots
\exp(t_{2,3}e_2)\right).	\notag
\end{align}
represents a group element of the unipotent subgroup of embedded
$\SO(n,n-1)$.
The matrix $\Xn{n}(\p)$ has the same structure,
$
\Xn{n}(\p)=\tXn{n-1}(\p)\cdot \An{n}(\p)$.
Denote by $\Ln{n}$ the diagonal matrix $\Ln{n}=e_{n+1,n+1}+\sum_{k=1}^n(\Ln{n}_{k,k}e_{k,k}+{\Ln{n}_{k,k}}^{-1}e_{\hat{k},\hat{k}})$, where
\beqq
\Ln{n}_{1,1}=t_1^2\prod_{k=1}^{n-1}t_{1,k}s_{1,k}\qquad \Ln{n}_{k,k}=\frac{t_{1,k}}{s_{1,k}},\qquad 1<k\leq n,	
\eeqq
and define the variables $\tilde{t}_{i,j}$ and $\tilde{s}_{i,j}$, $1< i<j\leq n$ by the relations
\beqq
	\tilde{t}_{i,j}=t_{i,j}\cdot\frac{s_{1,j-1}}{t_{1,j-1}}\frac{t_{1,j}}{s_{1,j}}, \qquad
	\tilde{s}_{i,j}=s_{i,j}\cdot\frac{s_{1,j-1}}{t_{1,j-1}}\frac{t_{1,j}}{s_{1,j}},\qquad
	\tilde{t}_i=t_i\frac{s_{1,n}}{t_{1,n}}
\eeqq
The induction step is given by the following lemma and proposition.
\begin{lemma}\label{lemma5} Nondiagonal matrix elements $U_{i,j}$ ($i\not=j$) of the matrix \rf{gl12} equal zero if $j>1$ and $i<2n+1$.
	Matrix element $U_{1,1}$ equals $t_1^2\prod_{j>1} t_{1,j}s_{1,j}$. Matrix element $U_{k,k}$ equals $\dfrac{s_{1,k}}{t_{1,k}}$ for $1<k\leq n$. Matrix element $U_{\hat{k},\hat{k}}$ equals $\dfrac{t_{1,k}}{s_{1,k}}$ for $1<k\leq n$. Element $U_{n+1,n+1}=1$.	
\end{lemma}
\begin{proposition}\label{prop6} {\em{a)}} The parameters $p_{1,j}$, $q_{1,j}$, $j=2,...,n$ and $p_1$ of $\An{n}(\p)$ are equal to
	\beqq
	 p_{1,j}=\frac{u_{1,j-1}}{u_{1,j}s_{1,j-1}}, \qquad
	q_{1,j}=\frac{u_{1,j-1}s_{1,j}}{u_{1,j}s_{1,j-1}t_{1,j}},\qquad p_1=\frac{u_{1,n}}{t_1s_{1,n}}
	\eeqq
	{\em{b)}} Cartan twist $T$ and the matrix $\tXn{n-1}(\p)$ satisfy the relation
	\beqq
	T \tXn{n-1}(\p)=\Ln{n}\cdot\left(\tXn{n-1}(-\tilde{\t})\ww{n-1}\right)_{0+}
	\eeqq
\end{proposition}
All the proofs repeat that of the previous sections \hfill{$\square$}

\subsection{Whittaker function} \label{section4.2}
Describe first the action of Cartan
generators on the right Whittaker vector. The vector
$\exp(\sum_{k=1}^n  x_k (e_{\hat{k},\hat{k}}-e_{k,k}))\cdot\wr_\mu$
is presented by the  function  $\exp(-\sum_k \mu_k x_k)\wwr_\mu(\t')$, where
\begin{equation}\label{b16} t'_{i,j}=t_{i,j}\cdot\exp\left(x_{j}-x_{j-1} \right),\qquad
	s'_{i,j}=s_{i,j}\cdot\exp\left(x_{j}-x_{j-1} \right),\qquad t'_i=t_i\cdot\exp (-x_n).
\end{equation}
Denote by $\wwlm_{\nu}(\bm{\ga})$ and $\wwrm_\mu(\bm{\ga})$ the Mellin transforms of the functions
$\wwl_{\nu}(\t)$ and $\wwr_\mu(\t)$,
\begin{equation*}\begin{split}& \wwrm_\mu(\bm{\ga})=\int\limits_{\t>0}\wwr_\mu(\t) \t^{\mathbf{\ga}}\frac{d\t}{\t}:=
\!\!\int\limits_{\substack{t_{i,j}>0,\\ s_{i,j}>0,t_i>0}}\!\! \wwr_\mu(\t)\prod_{i<j}t_{i,j}^{\ga_{i,j}}s_{i,j}^{\delta_{i,j}}t_i^{\gamma_i}\frac{dt_{i,j}}{t_{i,j}}\frac{ds_{i,j}}{s_{i,j}}\frac{dt_i}{t_i},\\
& \wwlm_{\nu}(\bm{\ga})=\int\limits_{\t>0}\t^\rho\wwl_{\nu}(\t) \t^{\mathbf{\ga}}\frac{d\t}{\t}:=
\!\!\int\limits_{\substack{t_{i,j}>0,\\ s_{i,j}>0,t_i>0}}\!\! \wwl_{\nu}(\t)\prod_{i<j}t_{i,j}^{\ga_{i,j}+j-i}s_{i,j}^{\delta_{i,j}+2n+1-i-j}t_i^{\gamma_i+n-i+1/2}\frac{dt_{i,j}}{t_{i,j}}
\frac{ds_{i,j}}{s_{i,j}}\frac{dt_i}{t_i}.\end{split}
\end{equation*}
 Here $\rho=\sum_{k=1}^n (n+\frac12-k)(\ve_{k,k}-\ve_{\hat{k},\hat{k}})$.
Due to \rf{b16} the action of Cartan subgroup on the right Whittaker vector $\wr_\mu$ transforms the function $\wwrm_\mu(\bm{\ga})$ to the product
$
\exp(-\sum_k \mu_k x_k)\cdot\exp \tilde{H}(\bm{x},\bm{\gamma})\wwrm_\mu(\bm{\ga})
$,
where
\beq\label{b17}
\tilde{H}(\bm{x},\bm{\gamma})=\sum_{i<j}(\gamma_{i,j}+\delta_{i,j})(x_{j}-x_{j-1}) -x_n\sum_{i}\gamma_i.\eeq

As before, we have
\beq\label{b18}
\wwrm_\mu(\bm{\gamma})=\prod_{1<i<j\leq n}\Gamma(\gamma_{i,j})\Gamma(\delta_{i,j})\prod_{k=1}^n\Gamma(\gamma_k).
\eeq
For the calculation of $\wwlm_{\nu}(\bm{\gamma})$ we pass in the integral
\beqq\int\limits_{t_{i,j},s_{i,j}>0,t_i>0}\prod_{i<j}t_{i,j}^{\ga_{i,j}+\nu'_i-\nu'_j}
s_{i,j}^{\delta_{i,j}+\nu'_i+\nu'_j}e^{ -p_{i,j}-q_{i,j}}\frac{dt_{i,j}}{t_{i,j}}\frac{ds_{i,j}}{s_{i,j}}\prod_it_i^{\ga_i+2\nu'_i}e^{-p_i}\frac{dt_i}{t_i}\eeqq
from the integration variables $t_{i,j}$,  $s_{i,j}$ and $t_i$ to $p_{i,j}$, $q_{i,j}$ and  $p_i$.
Here
$$\nu'=\nu+\rho,\qquad \nu'_i=\nu_i+n-i+\frac{1}{2}.$$
Substitution of \rf{b7} gives the relation
\beqq
\prod_it_i^{\gamma_i}\prod_{i<j}t_{i,j}^{\ga_{i,j}}
s_{i,j}^{\delta_{i,j}}t_i^{\gamma_i}=\prod_ip_i^{\varphi_i}\prod_ {i<j} p_{i,j}^{\varphi_{i,j}}q_{i,j}^{\psi_{i,j}}\left(p_{i,j}+q_{i,j}\right)^{\theta_{i,j}}\eeqq
where
\begin{equation}\label{b19}\begin{split}
\varphi_{i,j}&=-\delta_{i,j}+\sum_{k>i}\left(\gamma_{k,j+1}+\delta_{k,j+1}\right)-\sum_{k>i} \left(\gamma_{k,j}+\delta_{k,j}\right)\\
\psi_{i,j}&=\ \delta_{i,j}-\sum_{k\geq i}\left(\gamma_{k,j+1}+\delta_{k,j+1}\right)+\sum_{k>i}\left( \gamma_{k,j}+\delta_{k,j}\right)\\
\theta_{i,j}&=-\delta_{i,j}-\gamma_{i,j}+\delta_{i,j+1}+\gamma_{i,j+1},
\end{split}\end{equation}
for $1\leq i<j<n$, and
\begin{equation}\label{b20}
\begin{split}
\varphi_{i,n}& =-\delta_{i,n}-\sum_{k>i}\left( \gamma_{kn}+\delta_{kn}-\gamma_k\right),\\
\psi_{i,n}&= \delta_{i,n}-\ga_i+\sum_{k>i}(\ga_{kn}+\delta_{kn}-\ga_k),\\
\theta_{i,n}&=\ga_i-\delta_{i,n}-\ga_{i,n},\qquad \varphi_i=-\ga_i
\end{split}
\end{equation}
for $1\leq i<n$.
 Using   \rf{int} we present the multiple integral
$$
\int_{p_{i,j},q_{i,j},p_i>0}\prod_{i<j} p_{i,j}^{{\varphi}_{i,j}-\nu'_i+\n'u_j}q_{i,j}^{{\psi_{i,j}}-\nu'_i-\nu'_j}(p_{i,j}+q_{i,j})^{\theta_{i,j}}\frac{dp_{i,j}}{p_{i,j}}\frac{dq_{i,j}}{q_{i,j}}\prod_i p_i^{\varphi_i-2\nu'_i}\frac{dp_i}{p_i}
$$
as the product
$$ \prod_{k=1}^n\Gamma(-\ga_k-2\nu'_k)\prod_{{i<j}}\Gamma (\varphi_{i,j}-\nu'_i+\nu'_j)\Gamma(\psi_{i,j}-\nu'_i-\nu'_j)\frac{\Gamma(\varphi_{i,j}+\theta_{i,j}+\psi_{i,j}-2
	\nu'_i)}{\Gamma(\varphi_{i,j}+\psi_{i,j}-2\nu'_i)}.$$
Substitution of \rf{b19} and \rf{b20} into the product over $1\leq i<j<n$ of ratios of $\Gamma$ functions in the latter expression results, just as for $\SO(n,n)$
to the product
\begin{equation*}
\prod_{i=1}^{n-1}\frac{\Gamma(-\gamma_{i,i+1}-\delta_{i,i+1}-2\nu'_i)}{\Gamma(-\gamma_{i,n}-\delta_{i,n}-2\nu'_i)}.
\end{equation*}
while the same product over $j=n$, $i=1,...,n-1$   is
\beqq
\prod_{i=1}^{n-1}\frac{\Gamma(-\delta_{i,n}-\gamma_{i,n}-2\nu'_i)} {\Gamma(-\gamma_i-2\nu'_i)}
.
\eeqq
Thus we have the following expressions for the left Whittaker vector
\begin{equation*}
\begin{split}
\wwlm_{\nu}(\bm{\gamma})=\, &\Gamma(-\gamma_n-2\nu'_n)\prod_{i=1}^{n-1}\ {\Gamma(-\gamma_{i,i+1}-\delta_{i,i+1}-2\nu'_i)}\cdot
\\ 
& \prod_{1\leq i<j\leq n}\Gamma({\varphi}_{i,j}-\nu'_i+\n'u_j)\Gamma({\psi}_{i,j}-\nu'_i-\nu'_j).
\end{split}
\end{equation*}
which results after substitution $\nu'=i\lambda$ in the following expression for Whittaker wave function:
\begin{align*}
		&\Psi_{\bm{\lambda}}(\x)= \frac{e^{-i(\bm{\lambda},\x)}}{(2\pi i)^d}\int\limits_{C_n}
		\exp \tilde{H}(\bm{x},\bm{\gamma})
		\prod_{k=1}^{n-1}{\Gamma(\gamma_{k,k+1}+\delta_{k,i+1}+2i\lambda_k)}
		\Gamma(\gamma_k)d\gamma_k\cdot \notag\\
		&\!\!\prod_{k<l} \Gamma(-{\varphi}_{k,l}+i(\lambda_k-\lambda_l))\Gamma(-{\psi}_{k,l}+i(\lambda_k+\lambda_l))
		\Gamma(\gamma_{k,l})\Gamma(\delta_{k,l})d\gamma_{k,l}d\delta_{k,l}
		\Gamma(\gamma_n+2i\lambda_n)d\gamma_n
	\end{align*}
Here $\tilde{H}(\bm{x},\bm{\gamma})$ is given in \rf{b17}, $\psi_{k,l}$ and $\vf_{k,l}$ are given in \rf{b19} and \rf{b20},
$d=n^2$. The contour $C$ is a deformation of imaginary plane $\Re \ga_{i,j}>0$,\ $\Re \delta_{i,j}>0$,\ $\Re \ga_{i}>0$ into
the strip of analyticity of the integrand.

 The formula for $\Psinew_{\bm{\lambda}}(\x)$  can be  simplified by using the change of variables
\beq\label{b26}\begin{split} &\gamma_{i,j}=\tgamma_{i,j}-\tgamma_{i+1,j},\qquad
	 \delta_{i,j}=\tdelta_{i,j}-\tdelta_{i+1,j},\qquad 1\leq i<j \leq n,\\
&\gamma_i=\tgamma_i-\tgamma_{i+1},\qquad \qquad\qquad\qquad\qquad\qquad\qquad 1\leq i\leq n,	
	 	\end{split}\eeq
 	that is
\beqq \begin{split}
&\tgamma_{i,j}=\gamma_{i,j}+\ldots+\gamma_{j-1,j},\qquad \tdelta_{i,j}=\delta_{i,j}+\ldots+\delta_{j-1,j},
\qquad
\tgamma_{i}=\gamma_i+\ldots+\gamma_n	
	\end{split}
\eeqq	
In this variables we have
\beqq
\begin{split}
&\vf_{i,j}=\tgamma_{i+1,j+1}+\tdelta_{i+1,j+1}-\tgamma_{i+1,j}-\tdelta_{i,j},\qquad 1\leq i<j<n,\\	
&\psi_{i,j}=-\tgamma_{i,j+1}-\tdelta_{i,j+1}+\tgamma_{i+1,j}+\tdelta_{i,j},\qquad 1\leq i<j<n,\\
&\vf_{i,n}=-\tgamma_{i+1,n}-\tdelta_{i,n}+\tgamma_{i+1},\qquad
\psi_{i,n}=\tgamma_{i+1,n}+\tdelta_{i,n}-\tgamma_i.	\ i=1,..., n.\end{split}
\eeqq
so that each Gamma function in the integrand depends now not more than of four summands.
 The exponential factor is now
$
	H(\bm{x},\tilde{\bm{\gamma}})=\sum_{j=2}^n(\tgamma_{1,j}+\tdelta_{1,j})(x_{j}-x_{j-1}) -\tgamma_1 x_n.$.
Finally we have
\begin{theorem}\label{t3} The function $\Psinew_{\blambda}(\x)$ is given by the integral
	\begin{align*}
	&\Psinew_{\blambda}(\x)= \frac{e^{  -i(\x,\blambda)}}{(2\pi i)^\d}\int\limits_{C}
	\exp {H}(\bm{x},{\bm{\gamma}})\Gamma(\gamma_n+2i\lambda_n)\Gamma(\gamma_n)d\gamma_n\cdot \\
	&
	\ \, \prod_{k=1}^{n-1}{\Gamma(\gamma_{k,k+1}+\delta_{k,k+1}+2i\lambda_k)}
	\Gamma(\gamma_k-\gamma_{k+1})d\gamma_k\cdot \\
	&\!\!\prod_{1\leq k<l\leq n}\!\! \Gamma({\xi}_{k,l}+i(\lambda_k-\lambda_l))\Gamma({\eta}_{k,l}+i(\lambda_k+\lambda_l))
	\Gamma(\gamma_{k,l}-\gamma_{k+1,l})\Gamma(\delta_{k,l}-\delta_{k+1,l})d\gamma_{k,l}d\delta_{k,l}.
	\notag\end{align*}	
\end{theorem}
Here we set $\gamma_{k,l}=\delta_{k,l}=0$ if the pair $(k,l)$ does not satisfies the condition $1\leq k<l\leq n$, and $\gamma_k=0$ if $k=1$ or $k=n+1$, $\d=n^2$,
$
H(\bm{x},{\bm{\gamma}})=\sum_{j=2}^n(\gamma_{1,j}+\delta_{1,j})(x_{j}-x_{j-1}) -\gamma_1 x_n$, and
\beqq
\begin{split}
	&\xi_{i,j}=-\gamma_{i+1,j+1}-\delta_{i+1,j+1}+\gamma_{i+1,j}+\delta_{i,j},\qquad 1\leq i<j<n,\\	
	&\eta_{i,j}=\gamma_{i,j+1}+\delta_{i,j+1}-\gamma_{i+1,j}-\delta_{i,j},\qquad 1\leq i<j<n,\\
	&\xi_{i,n}=\gamma_{i+1,n}+\delta_{i,n}-\gamma_{i+1},\qquad
	\eta_{i,n}=-\gamma_{i+1,n}-\delta_{i,n}+\gamma_i, \qquad 1\leq i<n.	\end{split}
\eeqq
The integration cycle is a deformation of the imaginary plain
$\Re \gamma_{k,l}=\Re \delta_{k,l}=\Re \gamma_k=0$ into nonempty domain  $D\subset\C^{\d}$ of the analyticity of the integrand, which is described by the relations
$
\Re \gamma_{i,j}>0$, \ $\Re \delta_{i,j}>0$,\ $\Re\, \xi_{i,j}>0$,\ $\Re \eta_{i,j}>0$,
\ $1\leq i<j\leq n$, \ \ $\Re \gamma_k>0$,\ $k=1,...,n$.

\setcounter{equation}{0}
\section{$\Sp(2 n,\R)$}
\subsection{BZ transform}
The split real form of the group $\Sp(2n,\C)$ is the group $\Sp(2n,\R)\subset \SL(2n,\R)$ preserving skew-symmetric form $\sum_{i=1}^{n} (x_iy_{\hat{i}}-x_{\hat{i}}y_i)$, where $\hat{i}=2n+1-i$. Positive roots are elements $\ve_i\pm\ve_j$ for $1\leq i<j\leq n$ and $2\ve_i$, where
$\ve_k\in\h^*$, $k=1,...,n$ are  defined by the condition $\ve_k( e_{j,j}-e_{\hat{j}\hat{j}})=\delta_{j,k}$.
We denote Chevalley generator of the Lie algebra $\g=\so_{2,n}$ by $e_i$ and $f_i$, $i=1,\ldots,n$  Here
\begin{align*} &e_i=e_{i,i+1}-e_{2n-i,2n+1-i},&&  f_i=e_{i+1,i}-e_{2n+1-i,2n-i},&&
i=1,...,n-1,\\
& e_{n}=e_{n,n+1},&& f_{n}=e_{n+1,n}.
\end{align*}
Denote by $s_i$ $i=1,...,n$  the corresponding generators of the Weyl group and
their lifts to the group $\Sp(2n,\R)$ according to \rf{W2a}.
 The normal ordering of the system $\Delta_+$  copies that of $\SO(n+1,n)$:
\beqq
\begin{split}&2\ve_n,\ (\ve_{n-1}+\ve_n,\, \ve_{n-1},\ve_{n-1}-\ve_n),\ \cdots, \\
	&(\ve_{1}+\ve_{2},\ve_{1}+\ve_{3},\cdots,\ve_{1}+\ve_{n},2\ve_1,\ve_{1}-\ve_{n},\cdots, \ve_{1}-\ve_{3},\ve_{1}-\ve_{2}).
\end{split}
\eeqq
It corresponds to reduced decomposition of the longest element $w_0^{(n)}$ of the Weyl group $W$ literally the same as \rf{b2}:
$$w_0^{(n)}=s_n({s}_{n-1}s_ns_{n-1})\,(s_{n-2}s_{n-1}s_n{s}_{n-1}s_{n-2})\,
\cdots\,(s_{1}s_2...s_{n-1}{s}_{n}s_{n-1}...s_2s_1).
$$
 Note that as element of $\Sp(2n,\R)$ ,
$$\bar{w}^{(n)}_0=(-1)^{n+1}\left(\sum_{k=1}^{n}e_{k\hat{k}}-e_{\hat{k},k}\right)
$$

Denote by $t_{i,j}$ Lusztig parameter corresponding to the root $\ve_i-\ve_j$, and  by $s_{i,j}$ Lusztig parameter corresponding to the root $\ve_i+\ve_j$, and by $t_k$ the parameter, corresponding to $2\ve_k$. Here $1\leq i<j\leq n$, $k=1,...,n$. Then the group element $\Xn{n}(\t)$ has a form \rf{b3} so that
BZ transform for $\Sp(2n,\R)$ is the solution of the equation \rf{gl5}, where in $\Xn{n}(\p)$ we use the notation
$p_{i,j}$ for Lusztig parameter, corresponding to the root $\ve_i-\ve_j$,  $q_{i,j}$ for the Lusztig parameter, corresponding to the root $\ve_i+\ve_j$, and $p_i$ for the parameters, corresponding to $2\ve_i$, see \rf{b4a}.

Keep the notation \rf{b5a}. We have	
\begin{theorem} a)\label{theorem7} BZ map for $Sp(2n,\R)$ looks as
	\begin{equation}\label{c2}
	 p_{i,j}=r_{i,j}\cdot\prod_{k<i}\frac{t_{k,j-1}}{s_{k,j-1}}\prod_{k<i}\frac{s_{k,j}}{t_{k,j}},
	\qquad    	
	q_{i,j}=r_{i,j}\cdot\prod_{k<i}\frac{t_{k,j-1}}{s_{k,j-1}}\prod_{k\leq i}\frac{s_{k,j}}{t_{k,j}},\qquad
	p_i=\frac{u_{i,n}^2}{t_is_{i,n}^2}\prod_{k<i}	\frac{t^2_{kn}}{s^2_{kn}}        		
	\end{equation}				
	b) The twisting matrix $T\in H$ is $T=\sum_{k=1}^n\left(T_{k}e_{k,k}+T_k^{-1}e_{\hat{k},\hat{k}}\right)$, where
	\beqq
	 T_k=t_k\prod_{j<k}\frac{s_{j,k}}{t_{j,k}}\cdot\prod_{j>k}s_{k,j}t_{k,j}
	\eeqq
	c) BZ transform preserves the measure $\frac{d\t}{\t}$:
	$$\prod_{i<j}\frac{dt_{i,j}ds_{i,j}}{t_{i,j}s_{i,j}}\prod_k\frac{dt_k}{t_k}= \prod_{i<j}\frac{dp_{i,j}dq_{i,j}}{p_{i,j}q_{i,j}}\prod_k\frac{dp_k}{p_k} .$$				
\end{theorem}	
Again, we have the same remarks with the same proofs:

\noindent {\bf Remark 1}. BZ transform $\tau$ is involutive, that is the inverse relation
\begin{equation}\label{c4}
t_{i,j}=\tilde{r}_{i,j}\cdot\prod_{k<i}\frac{p_{k,j-1}}{q_{k,j-1}}\prod_{k<i}\frac{q_{k,j}}{p_{k,j}},
\qquad    	
s_{i,j}=\tilde{r}_{i,j}\cdot\prod_{k<i}\frac{p_{k,j-1}}{q_{k,j-1}}\prod_{k\leq i}\frac{q_{k,j}}{p_{k,j}},\qquad
t_i=\frac{v^2_{i,n}}{p_iq^2_{i,n}}\prod_{k<i}	\frac{p^2_{kn}}{q^2_{kn}}        		
\end{equation}
has the same form as \rf{d5}. Here
\begin{equation*}
v_{i,j}=p_{i,j}+q_{i,j},\qquad \tr_{i,j}=\frac{v_{i,j-1}}{v_{i,j}q_{i,j-1}}.
\end{equation*}
\medskip
\noindent {\bf Remark 2}. The Cartan twist $T$ can be as well  written  as
$T(\t)=T(\p^{-1})$, that is $T=\sum_{k=1}^n\left(T_{k}e_{k,k}+T_k^{-1}e_{\hat{k},\hat{k}}\right)$, where
\beqq
T_k=\frac{1}{p_k}\prod_{j<k}\frac{q_{j,k}}{p_{j,k}}\cdot\prod_{j>k}\frac{1}{q_{k,j}p_{k,j}}
\eeqq

\begin{corollary}\label{corollary5}
	The restriction of left Whittaker vector to $N_+$ is given by the function
	\begin{equation*}
	 \wwl_{\nu}(\t)=\t^{\nu}\cdot\exp\Big(-\sum_{i<j}(p_{i,j}+q_{i,j})-\sum_kp_k\Big)=\p^{-\nu}\cdot\exp\Big(-\sum_{i<j}(p_{i,j}+q_{i,j})-\sum_kp_k)\Big)
	\end{equation*}
	where $p_{i,j}$  and $p_k$ are given by \rf{c2}  and $\t^{\nu}$ for a weight $\nu=\sum_{k=1}^n \nu_k\ve_k$ means the product
	\beqq
	\t^{\nu}=\prod_{i<j} t_{i,j}^{\nu_i-\nu_j}s_{i,j}^{\nu_i+\nu_j}\prod_kt_k^{\nu_k}
	\eeqq
\end{corollary}
\bigskip

For inductive proof Theorem \ref{theorem7} we use the
factorization
$\Xn{n}(-\t)=\tXn{n-1}(-\t)\cdot \An{n}(-\t)$, where $\tXn{n-1}(-\t)$ and $\An{n}(-\t)$ are given by the expressions \rf{b13a} and \rf{b13b}.
Denote by $\Ln{n}$ the diagonal matrix $\Ln{n}=\sum_{k=1}^n(\Ln{n}_{k,k}e_{k,k}+{\Ln{n}_{k,k}}^{-1}e_{\hat{k},\hat{k}})$, where
\beqq
\Ln{n}_{1,1}=t_1\prod_{k=1}^{n-1}t_{1,k}s_{1,k}\qquad \Ln{n}_{k,k}=\frac{t_{1,k}}{s_{1,k}},\qquad 1<k\leq n,	
\eeqq
and define the variables $\tilde{t}_{i,j}$ and $\tilde{s}_{i,j}$, $1< i<j\leq n$ by the relations
\beqq
\tilde{t}_{i,j}=t_{i,j}\cdot\frac{s_{1,j-1}}{t_{1,j-1}}\frac{t_{1,j}}{s_{1,j}}, \qquad
\tilde{s}_{i,j}=s_{i,j}\cdot\frac{s_{1,j-1}}{t_{1,j-1}}\frac{t_{1,j}}{s_{1,j}},\qquad
\tilde{t}_i=t_i\frac{s^2_{1,n}}{t^2_{1,n}}
\eeqq
Then we have again
\begin{lemma}\label{lemma6} Nondiagonal matrix elements $U_{i,j}$ ($i\not=j$) of the matrix \rf{gl12} equal zero if $j>1$ and $i<2n+1$.
	Matrix element $U_{1,1}$ equals $t_1\prod_{j>1} t_{1,j}s_{1,j}$. Matrix element $U_{k,k}$ equals $\dfrac{s_{1,k}}{t_{1,k}}$ for $1<k\leq n$. Matrix element $U_{\hat{k},\hat{k}}$ equals $\dfrac{t_{1,k}}{s_{1,k}}$ for $1<k\leq n$. 	
\end{lemma}
\begin{proposition}\label{prop7} {\em{a)}} The parameters $p_{1,j}$, $q_{1,j}$, $j=2,...,n$ and $p_1$ of $\An{n}(\p)$ are equal to
	\beqq
	 p_{1,j}=\frac{u_{1,j-1}}{u_{1,j}s_{1,j-1}}, \qquad
	q_{1,j}=\frac{u_{1,j-1}s_{1,j}}{u_{1,j}s_{1,j-1}t_{1,j}},\qquad p_1=\frac{u^2_{1,n}}{t_1s^2_{1,n}}
	\eeqq
	{\em{b)}} Cartan twist $T$ and the matrix $\tXn{n-1}(\p)$ satisfy the relation
	\beqq
	T \tXn{n-1}(\p)=\Ln{n}\cdot\left(\tXn{n-1}(-\tilde{\t})\ww{n-1}\right)_{0+}
	\eeqq
\end{proposition}
These statement are proved in the same manner as in Section \ref{section2}. They are sufficient for inductive proof of Theorem \ref{theorem7}. \hfill{$\square$}

\subsection{Whittaker function}\label{section5.2}
We start again with the action of Cartan
generators on the right Whittaker vector. The vector
$\exp(\sum_{k=1}^n  x_k (e_{\hat{k},\hat{k}}-e_{k,k}))\cdot\wr_\mu$
is presented by the  function  $\exp(-\sum_k \mu_k x_k)\wwr_\mu(\t')$, where
\begin{equation}\label{c13a} t'_{i,j}=t_{i,j}\cdot\exp\left(x_{j}-x_{j-1} \right),\qquad
s'_{i,j}=s_{i,j}\cdot\exp\left(x_{j}-x_{j-1} \right),\qquad t'_i=t_i\cdot\exp (-2x_n).
\end{equation}
Denote by $\wwlm_{\nu}(\bm{\ga})$ and $\wwrm_\mu(\bm{\ga})$ the Mellin transforms of the functions
$\wwl_{\nu}(\t)$ and $\wwr_\mu(\t)$,
\begin{equation*}\begin{split}& \wwrm_\mu(\bm{\ga})=\int\limits_{\t>0}\wwr_\mu(\t) \t^{\mathbf{\ga}}\frac{d\t}{\t}:=
\!\!\int\limits_{\substack{t_{i,j}>0,\\ s_{i,j}>0,t_i>0}}\!\! \wwr_\mu(\t)\prod_{i<j}t_{i,j}^{\ga_{i,j}}s_{i,j}^{\delta_{i,j}}t_i^{\gamma_i}\frac{dt_{i,j}}{t_{i,j}}\frac{ds_{i,j}}{s_{i,j}}\frac{dt_i}{t_i},\\
& \wwlm_{\nu}(\bm{\ga})=\int\limits_{\t>0}\t^\rho\wwl_{\nu}(\t) \t^{\mathbf{\ga}}\frac{d\t}{\t}:=
\!\!\int\limits_{\substack{t_{i,j}>0,\\ s_{i,j}>0,t_i>0}}\!\! \wwl_{\nu}(\t)\prod_{i<j}t_{i,j}^{\ga_{i,j}+j-i}s_{i,j}^{\delta_{i,j}+2n+2-i-j}t_i^{\gamma_i+n-i+1}\frac{dt_{i,j}}{t_{i,j}}
\frac{ds_{i,j}}{s_{i,j}}\frac{dt_i}{t_i}.\end{split}
\end{equation*}
Here $\rho=\sum_{k=1}^n (n+1-k)(\ve_{k,k}-\ve_{\hat{k},\hat{k}})$.
Due to \rf{c13a} the action of Cartan subgroup on the right Whittaker vector $\wr_\mu$ transforms the function $\wwrm_\mu(\bm{\ga})$ to the product
$
\exp(-\sum_k \mu_k x_k)\cdot\exp \tilde
{H}(\bm{x},\bm{\gamma})\wwrm_\mu(\bm{\ga})
$,
where
\beq\label{c12}
	\tilde{H}(\bm{x},\bm{\gamma})=\sum_{i<j}(\ga_{i,j}+\delta_{i,j})(x_{j}-x_{j-1}) -2x_n\sum_{i}\gamma_i.\eeq
The right Whittaker vector is presented by the function $\wwrm_\mu(\bm{\gamma})$ of the form \rf{b18}.
For the calculation of $\wwlm_{\nu}(\bm{\gamma})$ we pass in the integral
\beqq\int\limits_{t_{i,j},s_{i,j}>0,t_i>0}\prod_{i<j}t_{i,j}^{\ga_{i,j}+\nu'_i-\nu'_j}
s_{i,j}^{\delta_{i,j}+\nu'_i+\nu'_j}e^{ -p_{i,j}-q_{i,j}}\frac{dt_{i,j}}{t_{i,j}}\frac{ds_{i,j}}{s_{i,j}}\prod_it_i^{\ga_i+2\nu'_i}e^{-p_i}\frac{dt_i}{t_i}\eeqq
from the integration variables $t_{i,j}$,  $s_{i,j}$ and $t_i$ to $p_{i,j}$, $q_{i,j}$ and  $p_i$.
Here
$$\nu'=\nu+\rho,\qquad \nu'_i=\nu_i+n-i+1.$$
 Substitution of \rf{c4} gives the relation
$$\prod_it_i^{\gamma_i}\prod_{i<j}t_{i,j}^{\ga_{i,j}}
s_{i,j}^{\delta_{i,j}}t_i^{\gamma_i}=\prod_ip_i^{\varphi_i}\prod_ {i<j} p_{i,j}^{\varphi_{i,j}}q_{i,j}^{\psi_{i,j}}\left(p_{i,j}+q_{i,j}\right)^{\theta_{i,j}}$$
where as before
\begin{equation}\label{c13}\begin{split}
\varphi_{i,j}&=-\delta_{i,j}+\sum_{k>i}\left(\gamma_{k,j+1}+\delta_{k,j+1}\right)-\sum_{k>i} \left(\gamma_{k,j}+\delta_{k,j}\right)\\
\psi_{i,j}&=\ \delta_{i,j}-\sum_{k\geq i}\left(\gamma_{k,j+1}+\delta_{k,j+1}\right)+\sum_{k>i}\left( \gamma_{k,j}+\delta_{k,j}\right)\\
\theta_{i,j}&=-\delta_{i,j}-\gamma_{i,j}+\delta_{i,j+1}+\gamma_{i,j+1},
\end{split}\end{equation}
for $1\leq i<j<n$, and
\begin{equation}\label{c14}
\begin{split}
\varphi_{i,n}& =-\delta_{i,n}-\sum_{k>i}\left( \gamma_{kn}+\delta_{kn}-2\gamma_k\right),\\
\psi_{i,n}&= \delta_{i,n}-2\ga_i+\sum_{k>i}(\ga_{kn}+\delta_{kn}-2\ga_k),\\
\theta_{i,n}&=2\ga_i-\delta_{i,n}-\ga_{i,n},\qquad \varphi_i=-\ga_i
\end{split}
\end{equation}
for $1\leq i<n$.
 Using the  \rf{int} we present the multiple integral
$$
\int_{p_{i,j},q_{i,j},p_i>0}\prod_{i<j} p_{i,j}^{{\varphi}_{i,j}-\nu'_i+\nu'_j}q_{i,j}^{{\psi_{i,j}}-\nu'_i-\nu'_j}(p_{i,j}+q_{i,j})^{\theta_{i,j}}\frac{dp_{i,j}}{p_{i,j}}\frac{dq_{i,j}}{q_{i,j}}\prod_i p_i^{\varphi_i-\nu'_i}\frac{dp_i}{p_i}
$$
as the product
$$ \prod_{k=1}^n\Gamma(-\ga_k-\nu_k)\prod_{{i<j}}\Gamma (\varphi_{i,j}-\nu'_i+\nu'_j)\Gamma(\psi_{i,j}-\nu'_i-\nu'_j)\frac{\Gamma(\varphi_{i,j}+\theta_{i,j}+\psi_{i,j}-2
	\nu'_i)}{\Gamma(\varphi_{i,j}+\psi_{i,j}-2\nu'_i)}.$$
Substitution of \rf{c13}  into the product over $1\leq i<j<n$ of ratios of $\Gamma$ functions in the latter expression results, just as for $\SO(n,n)$
to the product
\begin{equation*}
\prod_{i=1}^{n-1}\frac{\Gamma(-\gamma_{i,i+1}-\delta_{i,i+1}-2\n'u_i)}{\Gamma(-\gamma_{i,n}-\delta_{i,n}-2\nu'_i)}.
\end{equation*}
while due to \rf{c14} the same product over $j=n$, $i=1,...,n-1$   is
\beq\label{c16}
\prod_{i=1}^{n-1}\frac{\Gamma(-\delta_{i,n}-\gamma_{i,n}-2\nu'_i)}{\Gamma(-2\gamma_i-2\nu'_i)}
.
\eeq
Thus we have the following expressions for the left Whittaker vector
\beqq
\begin{split}
\wwlm_{\nu}(\bm{\gamma})=\Gamma(-\gamma_n-\nu'_n)\prod_{i=1}^{n-1}\ & {\Gamma(-\gamma_{i,i+1}-\delta_{i,i+1}-2\nu'_i)}\cdot
\frac{\Gamma(-\gamma_i-\nu'_i)}{\Gamma(-2\gamma_i-2\nu'_i)}\cdot\\
&\prod_{1\leq i<j\leq n}\Gamma({\varphi}_{i,j}-\nu'_i+\nu'_j)\Gamma({\psi}_{i,j}-\nu'_i-\nu'_j)
\end{split}\eeqq
and for the Whittaker wave function:
	\begin{align*}
	\Psinew_{\bm{\lambda}}(\x)&= \frac{e^{-i(\bm{\lambda},\x)}}{(2\pi i)^d} \int\limits_{C} \!\!
	\exp \tilde{H}(\bm{x},\bm{\gamma})
	\prod_{k=1}^{n-1}\frac{\Gamma(\gamma_{k,k+1}+\delta_{k,k+1}+2i\lambda_k)}
	{\Gamma(2\gamma_k+2i\lambda_k)}\cdot \prod_{k=1}^n{\Gamma(\gamma_k+i\lambda_k)}\Gamma(\gamma_k)d\gamma_k\cdot
	 \notag\\
	&
\!\!	\prod_{1\leq k<l\leq n} \!\!\Gamma(-{\varphi}_{k,l}+i(\lambda_k-\lambda_l))\Gamma(-{\psi}_{k,l}+
	i(\lambda_k+\lambda_l))
	\Gamma(\gamma_{k,l})\Gamma(\delta_{k,l})d\gamma_{k,l}d\delta_{k,l}\cdot 
	\end{align*}
Here $H(\bm{x},\bm{\gamma})$ is given in \rf{b17}, $\vf_{i,j}$ and $\phi_{i,j}$ are given in \rf{b19} and \rf{b20}, $d=n^2$. The contour $C$ is the same as in the previous section.

The change of variables \rf{b26} simplifies the formula for Whittaker wave functions. Now we have
\beqq
\begin{split}
	&\vf_{i,j}=\tgamma_{i+1,j+1}+\tdelta_{i+1,j+1}-\tgamma_{i+1,j}-\tdelta_{i,j},\qquad 1\leq i<j<n,\\	
	&\psi_{i,j}=-\tgamma_{i,j+1}-\tdelta_{i,j+1}+\tgamma_{i+1,j}+\tdelta_{i,j},\qquad 1\leq i<j<n,\\
	&\vf_{i,n}=-\tgamma_{i+1,n}-\tdelta_{i,n}+2\tgamma_{i+1},
	\psi_{i,n}=\tgamma_{i+1,n}+\tdelta_{i,n}-2\tgamma_i, \qquad 1\leq i<n,	\end{split}
\eeqq
 Finally we have
 \begin{theorem}\label{t4} The function $\Psinew_{\blambda}(\x)$ is given by the integral
 	\begin{align*}
 	&\Psinew_{\blambda}(\x)= \frac{e^{  -i(\x,\blambda)}}{(2\pi i)^\d}\int\limits_{C} \!\!
 	\exp {H}(\bm{x},\bm{\gamma})
 	\prod_{k=1}^{n-1}\frac{\Gamma(\gamma_{k,k+1}+\tdelta_{k,k+1}+2i\lambda_k)}
 	{\Gamma(2\gamma_k-2\gamma_{k+1}+2\lambda_k)}\cdot
 	\\
 	&
 	\prod_{k<l}\left( \Gamma({\xi}_{k,l}+i(\lambda_k-\lambda_l))\Gamma({\eta}_{k,l}+
 	i(\lambda_k+\lambda_l))\right.
 	\left.
 	\Gamma(\gamma_{k,l}-\gamma_{k+1,l})\Gamma(\delta_{k,l}-\delta_{k+1,l})d\gamma_{k,l}d\delta_{k,l}\right)\cdot \\
 	& \prod_{k=1}^n{\Gamma(\gamma_k-\gamma_{k+1}+i\lambda_k)}\Gamma(\gamma_k-\gamma_{k+1})d\gamma_k.
 	\end{align*}	
 \end{theorem}
 Here we set $\gamma_{k,l}=\delta_{k,l}=0$ if the pair $(k,l)$ does not satisfies the condition $1\leq k<l\leq n$, and $\gamma_k=0$ if $k=1$ or $k=n+1$, $\d=n^2$ is the dimension of the maximal unipotent subgroup of $\Sp(2n,\R)$,
 $
 H(\bm{x},{\bm{\gamma}})=\sum_{j=2}^n(\gamma_{1,j}+\delta_{1,j})(x_{j}-x_{j-1}) -2\gamma_1 x_n$,
 \beqq
 \begin{split}
 	&\xi_{i,j}=-\gamma_{i+1,j+1}-\delta_{i+1,j+1}+\gamma_{i+1,j}+\delta_{i,j},\qquad 1\leq i<j<n,\\	
 	&\eta_{i,j}=\gamma_{i,j+1}+\delta_{i,j+1}-\gamma_{i+1,j}-\delta_{i,j},\qquad 1\leq i<j<n,\\
 	&\xi_{i,n}=\gamma_{i+1,n}+\delta_{i,n}-2\gamma_{i+1},\qquad 1\leq i<n,\\
 	&\eta_{i,n}=-\gamma_{i+1,n}-\delta_{i,n}+2\gamma_i.	\end{split}
 \eeqq
 The integration cycle is a deformation of the imaginary plain
 $\Re \gamma_{k,l}=\Re \delta_{k,l}=\Re \gamma_k=0$ into nonempty domain  $D\subset\C^{\d}$ of the analyticity of the integrand, which is described by the relations
 \beqq
 \Re \gamma_{i,j}>0, \ \Re \delta_{i,j}>0,\ \ Re\, \xi_{i,j}>0,\ \Re \eta_{i,j}>0,
 \ 1\leq i<j\leq n, \ \ \Re \gamma_k>0,\ k=1,...,n.
 \eeqq

\setcounter{equation}{0}
\section{Mellin transforms of Whittaker functions} \label{section6}
 The presentations  of Whittaker functions  given in Theorems \ref{t1}, \ref{t2}, \ref{t3},\ref{t4}, have a form
 which is easy to interpret as an inverse Mellin transform. This enables one to write down precise expressions for direct Mellin transforms of Whittaker functions.
 \medskip

{\bf 1. $ \GL(n)$}
Using the notations of Section \ref{section2.2} we denote $z_k=\exp x_k$, $z_{k,k+1}=\exp(x_k-x_{k+1})$,
$\z=\{z_1,\ldots,z_n\}$.
Set \beqq
\Lambda_i=\frac{(n-1)\lambda_i-\sum_{j\not=i}\lambda_j}{n},\qquad \Lambda_{i,j}=\Lambda_i+\Lambda_{i+1}+\ldots \Lambda_j
\eeqq
and put $s_{k,k+1}=\tgamma_{1,n+1-k}+\Lambda_{1,k}$.
 Then the formula \rf{I1} can be written as follows
\beqq
	\Psi_{\bm{\lambda}}(\z)= \frac{(z_1\cdots z_n)^{-\frac{i(\lambda_1+\ldots +\lambda_n)}{n}}}{(2\pi i)^{n-1}}\!\!\!\!\!\!\!\int\limits_{\Re s_{k,k+1}=0+}\!\!\!\!z_{1,2}^{-s_{1,2}}\cdots z_{n-1,n}^{-s_{n-1,n}} M_{\bm{\lambda}}^n(s_{12},\ldots, s_{n-1,n})ds_{12}
	\cdots ds_{n-1,n}
	\eeqq
where
\beq\label{M2} \begin{split}M_{\bm{\lambda}}^n & (s_{12},\ldots, s_{n-1,n})=\frac{1}{(2\pi i)^{n-1)(n-2)/2}}\cdot\\ &\int\limits_{C}\prod_{k=1}^{n-1}\Gamma(s_{k,k+1}-
\gamma_{2,n+1-k}+i(\lambda_1-\lambda_{k+1}-\Lambda_{1,k}))\Gamma(s_{k,k+1}-\gamma_{2,n+1-k}-i\Lambda_{1,k})\\
&	\prod_{1<k<l\leq n}{\Gamma}(\gamma_{kl}-\gamma_{k+1,l}+i(\lambda_k-\lambda_{n+k-l+1}))
	\Gamma(\gamma_{kl}-\gamma_{k+1,l+1})d\gamma_{k,l}
\end{split}\eeq	
The integration contour $C$ is a  deformation of the imaginary plane $\Re \gamma_{ij}=0$ to the region of analyticity of the integrand. The function $M_{\bm{\lambda}}^n  (s_{12},\ldots, s_{n-1,n})$ is then equal to the Mellin transform of $\SL(n)$ Whittaker function
$\bar{\Psi}_{\bm{\lambda}}(\z)={\Psi}_{\bm{\lambda}}(\z)\cdot(z_1\cdots z_n)^{\frac{i(\lambda_1+\ldots +\lambda_n)}{n}}:$
\beqq M_{\bm{\lambda}}^n  (s_{12},\ldots, s_{n-1,n})=\int\limits_{z_{i,j}>0} \bar{\Psi}_{\bm{\lambda}}(\z)
z_{1,2}^{s_{1,2}-1}\cdots z_{n-1,n}^{s_{n-1,n}-1}dz_{12}\cdots dz_{n-1,n}
\eeqq
 The relation \rf{M2} can be taken as a starting point for iterative construction of the Mellin transform of $\SL(n)$ part of the Whittaker function 	$\Psi_{\bm{\lambda}}(\z)$.
\medskip

{\bf 2. $ \SO(n,n)$.} Using the notations of Section \ref{section3.2} denote  $z_n= \exp(x_{n-1}+x_n)$ and  $z_i=\exp(x_{i}-x_{i+1})$ for
 $1\leq i\leq n-1$.
 Set
 \beqq\begin{split}
 	&\Lambda_i=\lambda_1+\ldots +\lambda_i, \qquad 1\leq i=1\leq n-2,\\
 	&\Lambda_{n-1}=\frac{\lambda_1+\ldots +\lambda_{n-1}+\lambda_n}{2},\qquad
 	\Lambda_{n}=\frac{\lambda_1+\ldots +\lambda_{n-1}-\lambda_n}{2},\\
& s_{n-1}=\tgamma_{1,n}+i\Lambda_{n-1},\qquad
 s_{n}=\tdelta_{1,n}+i\Lambda_{n},\\
 & s_i=\tgamma_{1,i+1}+\tdelta_{1,i+1}+i\Lambda_i, \qquad
  t_i=\tgamma_{1,i+1}-\tdelta_{1,i+1},  \qquad 1\leq i\leq n-2,\\
 &\ts_j=s_i-i\Lambda_j=\tgamma_{1j}+\tdelta_{1,j}, \qquad 1\leq j\leq n.
 \end{split}\eeqq
 In this notations the formula of Theorem \ref{t2} for the Whittaker function is written as the inverse Mellin transform:
 \beqq
\Psi_{\bm{\lambda}}(\z)=\frac{1}{(2\pi i)^n} \int\limits_{\Re s_i=0+}
z_1^{-s_1}\cdots z_n^{-s_n} M_{\bm{\lambda}}^n(s_{1},\ldots, s_{n})ds_{1}
\cdots ds_{n},\eeqq
where
\beqq \begin{split}
& M_{\bm{\lambda}}^n(s_{1},\ldots, s_{n})=\frac{1}{(2\pi i)^{(n-1)^2}}\int\limits_{C_n}
\frac{\Gamma(\ts_1+2i\lambda_1)\Gamma(\ts_{n-1}-\gamma_{2,n})\Gamma(\ts_n-\delta_{2,n})}{\Gamma(\ts_{n-1}+\ts_n-\gamma_{2,n}-\delta_{2,n}+2i\lambda_1)}
\\
&\prod_{k=1}^{n-2}\Gamma\left(\frac{\ts_k+t_k}{2}-\gamma_{2,k+1}\right)
\Gamma\left(\frac{\ts_k-t_k}{2}-\delta_{2,k+1}\right)\cdot L_{\bm{\lambda}}^n \,\prod_{k=1}^{n-2}dt_k\prod_{2\leq k<l\leq n}
d\gamma_{k,l}d\delta_{k,l}\end{split}\eeqq
The contour $C$ is a deformation of the  imaginary plane to the region of analyticity of the integrand. Here
\beqq\begin{split}
L_{\bm{\lambda}}^n&=\prod_{k=2}^{n-2}\frac{\Gamma(\gamma_{k,k+1}+\delta_{k,k+1}+2i\lambda_k)}{\Gamma(\gamma_{k,n}	+\delta_{k,n}-\gamma_{k+1,n}-\delta_{k+1,n}+
	2i\lambda_k)}\prod_{1<k<l\leq n}\Gamma(\gamma_{k,l}-\gamma_{k+1,l})\Gamma(\delta_{k,l}-\delta_{k+1,l})\\
& \prod_{k=2}^n\Gamma(\gamma_k-\gamma_{k+1})\prod_{1\leq k<l\leq n} \Gamma({\xi}_{k,l}+i(\lambda_k-\lambda_l))\Gamma({\eta}_{k,l}+i(\lambda_k+\lambda_l))
\end{split}\eeqq
with
\beqq
\begin{split}
	&\xi_{i,j}=-\gamma_{i+1,j+1}-\delta_{i+1,j+1}+\gamma_{i+1,j}+\delta_{i,j},\qquad i<j<n,\\	
	&\eta_{i,j}=\gamma_{i,j+1}+\delta_{i,j+1}-\gamma_{i+1,j}-\delta_{i,j},\qquad i<j<n,\\
	&\xi_{i,n}=\gamma_{i,n}-\delta_{i+1,n},\qquad
	\eta_{i,n}=\delta_{i,n}-\gamma_{i+1,n},\qquad i<n,	\end{split}
\eeqq
for $i>1$ and
\beqq
\begin{split}
	&\xi_{1,j}=-\gamma_{2,j+1}-\delta_{2,j+1}+\gamma_{2,j}+\frac{\ts_{j-1}-t_{j-1}}{2},\qquad 1<j<n,\\	
	&\eta_{1,j}=\ts_j-\gamma_{2,j}-\frac{\ts_{j-1}+t_{j-1}}{2},\qquad 1< j<n,\\
	&\xi_{1,n}=\ts_{n-1}-\delta_{2,n},\qquad
	\eta_{1,n}=\ts_n-\gamma_{2,n}.	\end{split}
\eeqq
\medskip

{\bf 3. $\SO(n+1,n)$.} Using the notations of Section \ref{section4.2} denote  $z_n= \exp(x_n)$ and  $z_i=\exp(x_{i}-x_{i+1})$ for
$1\leq i\leq n-1$.
Set
\beqq\begin{split}
	&\Lambda_i=\lambda_1+\ldots +\lambda_i, \qquad 1\leq i=1\leq n,\\
	& s_i=\gamma_{1,i+1}+\delta_{1,i+1}+i\Lambda_i, \qquad
	t_i=\gamma_{1,i+1}-\delta_{1,i+1},  \qquad 1\leq i\leq n-1,\\
	& s_n=\gamma_1+i\Lambda_n,\qquad
	\ts_i=s_i-i\Lambda_i \qquad 1\leq i\leq n.
\end{split}\eeqq
In this notations the formula of Theorem \ref{t3} for the Whittaker function looks as follows:
\beqq
\Psi_{\bm{\lambda}}(\z)= \frac{1}{(2\pi i)^n}\int\limits_{\Re s_i=0+}
z_1^{-s_1}\cdots z_n^{-s_n} M_{\bm{\lambda}}^n(s_{1},\ldots, s_{n})ds_{1}
\cdots ds_{n},\eeqq
where
\beqq
\begin{split}
\Psi_{\bm{\lambda}}(\x)=&\frac{1}{(2\pi i)^{n(n-1)}} \int\limits_{C}
\prod_{k=1}^{n-1}\Gamma\left(\frac{\ts_k+t_k}{2}-\gamma_{2,k+1}\right)
\Gamma\left(\frac{\ts_k-t_k}{2}-\delta_{2,k+1}\right)\cdot\\
& \Gamma(\ts_1+2i\lambda_1)\Gamma(\ts_n-\gamma_2)\cdot L_{\bm{\lambda}}^n \prod_{k=1}^{n-1}dt_k\prod_{2\leq k<l\leq n}
d\gamma_{k,l}d\delta_{k,l}\prod_{k=2}^nd\gamma_k\end{split}\eeqq
The contour $C$ is a deformation of imaginary plane to the region of analyticity  of the integrand.
Here
\begin{align*}
L_{\bm{\lambda}}^n &=\Gamma(\gamma_n+2i\lambda_n)
\prod_{k=2}^{n-1}{\Gamma(\gamma_{k,k+1}+\delta_{k,k+1}+2i\lambda_k)}
\prod_{k=2}^{n}\Gamma(\gamma_k-\gamma_{k+1})\cdot \\
\prod_{1\leq k<l\leq n} \Gamma({\xi}_{k,l}+&i(\lambda_k-\lambda_l))\Gamma(\eta_{k,l}+i(\lambda_k+\lambda_l))
\prod_{2\leq k<l\leq n}\Gamma(\gamma_{k,l}-\gamma_{k+1,l})\Gamma(\delta_{k,l}-\delta_{k+1,l})
\notag\end{align*}
with
\beqq
\begin{split}
	&\xi_{i,j}=-\gamma_{i+1,j+1}-\delta_{i+1,j+1}+\gamma_{i+1,j}+\delta_{i,j},\qquad i<j<n,\\	
	&\eta_{i,j}=\gamma_{i,j+1}+\delta_{i,j+1}-\gamma_{i+1,j}-\delta_{i,j},\qquad i<j<n,\\
	&\xi_{i,n}=\gamma_{i,n}-\delta_{i+1,n},\qquad
	\eta_{i,n}=\delta_{i,n}-\gamma_{i+1,n}\qquad i<n,	\end{split}
\eeqq
for $i>1$ and
\beqq
\begin{split}
	&\xi_{1,j}=-\gamma_{2,j+1}-\delta_{2,j+1}+\gamma_{2,j}+\frac{\ts_{j-1}-t_{j-1}}{2},\qquad 1<j<n,\\	
	&\eta_{1,j}=\ts_j+\gamma_{2,j}-\frac{\ts_{j-1}+t_{j-1}}{2},\qquad 1< j<n,\\
	&\xi_{1,n}=\frac{\ts_1-t_1}{2}-\gamma_2+\gamma_{2,n},\qquad
	\eta_{1,n}=-\frac{\ts_1-t_1}{2}+\ts_1-\gamma_{2,n}.	\end{split}
\eeqq
\medskip

{\bf 4. $\Sp(2n,\R)$.}  Using the notations of Section \ref{section5.2} denote again $z_n= \exp(x_n)$ and  $z_i=\exp(x_{i}-x_{i+1})$ for
$1\leq i\leq n-1$.
Set
\beqq\begin{split}
	&\Lambda_i=\lambda_1+\ldots +\lambda_i, \qquad 1\leq i=1\leq n,\\
	& s_i=\gamma_{1,i+1}+\delta_{1,i+1}+i\Lambda_i, \qquad
	t_i=\gamma_{1,i+1}-\delta_{1,i+1},  \qquad 1\leq i\leq n-1,\\
	& s_n=2\gamma_1+i\Lambda_n,\qquad
	\ts_i=s_i-i\Lambda_i \qquad 1\leq i\leq n.
\end{split}\eeqq
In this notations the formula of Theorem \ref{t4}  looks as follows:
\beqq
\Psi_{\bm{\lambda}}(\z)= \frac{1}{(2\pi i)^n}\int\limits_{\Re s_i=0+}
z_1^{-s_1}\cdots z_n^{-s_n} M_{\bm{\lambda}}^n(s_{1},\ldots, s_{n})ds_{1}
\cdots ds_{n},\eeqq
where
\beqq
\begin{split}
&	\Psi_{\bm{\lambda}}(\x)= \frac{1}{(2\pi i)^{n(n-1)}}\int\limits_{C}
	\prod_{k=1}^{n-1}\Gamma\left(\frac{\ts_k+t_k}{2}-\gamma_{2,k+1}\right)
	\Gamma\left(\frac{\ts_k-t_k}{2}-\delta_{2,k+1}\right)\cdot\\
&	\frac{\Gamma(\ts_n-\gamma_2+i\lambda_1)}{\Gamma(\ts_n-2\gamma_2+2i\lambda_1)}
	 \Gamma(\ts_1+2i\lambda_1)\Gamma(\ts_n-\gamma_2)\cdot L_{\bm{\lambda}}^n \prod_{k=1}^{n-1}dt_k\prod_{2\leq k<l\leq n}
d\gamma_{k,l}d\delta_{k,l}\prod_{k=2}^nd\gamma_k\end{split}\eeqq
The contour $C$ is a deformation of the imaginary to the region of analyticity of the integrand.
Here
\begin{align*}
L_{\bm{\lambda}}^n &=\Gamma(\gamma_n+2i\lambda_n)
\prod_{k=2}^{n-1}{\Gamma(\gamma_{k,k+1}+\delta_{k,k+1}+2i\lambda_k)}
\prod_{k=2}^{n}\Gamma(\gamma_k-\gamma_{k+1})\cdot \\
\prod_{1\leq k<l\leq n} \Gamma(\xi_{k,l}+&i(\lambda_k-\lambda_l))\Gamma(\eta_{k,l}+i(\lambda_k+\lambda_l))
\prod_{2\leq k<l\leq n}\Gamma(\gamma_{k,l}-\gamma_{k+1,l})\Gamma(\delta_{k,l}-\delta_{k+1,l})\cdot
\\ &\Gamma(\gamma_n+i\lambda_n)
\prod_{k=2}^{n-1}
\frac{\Gamma(\gamma_k-\gamma_{k+1}+i\lambda_k)}{\Gamma(\gamma_k-\gamma_{k+1}+2i\lambda_k)}
\notag\end{align*}
with
\beqq
\begin{split}
	&\xi_{i,j}=-\gamma_{i+1,j+1}-\delta_{i+1,j+1}+\gamma_{i+1,j}+\delta_{i,j},\qquad i<j<n,\\	
	&\eta_{i,j}=\gamma_{i,j+1}+\delta_{i,j+1}-\gamma_{i+1,j}-\delta_{i,j},\qquad i<j<n,\\
	&\xi_{i,n}=\gamma_{i,n}-\delta_{i+1,n},\qquad
	\eta_{i,n}=\delta_{i,n}-\gamma_{i+1,n}\qquad i<n,	\end{split}
\eeqq
for $i>1$ and
\beqq
\begin{split}
	&\xi_{1,j}=-\gamma_{2,j+1}-\delta_{2,j+1}+\gamma_{2,j}+\frac{\ts_{j-1}-t_{j-1}}{2},\qquad 1<j<n,\\	
	&\eta_{1,j}=\ts_j+\gamma_{2,j}-\frac{\ts_{j-1}+t_{j-1}}{2},\qquad 1< j<n,\\
	&\xi_{1,n}=\frac{\ts_1-t_1}{2}-2\gamma_2+\gamma_{2,n}, \qquad
	\eta_{1,n}=-\frac{\ts_1-t_1}{2}+\ts_1-\gamma_{2,n}.	\end{split}
\eeqq

\setcounter{equation}{0}
\def\n{\hat{n}}
\section*{Appendix}
\appendix
\section{Proof of Lemma \ref{lemma3}}
   The proof essentially consists of calculation of the product of two matrices,
$U=\tilde{A}\cdot\left(\An{n}(\p)\right)^{-1}$, where $\tilde{A}=(\wb_0^{(n-1)})^{-1}A^{(n)}(-\t)\wb_0^{(n)}$ with matrix elements
$$\tilde{A}_{i,j}=\left\{ \begin{array}{cc} (-1)^{n-j}\An{n}_{1,n+1-j}(-\t),& i=1,\\
(-1)^{i+j}\An{n}_{n+2-i,n+1-j}(-\t)& i>1, \end{array}\right.$$
and
\begin{equation*}\left(A^{(n)}\right)^{-1}(\p)=\begin{pmatrix}1&-p_{1,n}&0&...&...&0\\
0&1&-p_{1,n-1}&...&...&0\\0&&...&...&&0\\0&0&0&...&1&-p_{1,2}\\0&0&0&...&0&1
\end{pmatrix}
\end{equation*}
so that
$\left(A^{(n)}\right)^{-1}_{i,j}(\p)=-p_{1\i}$, if $j=i+1$, and  $\left(A^{(n)}\right)^{-1}_{i,j}(\p)=0$, if $j>i+1$. We then have
\begin{align*} & U_{1,j}=(-1)^{n-j}\An{n}_{1,n+1-j}(-\t)+(-1)^{n-j+1}\An{n}_{1,n+2-j}(-\t)\cdot (-p_{1,\j+1}),&& i=1,\\
& U_{i,j}=(-1)^{i+j}\An{n}_{n+2-i,n+1-j}(-\t)+(-1)^{i+j+1}
\An{n}_{n+2-i,n+2-j}(-\t)\cdot(-p_{1,\j+1}),&& i>1
\end{align*}	
so that the relations \rf{gl10} and \rf{gl24} imply the equalities
\begin{equation*}\begin{split}& U_{1,1}=t_{1,2}t_{1,3}\cdots t_{1,n}=\frac{1}{p_{1,2}p_{1,3}\cdots p_{1,n}},\qquad U_{ii}=p_{1,\i+1}=\frac{1}{t_{1,i}}, \ \ i>1,\\
&	 U_{i,j}=0, \ i\not=j,\ j>1.
\end{split}
\end{equation*} \hfill{$\square$}
\setcounter{equation}{0}
\section{Proof of Lemma \ref{lemma4}} This proof is a analogous to that of Lemma \ref{lemma3} with more technical details which differ for $n$ even and odd. Assume first that $n$ is even. Denote  for simplicity of notations entries of the matrix $\An{n}(-\t)$ by $A_{i,j}$, entries of the matrix
$\An{n}(\p)$ by $C_{i,j}$, variables $t_{1,j}$, $s_{1,j}$, $u_{1,j}$, $p_{1,j}$, $q_{1,j}$ $v_{1,j}$ by $t_j$, $s_j$, $u_j$, $p_j$, $q_j$, and $v_j$ correspondingly.  The upper triangular matrix $\An{n}(-\t)$ has a natural block structure with matrix coefficients equal to
\begin{align}\notag & A_{i,j}=(-1)^{i+j}s_{i+1}s_{i+2}\cdots s_{j-1}u_j,&& 1\leq i\leq j\leq n,\\
&\notag A_{\i\j}=t_{j+1}t_{j+2}\cdots t_{i-1}u_i,&& 1\leq j\leq i\leq n,\\
&\label{D4} A_{i,\j}=(-1)^i(s_ns_{n-1}\cdots s_{i+1})\cdot(t_nt_{n-1}\cdots t_{j+1})&& 1\leq i<n,\ 1\leq j\leq n,\\
&\notag A_{n,\n}=0,\ \  A_{n,\j}=(-1)^i(s_nt_{n-1}\cdots t_{j+1}),&& \ 1\leq j< n,\\
&\notag A_{\i,j}=0,&& \ 1\leq i,j\leq n
\end{align}
Using \rf{d4} we can express the matrix $\tilde{A}= \left(\wb_0^{(n-1)}\right)^{-1}A^{(n)}(-\t)\wb_0^{(n)}$ as the sum
\begin{equation*}\begin{split}-\tilde{A}=&\sum_{i,k: i\not=1,n}\left(   A_{ik}e_{\i\k}+ A_{\i\k}e_{ik}+ A_{i\k}e_{\i   k}      \right)+\\& \sum_k\left( A_{1,k}e_{1\k}+ A_{1\k}e_{1,k}+ A_{n\k}e_{nk}+ A_{\n\k}e_{\n k}\right)+e_{\hat{1} 1}+e_{n\n}.
\end{split}
\end{equation*}
Next, for any matrix $X\in\SO(n,n)$ we have the relation
 $$\left(X^{-1}\right)_{i,j}=X_{\j \i}$$
so the matrix coefficients $U_{i,j}$ of the matrix
 $U=	\left(\ww{n-1}\right)^{-1}\An{n}(-\t)\ww{n}\left(\An{n}(\p)\right)^{-1}$ are
\begin{align*}&U_{i,\j}=-\sum_k A_{\i\k} C_{j \k},&
&U_{\i j}=-\sum_k A_{i\k} C_{\j\k},\\
&U_{i,j}=-\sum_k A_{\i\k} C_{\j \k},&
&U_{\i \j}=-\sum_k\left( A_{i\k} C_{j\k}+ A_{ik} C_{j,k}\right).
\end{align*}
and
\begin{align*} &U_{1\j}=\sum_k\left( A_{1,k} C_{j,k}+ A_{1\k} C_{j\k}\right),&&
U_{n\j}= C_{jn}+\sum_k A_{n\k} C_{j\k},
&&U_{\hat{1} j}= C_{\j\hat{1}},\\
&U_{1,j}=\sum_k A_{1\k} C_{\j\k},&&U_{nj}=\sum_k A_{n\k} C_{\j\k},
&&U_{\hat{1}\j}= C_{j\hat{1}}, \\	
&U_{\n j}=\sum_k A_{\n\k} C_{\j\k},&& U_{\n\j}=\sum_k A_{\n\k} C_{j\k},\end{align*}
Then the proof reduces to the check of identities
\begin{align}\label{D2}&\sum_k A_{\i\k} C_{j \k}=0,&&
&&\sum_k A_{i\k} C_{\j\k}=0,\\  &\sum_k A_{\i\k} C_{\j \k}=-\delta_{i,j}\frac{s_{i}}{t_{i}},&& \label{D6}
&&\sum_k\left( A_{i\k} C_{j\k}+ A_{ik}C_{j,k}\right)=-\delta_{i,j}
\frac{t_{i}}{s_{i}}
\end{align}
for $i,j\not=1,n$, and of special cases
\begin{equation*}
\begin{split}
&	U_{1,j}=U_{1\j}=U_{j\hat{1}}=U_{\j\hat{1}}=0, \qquad j\not=1, \qquad U_{n,n}=-\frac{s_n}{t_n}, \ U_{\n\n}=-\frac{t_n}{s_n}\\
& U_{nj}=0,\ j\not=1,n,\ U_{n\j}=0,  U_{\n\j}=0,\ j\not=n,\  U_{\n\j}=0, \ \j\not=n,\ \end{split}\end{equation*}

 We check here several of them. First (\ref{D2}). According to (\ref{d17}),
\begin{equation}
\label{D7} p_i=-\frac{u_{i-1}}{u_{i}s_{i-1}},\qquad q_i=-\frac{u_{i-1}s_i}{u_{i}t_is_{i-1}},\qquad v_i=-\frac{u_{i-1}}{t_is_{i-1}}	
\end{equation}
Here we assume $s_1=t_1=u_1=1$. The  relations (\ref{D2}) reduce to
\begin{equation}\label{D9} 1+u_ip_i\left(1+p_{i-1}t_{i-1}+p_{i-2}p_{i-1}t_{i-1}t_i+\cdots+ (p_2\cdots  p_{i-1})(t_2\cdots t_{i-1})\right)=0
\end{equation}
 Substitute (\ref{D7}):
\begin{equation*}\begin{split}&1+\left(-\frac{u_{i-1}}{s_{i-1}}+\frac{u_{i-2}t_{i-1}}{s_{i-1}s_{i-2}}+\cdots +(-1)^{i-1}\frac{t_{i-1}\cdots t_2}{s_{i-1}\cdots s_1}\right)=\\&
1-\left(1+\frac{t_{i-1}}{s_{i-1}}\right)+\frac{t_{i-1}}{s_{i-1}}\left(1+\frac{t_{i-2}}{s_{i-2}}\right)+
\cdots+(-1)^{i-1}\frac{t_{i-1}\cdots t_2}{s_{i-1}\cdots s_1}=0\end{split}
\end{equation*}
Compute the diagonal entries of $U$ for $i\not=1$:
\begin{equation*}\begin{split} U_{ii}=1+u_iv_i(1+t_{i-1}p_{i-1}+\cdots +(t_2\cdots t_{i-1})(p_2\cdots p_{i-1}))
\end{split}\end{equation*}
By the previous calculation, $1+t_{i-1}p_{i-1}+\cdots+ (t_2\cdots t_{i-1})(p_2\cdots p_{i-1})=-1/(u_ip_i)$, so that
$$U_{ii}=1-\frac{u_iv_i}{u_ip_i}=-\frac{s_i}{t_i}$$
Vanishing of nondiagonal entries $U_{i,j}$ for $i,j\not=1,n$,  as well as of $U_{1,j}$ for $j\not=1$ also uses the equality (\ref{D9}). Next, the element $U_{n,n}$ equals to the sum
$$\sum_{k} A_{n\k}C_{\n\k}=0+s_nv_n(1+t_{n-1}p_{n-1}+\cdots +(t_2\cdots t_{n-1})(p_2\cdots p_{n-1}).$$
Again use (\ref{D9}) and get
$$ -\frac{s_nv_n}{p_nu_n}=-\frac{s_n}{t_n}.$$

 Calculations for odd $n$ are analogous with slightly different initial description of matrix elements $A_{i,j}$ of the matrix $\An{n}(-\t)$:
\begin{align*}\notag & A_{i,j}=(-1)^{i+j}s_{i+1}s_{i+2}\cdots s_{j-1}u_j,&& 1\leq i\leq j\leq n,\\
&\notag A_{\i\j}=t_{j+1}t_{j+2}\cdots t_{i-1}u_i,&& 1\leq j\leq i\leq n,\\
& A_{i,\j}=(-1)^i(s_ns_{n-1}\cdots s_{i+1})\cdot(t_nt_{n-1}\cdots t_{j+1})&& 1\leq i<n,\ 1\leq j\leq n,\\
&\notag A_{n,\n}=0,\ \  A_{n,\j}=(-1)^i(t_nt_{n-1}\cdots t_{j+1}),&& \ 1\leq j< n,\\
&\notag A_{\i,j}=0,&& \ 1\leq i,j\leq n
\end{align*}
\hfill{$\square$}
\setcounter{equation}{0}
\section{$\GL(3)$ example}
Theorem \ref{t1} for $\GL(3)$ reads as follows:
 \begin{align} \notag
    \Psi_{\bm{\lambda}}&(\x)= \frac{e^{-i(\lambda_1x_1+\lambda_2x_2+\lambda_3x_3)}}{(2\pi i)^3}\int\limits_{\Re\bm{\ga}=0+}
    \exp\big({\ga}_{1,3}(x_{2}-x_{1})+(\gamma_{1,2}+\gamma_{2,3})(x_3-x_2)\big)\cdot
    \\ & \label{C1}
    \Gamma(\gamma_{1,2}+\gamma_{2,3}+i(\lambda_1-\lambda_3))\Gamma(\gamma_{2,3}+i(\lambda_2-\lambda_3))\Gamma(\gamma_{1,3}-\gamma_{2,3}+i(\lambda_1-\lambda_2))\\ \notag
   & \Gamma(\gamma_{1,2})\Gamma(\gamma_{1,3})\Gamma(\gamma_{2,3})d\gamma_{1,2}d\gamma_{1,3}d\gamma_{2,3}
    \end{align}
Using the notation
    $$\Lambda_1=\frac{2\lambda_1-\lambda_2-\lambda_3}{3},\qquad \Lambda_2=\frac{2\lambda_2-\lambda_1-\lambda_3}{3},\qquad
    \Lambda_3=\frac{2\lambda_3-\lambda_1-\lambda_2}{3}.$$
we rewrite \rf{C1}   as
\begin{align*} \notag
\Psi_{\bm{\lambda}}&(\x)= \frac{e^{-i\sum_{k,l=1}^3\frac{\lambda_k x_l}{3}}}{(2\pi i)^3}\int\limits_{\Re\bm{\ga}=0+}
\exp\big(({\ga}_{1,3}+i\Lambda_1)(x_{2}-x_{1})+(\gamma_{1,2}+\gamma_{2,3}-i\Lambda_3)(x_3-x_2)\big)\cdot
\\& \notag
\Gamma(\gamma_{1,2}+\gamma_{2,3}+i(\lambda_1-\lambda_3))\Gamma(\gamma_{2,3}+i(\lambda_2-\lambda_3))\Gamma(\gamma_{1,3}-\gamma_{2,3}+i(\lambda_1-\lambda_2))\\ 
& \Gamma(\gamma_{1,2})\Gamma(\gamma_{1,3})\Gamma(\gamma_{2,3})d\gamma_{1,2}d\gamma_{1,3}d\gamma_{2,3}
\end{align*}
Denote $q_1=\exp(x_1-x_2)$, $q_2=\exp(x_2-x_3)$, $s_1=\gamma_{1,3}+i\Lambda_1$, $s_2=\gamma_{1,2}+\gamma_{2,3}-i\Lambda_3$. Then
\rf{c2} can be presented as the inverse Mellin transform
$$ \Psi_{\bm{\lambda}}(\x)=\frac{e^{-i\sum_{k,l=1}^3\frac{\lambda_k x_l}{3}}}{(2\pi i)^2}\int_{\Re s_i=0+}q_1^{-s_1}q_2^{-s_2}
M(s_1,s_2)ds_1ds_2,$$
where
\begin{equation*}
\begin{split}
M(s_1,s_2)=
& \Gamma(s_1-i\Lambda_1)\Gamma(s_2+i\Lambda_1))\cdot\\
 \frac{1}{2\pi i}&\int_{\Re \gamma_{2,3}=0} \Gamma(\gamma_{2,3}+i(\lambda_2-\lambda_3))\Gamma(\gamma_{2,3})\Gamma(s_1-\gamma_{2,3}-i\Lambda_2)
 \Gamma(s_2-\gamma_{2,3}+i\Lambda_3)d\gamma_{2,3}.
 \end{split}
 \end{equation*}
Applying first Barnes lemma, we arrive to Bump \cite{Bu} formula
\begin{equation} \label{Bump} M(s_1,s_2)=\frac{\prod_{j=1}^3\Gamma(s_1-i\Lambda_j)\Gamma(s_2+i\Lambda_j)}{\Gamma(s_1+s_2)}
 \end{equation}
Note that Bump formula \rf{Bump} can be also derived from 'Gelfand-Tsetlin' presentation of $\GL(n,\R)$ Whittaker
function studied in \cite{GKL}. The derivation uses an integral calculated by de Branges and Wilson \cite{Br,W}
 \section*{Acknowledgements}
 The authors thank  V.Spiridonov and A.Shapiro for interesting discussions and A.Mironov for the help with the Maple package.
 The research  of the first author was supported  by RFBR grant 18-01-00460 used to obtain the results presented in sections 5,6,7.  The second author appreciates the support of RSF grant, project 16-11-10316 used  to obtain the results presented in sections 2,3,4.

\bibliographystyle{12}
\bibliographystyle{amsalpha}
		
	\end{document}